\documentclass[a4paper,leqno,11pt,twoside]{amsart}

\usepackage{amsmath,amsthm,amssymb,amsfonts,bm}
\usepackage[all]{xy}
\UseTips

\swapnumbers
\theoremstyle{definition} \newtheorem{defn}{Definition}[section]
                          \newtheorem{rmk}[defn]{Remark}
                          \newtheorem{properties}[defn]{Properties}
                          \newtheorem{wittdefn}[defn]{Definition of the Witt ring}
                          \newtheorem{fvt}[defn]{Maps: Frobenius, Verschiebung, Teichm\"uller}
                          \newtheorem{adw}[defn]{Alternative description of the Witt vectors}
                          \newtheorem{constr}[defn]{Construction of a Witt complex}
                          \newtheorem{ex}[defn]{Example}
                          \newtheorem{defn-prop}[defn]{Definition-Proposition}
                          \newtheorem{smox}[defn]{Some morphisms on $X_n$}
                          \newtheorem{pbcd}[defn]{Pull-back of Cartier Divisors}
                          \newtheorem{fpb}[defn]{Flat Pull-back of Cycles}
                          \newtheorem{pf}[defn]{Push-forward of Cycles}
                          \newtheorem{ep}[defn]{Exterior Product}
                          \newtheorem{icd}[defn]{Intersecting with Cartier Divisors}
                          \newtheorem{mkt}[defn]{Milnor $K$-Theory}

\theoremstyle{plain}      \newtheorem{thm}[defn]{Theorem}
                          \newtheorem*{thmwn}{Theorem}
                          \newtheorem{lem}[defn]{Lemma}
                          \newtheorem{cor}[defn]{Corollary}
                          \newtheorem{prop}[defn]{Proposition}

\newcommand{\eq}[2]{\begin{equation}\label{#1}#2 \end{equation}}
\newcommand{\ml}[2]{\begin{multline}\label{#1}#2 \end{multline}}

\newcommand{\surj}{\twoheadrightarrow}
\newcommand{\inj}{\hookrightarrow}
\newcommand{\red}{{\rm red}}

\newcommand{\sgn}{{\rm sgn}}

\newcommand{\Div}{{\rm div}}
\newcommand{\Hom}{{\rm Hom}}
\newcommand{\Spec}{{\rm Spec \,}}
\newcommand{\im}{{\rm im}}

\newcommand{\Tr}{{\rm Tr}}
\newcommand{\Nm}{{\rm Nm }}
\newcommand{\Res}{{\rm Res}}

\newcommand{\ord}{{\rm ord}}
\newcommand{\va}{{\rm v}}

\newcommand{\id}{{\rm id}}

\newcommand{\BW}[2]{\mathbb{W}_{#1}({#2})}
\newcommand{\BWC}[3]{\mathbb{W}_{#1}\Omega^{#2}_{#3}}
\newcommand{\W}[2]{\mathrm{W}_{#1}({#2})}
\newcommand{\WC}[3]{\mathrm{W}_{#1}\Omega^{#2}_{#3}}

\newcommand{\Th}[3]{\mbox{$\mathrm{TH}^{#1}(#2,#1;#3)$}}
\newcommand{\Ch}{\mathrm{CH}}
\newcommand{\gh}{\mathrm{gh}}
\newcommand{\Gh}{\mathrm{Gh}}
\newcommand{\V}{\mathrm{V}}
\newcommand{\F}{\mathrm{F}}
\newcommand{\R}{\mathrm{R}}

\newcommand{\bb}{\mathbb}

\newcommand{\Zy}{\mathrm{Z}}

\newcommand{\kaydot}{{\bm{\cdot}}}
\newcommand{\st}{\bm{*}}

\newcommand{\sD}{{\mathcal D}}

\newcommand{\sF}{{\mathcal F}}

\newcommand{\sI}{{\mathcal I}}

\newcommand{\sO}{{\mathcal O}}
\newcommand{\sP}{{\mathcal P}}

\newcommand{\sV}{{\mathcal V}}

\newcommand{\A}{{\mathbb A}}

\newcommand{\G}{{\mathbb G}}

\newcommand{\N}{{\mathbb N}}
\renewcommand{\P}{{\mathbb P}}
\newcommand{\Q}{{\mathbb Q}}

\newcommand{\Z}{\mathbb{Z}}

\begin{document}

\title[Additive Chow Groups and the De Rham-Witt Complex]{Additive Chow Groups with higher Modulus\\ and \\ the generalized De Rham-Witt Complex}
\author{Kay R\"ulling}  
\address{Universit\"at Duisburg-Essen, Essen, FB6, Mathematik, 45117 Essen, Germany}
\email{kay.ruelling@uni-essen.de}
\date{March 10, 2005}
\begin{abstract}
Bloch and Esnault defined additive higher Chow groups with modulus $(m+1)$ on the level of zero cycles over a field $k$, denoted by
$\Th{n}{k}{m}$, $n,m\ge 1$. They prove $\Th{n}{k}{1}\cong\Omega^{n-1}_{k/\Z}$. In this paper we generalize their result and obtain an isomorphism
\[\Th{n}{k}{m}\cong\BWC{m}{n-1}{k},\]
where $\BWC{-}{\kaydot}{k}$ is the generalized de Rham-Witt complex of Hesselholt-Madsen, generalizing the $p$-typical de Rham-Witt complex
of Bloch-Deligne-Illusie. 

Before we can prove this theorem we have to generalize some classical results to the de Rham-Witt complex. We give a 
construction of the generalized de Rham-Witt complex for $\Z_{(p)}$-algebras analogous to the construction in the 
$p$-typical case. We construct a trace 
$\Tr :\BWC{-}{\kaydot}{L}\to\BWC{-}{\kaydot}{k}$ for finite field extensions $L\supset k$ and if $K$ is the function field of a smooth projective curve
$C$ over $k$ and $P\in C$ is a point we define a residue map $\Res_P: \BWC{-}{1}{K}\to\BW{-}{k}$, which satisfies $\sum_{P\in C}\Res_P(\omega)=0$, for
all $\omega\in \BWC{-}{1}{K}$. 

\end{abstract}
\maketitle

\tableofcontents

\pagenumbering{arabic}

%
%
%
%
%
%
%

\section{Introduction}
Let $X$ be an equidimensional scheme of finite type over a field $k$ and write $\Delta_n=\Spec k[t_0,\ldots,t_n]/(\sum_{i=0}^nt_i-1)$. 
In \cite{Bl86} Bloch introduced higher Chow groups $\Ch^p(X,n)$, generalizing the Chow groups of $X$ (i.e. $\Ch^p(X,0)=\Ch^p(X)$).
Roughly speaking, these are defined by considering the quotient of $p$-codimensional cycles in $X\times\Delta_n$ 
in suitable good position modulo the boundary of $p$-codimensional cycles in $X\times\Delta_{n+1}$, where the boundary is given by intersecting
with the faces $(t_i=0)$ and then take the alternating sum. This construction is usually referred to as the simplicial definition of the higher Chow groups.
There is also a cubical one, which mainly differs by taking $(\P^1\setminus\{1\})^n$ instead of $\Delta_n$ (see \cite{To92}) and one can show that
these two definitions of Bloch's higher Chow groups coincide. One of the rare cases, where the higher Chow groups  can be computed, which means
to give a presentation in terms of generators and relations, is the case $X=\Spec k$ and $p=n$, i.e. formal sums of points in $\Delta_n$ 
(resp. $(\P^1\setminus\{0,1,\infty\})^n$) modulo the boundary of formal sums of curves in good position in $\Delta_{n+1}$ 
(resp. $(\P^1\setminus\{1\})^{n+1}$). In this case Nesterenko-Suslin and Totaro proved (see \cite{NeSu89}, \cite{To92})
\[\Ch^n(k,n)\cong K^M_n(k),\]
where $K^M_n(k)$ are the degree $n$ elements in the Milnor ring of $k$.

We observe that one could replace $\Delta_n$ in the definition of the higher Chow groups by $\Spec k[t_0,\ldots, t_n]/(\sum t_i-\lambda)$, 
for any $\lambda\in k^\times$. In \cite{BlEs03a} Bloch and Esnault investigated the degenerated case $\lambda=0$. They obtain a theory
of additive higher Chow groups, $\text{SH}^p(X,n)$, and prove in particular
\[\text{SH}^n(k,n)\cong\Omega^{n-1}_{k/\Z},\]
using a presentation of $\Omega^{n-1}_{k/\Z}$ as a quotient of the anti commutative graded ring $k\otimes_\Z\bigwedge^* k^\times$ modulo
the graded ideal generated by $a\otimes a+ (1-a)\otimes(1-a)$. 

In \cite[Section 6.]{BlEs03} Bloch and Esnault construct a cubical version of the additive higher Chow groups (so far only for fields and on the level
of zero cycles). Since these groups are our main object of study, we give a more precise definition. 
Denote by $\Zy_d(X)$ the group of $d$-dimensional cycles on $X$ and write 
\[X_n=\G_m\times(\P^1\setminus\{1\})^n,\quad \text{ with  coordinates } (x,y_1,\ldots, y_n).\]
Define $\Zy_1(X_n;1)$ to be the subgroup of $\Zy_1(X_n),$ which is freely generated by 
1-dimensional subvarieties $C\subset X_n$, $C\not\subset \bigcup_{i}(y_i=0,\infty)$ satisfying the following properties
         \begin{enumerate}
          \item[(a)] ({\em Good position}) $(y_i=j).[C]\in\Zy_0(X_{n-1}\setminus \bigcup_{i}(y_i=0,\infty))$, for $i=1,\ldots,n$, $j=0,\infty$. 
          \item[(b)] ({\em Modulus $2$ condition}) If $\nu: \widetilde{C}\to \P^1\times \bigl(\P^1\bigr)^n$ is the normalization of the 
                     compactification of $C$, then
              \eq{1}{2[\nu^*(x=0)]\le\sum_{i=1}^n[\nu^*(y_i=1)]\quad\text{in } \Zy_0(\widetilde{C}).} 
         \end{enumerate}

There is a map $\partial=\sum_{i=1}^n (-1)^i(\partial^0_i-\partial^\infty_i):\Zy_1(X_n;1)\to\Zy_0(X_{n-1}\setminus Y_{n-1})$ 
and Bloch and Esnault define
\[\text{TH}^n(k,n)=\Th{n}{k}{1}=\frac{\Zy_0(X_{n-1}\setminus Y_{n-1})}{\partial\Zy_1(X_n;1)}.\]  
They show by similar methods as before, assuming $1/6\in k$,
\[\Th{n}{k}{1}\cong \Omega^{n-1}_{k/\Z}.\]
Now the natural question is, what happens, if we replace the modulus $2$ condition in (b), by a modulus $(m+1)$ condition, 
i.e. replace the $2$ in (\ref{1})
by $(m+1)$ and $\Th{n}{k}{1}$ by $\Th{n}{k}{m}$? The answer to this question, which is given in this thesis, is motivated by the following considerations. 
Denote by $\text{Pic}(\A^1_k, (m+1)\{0\})$ the equivalence classes of divisors in $\A^1$, supported on $\G_m$. Two divisors $D$, $D'$ being equivalent
iff $D-D'=\Div(f/h)$, for some $f,h\in 1+tk[t]$ with $f-h\in t^{m+1}k[t]$, this is $f/h\equiv 1 $ mod $(m+1)\{0\}$ in the language of
\cite[Chapter III, \S 1]{Se88}. Now given such functions $f$ and $h$ we can define a curve $C\subset \G_m\times\P^1\setminus\{1\}$ by
the equation $h(x)y-f(x)=0$ and one easily checks that this curve satisfies the modulus $(m+1)$-condition. 
Thus we obtain a well defined and surjective map
\eq{2}{\text{Pic}(\A^1,(m+1)\{0\})\longrightarrow \Th{1}{k}{m}.}
Now it is quiet reasonable to believe this map to be an isomorphism (at least for $m=1$, both sides equal $k$, by the above). But we may identify
$\text{Pic}(\A^1,(m+1)\{0\})$ with the group $\left(\frac{1+tk[t]}{1+t^{m+1}k[t]}\right)^\times$, which in turn may be identified with the ring 
of generalized Witt vectors of length $m$ over $k$, $\BW{m}{k}$, (see \cite{Bl78}). Hence we hope to obtain an isomorphism
\[\Th{1}{k}{m}\cong \BW{m}{k}.\]
If this is true, then the groups $\Th{n}{k}{m}$ should in general give something, which generalizes the absolute K\"ahler differentials on the one hand
and the big Witt rings of finite length on the other. The natural suspect for this is the generalized de Rham-Witt complex of Hesselholt-Madsen,
$\BWC{m}{\kaydot}{k}$, which generalizes the $p$-typical de Rham-Witt complex of Bloch-Deligne-Illusie. And indeed the main theorem of this thesis is
\begin{thmwn}
Let $k$ be a field of characteristic $\neq 2$. Then we have for all $n,m\ge 1$ an isomorphism
\[\Th{n}{k}{m}\cong\BWC{m}{n-1}{k}.\]
\end{thmwn}
The first part of the proof, namely to define the map from the additive Chow groups to the de Rham-Witt complex, is analogous to the proof of the case 
$m=1$ by Bloch and Esnault. To verify that the map $\Th{n}{k}{m}\to \BWC{m}{n-1}{k}$ is well defined, we use a reciprocity law 
(as it was done in \cite{NeSu89}, \cite{To92}, \cite{BlEs03a} and \cite{BlEs03}). 
Here it is a "sum-of-residues-equal-zero" theorem, which we prove in section 3, generalizing
the well known residue formula for differentials on smooth projective curves. To obtain the inverse map we use the universality of the de Rham-Witt
complex, i.e. we equip the additive Chow groups with a structure of such a complex and then the universality of the de Rham-Witt complex yields a map 
$\BWC{m}{n-1}{k}\to \Th{n}{k}{m}$, which is inverse to the map we constructed first.

Now we want to give a more detailed description of the structure of this thesis.

In section 1 we recall the definition and the main properties of the (big) Witt vectors over a ring for arbitrary truncation sets $S\subset \N$, 
following \cite{Be66}. (Recall, $S\subset \N$ is  a truncation set iff for each $n\in S$ all the divisors of $n$ are contained in $S$.) 
We give an alternative description via formal power series (see \cite{Bl78}) and construct a trace 
$\Tr_{A/B}: \BW{S}{B}\to\BW{S}{A}$ 
for $A$-algebras $B$, with $B$ a free $A$-module of finite rank, via the norm on the formal power series (following hints by Bloch and Faltings).

In section 2 we give the definition of a generalized Witt complex over a ring $A$ (following Hesselholt-Madsen) as a contravariant functor from
the category of differential graded algebras, equipped with Frobenius maps $\F_n$ and Verschiebung maps $\V_n$ ($n\in\N$) satisfying certain relations. 
The de Rham-Witt complex over $A$ is
by definition the initial object in the category  of Witt complexes over $A$. Using category theory Hesselholt-Madsen showed in \cite{HeMa01} that 
an initial object always exists. We give a more concrete construction for $A$ a $\Z_{(p)}$-algebra, $p\neq 2$ a prime. Following a method of Illusie
(see \cite{Il79} or \cite{HeMa04}) we first construct a de Rham-Witt complex without Frobenius as a quotient of the absolute K\"ahler differentials
over the Witt vectors and then we use the results from the $p$-typical case (\cite{HeMa04}) to obtain Frobenius maps. 
In the rest of the section we give some standard
properties of the de Rham-Witt complex, which easily follow from the corresponding results in the $p$-typical situation.  

The first aim of section 3 is to define a trace $\Tr_{L/k}:\BWC{S}{n}{L}\to \BWC{S}{n}{k}$ for arbitrary finite field extensions $L\supset k$ and finite
truncation sets $S$, which generalizes the trace on the K\"ahler differentials (see \cite{Ku64}, \cite{Ku86}). This is done by defining it separately
for separable field extensions and purely inseparable ones of degree $p$ and a combination of these two cases in general. If $L\supset k$ is separable,
then we can write $\BWC{S}{n}{L}=\BW{S}{L}\otimes_{\BW{S}{k}}\BWC{S}{n}{k}$ and the trace is defined as $\Tr\otimes\id$, with $\Tr$ being the trace
constructed in section 1. If $L\supset k$ is purely inseparable of degree $p$, the trace of a Witt differential $\omega\in\BWC{S}{n}{L}$ is defined by 
the following procedure:
first lift $\omega$ to the level $S\cup pS$, then multiply the lift $\tilde{\omega}$ with $p$, then there is an element 
$\alpha$ in $\BWC{S\cup pS}{n}{k}$, which maps to $p\tilde{\omega}$ and then the trace is defined by $\Tr(\omega)=\alpha_{|S}$.
Afterwards we define, for a point $P$ on a smooth projective curve $C$, a residue symbol for Witt 1-forms over the function field of $C$, 
generalizing the residue symbol on the K\"ahler differentials (and also the one defined by Witt in \cite{Wi36}, see also \cite{AnRo04}). 
We prove the "sum-of-residues-equal-zero" theorem, using the classical method (use a trace formula to reduce the proof to $\P^1$ and then
calculate explicitly). Finally we generalize the definition of the residue symbol and the residue theorem to higher Witt forms.

In section 4 we finally arrive at our object of study, namely the additive cubical Chow groups. We give the definition, as it was done above
and show that we have a well defined map 
\[\theta: \Th{n}{k}{m}\to\BWC{m}{n-1}{k},\quad [P]\mapsto \Tr_{k(P)/k}(\psi_{n-1}(P)),\]
 with
\[\psi_{n-1}=\frac{1}{[x]}\frac{d[y_1]}{[y_1]}\ldots\frac{d[y_{n-1}]}{[y_{n-1}]}\in \BWC{m}{n-1}{k(x,\underline{y})}.\]
We prove that the map is well defined by generalizing the arguments in \cite{BlEs03}: We first show that, if $C$ satisfies the modulus
condition, then $\psi_n(\widetilde{C})$ has no poles in $(x=0)$ (where $\widetilde{C}$ is the normalization of the compactification of $C$). It follows
that 
\[\sum_P\Res_P(\psi_n(\widetilde{C}))=\theta(\partial[C]).\]
Hence $\theta$ is well defined by the residue theorem of section 3.

Afterwards we equip the additive Chow groups with the structure of a (restricted) Witt complex, where the multiplication is induced by taking exterior 
products of cycles and then pushing forward via the multiplication map on $\G_m$. The differential comes from pushing forward with the diagonal embedding
$\G_m\to\G_m\times\G_m$ and the Frobenius $\F_r$ (resp. the Verschiebung $\V_r$) is induced by pushing forward (resp. pulling back) via the map 
$\G_m\to\G_m$, $x\mapsto x^r$. We use the results of Nesterenko-Suslin and Totaro to show that these maps satisfy the relations, which should be 
satisfied in a Witt complex and then the universality of the de Rham-Witt complex gives a map, inverse to $\theta$. This finishes the proof.

In the appendix, finally, we list some results from \cite{Fu84}, concerning push-forward and pull-back of cycles as well as intersecting cycles with 
Cartier divisors.  \\
\\
{\bf Acknowledgments}
I am very grateful to H\'el\`ene Esnault for giving me the subject of this  thesis and for her excellent guidance and support, not only during the
work on this problem, but from the very beginning of my mathematical studies. 

I want to express my gratitude to Lars Hesselholt for explaining the Witt vectors and the generalized de Rham-Witt complex to me.

I would like to thank Spencer Bloch and Gerd Faltings for several hints, in particular, to construct the trace for the Witt vectors via the norm map.

I want to thank Stefan Kukulies for the many discussions concerning the foundations of algebraic geometry, Andre Chatzistamatiou for the 
careful reading of the introduction and Wioletta Syzdek for the supply with chocolate.

Finally I want to thank my parents for making untroubled studies possible.

This paper is my PhD thesis written at the university of Essen financially supported by the Graduiertenkolleg 
{\em Mathematische und ingenieurwissenschaftliche Methoden f\"ur sichere Daten\"ubertragung und Informationsvermittlung}.

\numberwithin{equation}{defn}
\newpage

%
%
%
%
%
\section{Witt Vectors}
In this section we will summarize the definition and properties of the Witt vectors. See \cite{Wi36} for 
the original paper, \cite{Be66} for a detailed treatment, \cite{Bl78} for the description via power series 
and \cite{HeMa01} for the language we are using.\\
\\
Rings are always assumed to be commutative and unital. The natural numbers do not contain $0$, i.e. $\N=\{1,2,3,\ldots\}$.\\
\\
A subset of the natural numbers $S\subset \N$ is  called a \emph{truncation set} if for any number $n\in S$
all its divisors are in $S$ too. Examples of truncation sets are $\N$ itself, $\{1,2\ldots,m\}$, for any $m\in \N$, 
$\sP=\{1, p, p^2,\ldots \}$, for a prime $p$, and
\[<n>=\{d | d \text{ divides } n\}.\]
If $S$ is a truncation set and $n\in \N$ we get a truncation set
\[S/n=\{d\in \N|nd\in S\}\subset S.\]
We have $(S/n)/m=S/(nm)$. Notice, that $S/n$ is not empty iff $n\in S$. 
We denote by $J$ the \emph{category of truncation sets} with inclusions as morphisms. It is a filtered category. For $S\in J$ we write
$J_S$ for the restricted category whose objects are subsets of $S$. 

\begin{wittdefn}\label{0.1}
Let $A$ be a ring and $S\in J$ a truncation set. The \emph{(big) Witt ring} is by definition the ring $\BW{S}{A}$ whose
elements are of the form $(w_s)_{s\in S}$, for $w_s\in A$, and whose ring structure is uniquely determined by the 
condition that the \emph{ghost map}
\[\gh:\BW{S}{A}\longrightarrow A^S,\quad \gh_s(w)=\sum_{d|s}dw_d^{s/d}\]
is a natural transformation of functors of rings. (The proof of the existence of such a ring structure reduces to
$A=\Z[x_1,x_2,\ldots]$ and follows then from the lemma below.) 
An element of $\BW{S}{A}$ is called a \emph{ Witt vector}. We write
\[\BW{}{A}:=\BW{\N}{A},\quad \BW{m}{A}=\BW{\{1,\ldots,m\}}{A}\]
and if a prime $p$ is fixed 
\[\W{}{A}:=\BW{\{1,p,p^2,\ldots\}}{A},\quad \W{n}{A}:=\BW{\{1,p,\ldots, p^{n-1}\}}{A},\]
which we call the \emph{$p$-typical Witt ring} and the \emph{$p$-typical Witt ring of length $n$}, respectively.    
\end{wittdefn}

The following lemma is taken from lecture notes from Michael Hopkins (written by Matthew Ando); knowing the formulation, the proof
is an easy exercise.
\begin{lem}\label{0.2}
Write $A=\Z[x_1,x_2,\ldots]$ and, for any prime $p$, let $\varphi_p: A\longrightarrow A$ be the map defined by
$\varphi_p(\sum a_Ix^I)=\sum a_I x^{pI}$. Let $S$ be a truncation set and $c=(c_s)_{s\in S}\in A^S$. Then $c$ is in the 
image of the ghost map if and only if for all $s\in S$ and all primes $p$,
\[c_s\equiv \varphi_p(c_{s/p}) \mod p^{v_p(s)}A.\]
\end{lem}

It is time for a first list of Properties, which follow immediately from the definition.
\pagebreak[2]

\begin{properties}\label{0.3}
Let $A$ be a ring and $S\in J$. 
\begin{enumerate}
\item If $\varphi :A\longrightarrow B$ is a homomorphism of rings, then
\[\BW{S}{A}\longrightarrow \BW{S}{B},\quad (w_s)_{s\in S}\mapsto (\varphi(w_s))_{s\in S}\]
is a homomorphism of rings, which will be also denoted by $\varphi$. If $\varphi$ is surjective (respectively injective) on 
the rings, then the same holds on the level of Witt vectors.
\item $\BW{\{1\}}{A}=A$.
\item If $A$ has no $\Z$-torsion, then the ghost map $\gh : \BW{S}{A}\longrightarrow A^S$ is injective, for all truncation 
        sets $S$.
\item If $A$ is a $\Q$-algebra, then the ghost map is an isomorphism, for all truncation sets $S$.
\item If $T\subset S$ are truncation sets, then there is a natural restriction map
\[\R^S_T:\BW{S}{A}\longrightarrow\BW{T}{A},\quad (w_s)_{s\in S}\mapsto (w_s)_{s\in T}.\]
It is a surjective homomorphism of rings. If $w\in\BW{S}{A}$, we also write $\R^S_T(w)=w_{|T}$.  
\end{enumerate}
\end{properties}

\begin{fvt}\label{0.4}
Let A be a ring and $S\in J$. For each $n\in\N$ there is a natural map
\[\F_n:\BW{S}{A}\longrightarrow \BW{S/n}{A}\]
which is indirectly defined by
\[\gh_s(\F_n(w))=\gh_{sn}(w).\] 
The existence follows from lemma \ref{0.2}. We call $\F_n$ the \emph{($n$-th) Frobenius}. 

The \emph{($n$-th) Verschiebung} 
\[\V_n:\BW{S/n}{A}\longrightarrow \BW{S}{A}\]
is given by
\[(\V_n(w))_s= \begin{cases}
                   w_{s/n} & \text{if } n|s\\
                      0    & \text{else}
               \end{cases}.\]  

And finally, the \emph{Teichm\"uller map} 
\[[-]=[-]_S: A\longrightarrow \BW{S}{A}\]
is defined by
\[[a]=(a,0,0,\ldots).\] 
\end{fvt}

\begin{properties}\label{0.5}
Let $A$ be a ring, $S\in J$ and $n\in \N$. Then
\begin{enumerate}
\item \[\gh_s(\V_n(w))=\begin{cases}
                         n\gh_{s/n}(w) & \text{if } n|s\\
                                0      & \text{else}
                        \end{cases},\quad s\in S\]
        for all $w\in\BW{S/n}{A}$.
\item$\gh_s([a]_S)=a^s$,  for all $s\in S$ and $a\in A$.
\item The Teichm\"uller map is multiplicative. Furthermore, $[1]_S$ is the $1$ in $\BW{S}{A}$ and $[0]_S$ is the 0.
\item The Frobenius $\F_n:\BW{S}{A}\longrightarrow\BW{S/n}{A}$ is a unital ring homomorphism. We have
       $\F_1=\id$ and $\F_m\circ\F_n=\F_{mn}$, for all $m$, $n$. If $A\to B$ is a ring homomorphism, then the Frobenius maps
       on $\BW{S}{A}$ and on $\BW{S}{B}$ form the obvious commutative diagram.
\item Write ${\F_n}_*\BW{S/n}{A}$ for $\BW{S/n}{A}$ considered as $\BW{S}{A}$-module via $\F_n$. Then
       $\V_n:{\F_n}_*\BW{S/n}{A}\longrightarrow \BW{S}{A}$ is $\BW{S}{A}$-linear. We have
       $\V_1=\id$ and $\V_m\circ\V_n=\V_{mn}$, for all $m$, $n$. If $A\to B$ is a ring homomorphism, then the 
       Verschiebung maps on $\BW{S}{A}$ and on $\BW{S}{B}$ form the obvious commutative diagram.
\end{enumerate}
\end{properties}
\begin{proof} (i) and (ii) are easy. The other statements reduce to $A=\Z[x_1,x_2,\ldots]$ and follow then by the 
              injectivity of the ghost map from (i), (ii) and the definition, respectively.
\end{proof}

\begin{rmk}\label{0.6}
If $p$ is a fixed prime number, it is customary to write $\V :\W{n}{A}\longrightarrow\W{n+1}{A}$ instead of 
$V_p$ and then, of course, by (v) above, $\V^n$ instead of $\V_{p^n}$. Similarly, we write 
$\F:\W{n}{A}\longrightarrow \W{n-1}{A}$ instead of $\F_p$ and $\F^n$ for $\F_{p^n}$. 
\end{rmk}

We give a final list of formulas.

\begin{properties}\label{0.7}
Let $A$ be a ring and $S\in J$. Then
\begin{enumerate}
\item For $w=(w_s)_{s\in S}\in \BW{S}{A}$ we have
       \[w=\sum_{n\in S}\V_n([w_n]_{S/n}).\]
\item If $(m,n)=1$, then $F_m\circ\V_n=\V_n\circ\F_m$. 
\item $\F_n\circ\V_n=n$.
\item $\F_m(\V_n([a]))=(m,n)\V_{\frac{n}{(m,n)}}([a]^{m/(m,n)})$, for all $a\in A$.
\item $\V_n([a])\V_r([b])=(n,r)\V_{\frac{nr}{(n,r)}}([a]^{r/(n,r)}[b]^{n/(n,r)})$, for all $a,b\in A$.
\item $[a]\V_n(w)=\V_n([a]^nw)$, for all $a\in A$ and $w\in \BW{S}{A}$.
\end{enumerate}
\begin{proof}
Same as in \ref{0.5}.
\end{proof}
\end{properties}

\begin{rmk}\label{0.8}
Notice, that we can do all this also if $A$ is not unital. Especially it makes sense to talk about 
the Witt ring of an ideal.
\end{rmk}

\begin{lem}[\cite{He04}, Lemma 1.2.1.]\label{0.8.1}
Let $S$ be a truncation set, $A$ a ring and $I\subset A$ an ideal. Then $\BW{S}{I}\subset\BW{S}{A}$ is an ideal
and we have an isomorphism
\[\BW{S}{A}/\BW{S}{I}\stackrel{\simeq}{\longrightarrow}\BW{S}{A/I}.\] 
\end{lem}
\begin{proof}
The proof is the one of Hesselholt (there for $p$-typical Witt vectors). It is enough to prove the isomorphism for finite truncation sets.
For $S=\{1\}$ it is trivial. Now let $S$ be any finite truncation set with at least two elements
and assume the assertion is proved for all $S_0\varsubsetneq S$. Take $n\in S$ maximal, then $S/n=\{1\}$ and we get the following
diagram with exact columns\\
\[\xymatrix{     &      0\ar[d]     &     0\ar[d]      &     0\ar[d]     &       \\
           0\ar[r] & I\ar[r]\ar[d]^{\V_n} &  A\ar[r]\ar[d]^{\V_n} &  {A/I}\ar[r]\ar[d]^{\V_n} & 0\\
           0\ar[r] & \BW{S}{I}\ar[r]\ar[d]^\R & \BW{S}{A}\ar[r]\ar[d]^\R & \BW{S}{A/I}\ar[r]\ar[d]^\R   & 0\\
           0\ar[r] & \BW{S\setminus \{n\}}{I}\ar[r]\ar[d] & \BW{S\setminus \{n\}}{A}\ar[r]\ar[d] & \BW{S\setminus \{n\}}{A/I}\ar[r]\ar[d] & 0\\
               &      0            &    0             &    0           & }\]
Now the top and the bottom row are exact, hence the middle one is exact too.
\end{proof}

\begin{adw}\label{0.10}
There is another description of the Witt vectors, which is in some sense more intuitive. Let $A$ be a ring. We write
\[\Gamma(A)=(1+TA[[T]])^\times.\]
Any $P\in\Gamma(A)$ can be uniquely written as a product
\[P(T)=\prod_{n\ge 1}(1-w_nT^n),\quad w_n\in A.\]
Therefore, we get a bijection
\[w:\Gamma(A)\longrightarrow \BW{}{A},\quad \prod(1-w_nT^n)\mapsto (w_n)_{n\in\N}\]
Now we define $\gamma:\Gamma(A)\longrightarrow A^\N$ to be the composition of the maps
\[\Gamma(A)\stackrel{-T\frac{d}{d T}\log}{\longrightarrow}TA[[T]]\stackrel{\simeq}{\longrightarrow} A^\N,\]
where the latter map is given by $\sum a_n T^n\mapsto (a_n)_{n\in \N}$. The map $\gamma$ is a group homomorphism
and functorial in $A$. One checks
\[\gh(w(P))=\gamma(P).\] 
It follows from the definition of the Witt ring that $w:\Gamma(A)\longrightarrow \BW{}{A}$ is an isomorphism
of groups. For $S\in J$, the kernel of the composition with the restriction $R^\N_S\circ w$ is given by
\[I_S=\{\prod_{n\not\in S}(1-w_nT^n)| w_n\in A\},\]
in particular $I_{\{1,2,\ldots,m\}}=(1+T^{m+1}A[[T]])^\times$. Denoting by $\Gamma_S(A)$ the quotient $\Gamma(A)/I_S$, we
get an isomorphism of groups, which we will call $w$ again
\[w:\Gamma_S(A)\stackrel{\simeq}{\longrightarrow}\BW{S}{A}.\]
The ring structure of $\BW{S}{A}$ induces one on $\Gamma_S(A)$ and we thus obtain an alternative description of the Witt
ring. The  Verschiebung and Teichm\"uller maps act as follows
\begin{itemize}
\item $[a]=w(1-aT)$
\item $V_n([a])=w(1-aT^n)$
\end{itemize} 
and the properties \ref{0.7}, (iv) and (v) determine the multiplication and the Frobenius maps.

{From} now on we will not distinguish between these two descriptions (so forget about the $\Gamma$-notation).
\end{adw}

\begin{rmk}\label{0.11}
\begin{enumerate}
\item If $p$ is a prime, it is now easy to see,  
\[\V_p\circ\F_p=p \quad\text{iff}\quad pA=0\]
and an integer $n\in \Z$ is a unit in $\BW{S}{A}$ iff it is a unit in $A$.
\item The Frobenius may be defined in a more natural way. Indeed, it comes out that 
$\F_n:\BW{}{A}\longrightarrow\BW{}{A}$ is the map induced by the norm map $\Nm: A[[T^{1/n}]]\longrightarrow A[[T]]$.
\item Notice that the above description of the Witt vectors is not restricted to  rings with $1$. Indeed, if $A$ is not unital, we may
      define the group $(1+TA[[T]])^\times$, with $1$ a formal symbol and group operation
      \[(1+TP(T))\cdot(1+TQ(T))=1+(P(T)+Q(T))T+P(T)Q(T)T^2.\]
      One easily checks that this gives a group and that the above discussion still works in this case.
\end{enumerate}
\end{rmk}

Next we want to define a trace map. The treatment via the norm map on the ring of power series is due to hints by Bloch and Faltings.

\begin{prop}\label{0.12}
Let $A$ be a ring and $B$ an $A$-algebra, which is a free $A$-module of finite rank. Let $S$ be a truncation set. Then 
the norm map $\Nm: B[[T]]\longrightarrow A[[T]]$ induces a trace map
\[\Tr=\Tr_{B/A,S}: \BW{S}{B}\longrightarrow\BW{S}{A},\]
(usually we will just write $\Tr$, but the notations $\Tr_{B/A}$ or $\Tr_S$ may also occur, depending on the point we want to stress)
which satisfies the following properties
\begin{enumerate}
\item $\gh\circ\Tr=\Tr\circ \gh$, where the latter trace is the usual one $\Tr: B^S\longrightarrow A^S$.
\item $\Tr$ is $\BW{S}{A}$-linear, it commutes with restriction and we have
          \[  \Tr\circ\F_n=\F_n\circ \Tr \quad \text{and}\quad \Tr\circ \V_n=\V_n\circ\Tr.\]
\item We have
       \[\Tr\biggl(\sum_{n\in S}\V_n([a_n])\biggr)=\sum_{n\in S}\Tr\left(\V_n([a_n])\right).\]
\item If $C$ is a $B$-algebra, which is a free $B$-module of finite rank, then
       \[\Tr_{C/A}=\Tr_{B/A}\circ\Tr_{C/B}.\]
\item If $A\longrightarrow A'$ is a ring homomorphism and $B'=B\otimes_A A'$, then we have the following commutative diagram\\
       \[\xymatrix{ \BW{S}{B}\ar[d]_\Tr\ar[r] & \BW{S}{B'}\ar[d]^\Tr\\
                    \BW{S}{A}\ar[r]             &    \BW{S}{A'}.             }\]
\item If we have $A$-algebras $B_i$, $i=1,\ldots,r$, which are free $A$-modules of finite rank and $B=\prod_{i=1}^rB_i$, then 
      $\BW{S}{B}=\prod_{i=1}^r\BW{S}{B_i}$ and  
            \[\Tr_{B/A}(w)=\sum_{i=1}^r\Tr_{B_i/A}(w_i),\]
     for $w=(w_1,\ldots,w_r)\in \BW{S}{B}$.

\end{enumerate}
\end{prop}
\begin{proof} Since $B$ is a free $A$-module of finite rank, $B[[T]]$ is a free $A[[T]]$-module of finite rank. Thus we have a norm map. Let 
 $e=\{e_1,\ldots,e_d\}$ be a basis of $B$ over $A$ and $\epsilon_i\in A$ such that $1=\sum_{i=1}^d\epsilon_i e_i$. 
Now take a $1$-power series over $B$, $f=1+T^mb(T)\in (1+T^mB[[T]])^\times$ and write $b(T)=\sum_{i=1}^d a_i(T)e_i$. If we denote by $E_i$ the 
matrix of $e_i$ in the basis $e$, the matrix of $f$ is given by 
\[F=\sum(\epsilon_i+T^ma_i(T))E_i=\id+T^m\sum a_i(T)E_i.\]
Calculating the determinant via the Leibniz rule, we obtain 
\[\Nm(f)\in(1+T^mA[[T]])^\times.\]
Thus, for all $m\in\N$, $\Nm$ induces a well defined group homomorphism (denoted by $\Nm$ again)
\[\Nm: (1+TB[[T]])^\times/(1+T^mB[[T]])^\times\longrightarrow (1+TA[[T]])^\times/(1+T^mA[[T]])^\times.\]
Now, since $(1+TB[[T]])^\times=\varprojlim \left((1+TB[[T]])^\times/(1+T^mB[[T]])^\times\right)$, we see that the map 
\[\Nm:(1+TB[[T]])^\times\longrightarrow (1+TA[[T]])^\times\]
satisfies 
\begin{equation}\label{0.12.1} \Nm\biggl(\prod_{i\ge 1}(1-a_iT^i)\biggr)=\prod_{i\ge 1}\Nm(1-a_iT^i).\end{equation}
(Which was not a priori clear, since the product is infinite.)
Write $I_S(B)=\{\prod_{j\not\in S}(1-b_jT^j)| b_j\in B\}$ and similarly for $A$. We want to show $\Nm(I_S(B))\subset I_S(A)$.
By (\ref{0.12.1}) it is enough to show 
\[\Nm(1-bT^n)\in I_S(A),\quad \text{for } b\in B \text{ and } n\not\in S.\]
For this we write
\[\Nm(1-bT^n)=\prod(1-\alpha_jT^j),\quad \alpha_j\in A\]
and the claim follows, if we show $\alpha_j=0$, for all $j$ with $n\not|j$. Assume this is not the case, 
i.e. $r:=\inf\{j\,|\,\, n\not|j \text{ and } \alpha_j\neq 0 \}$ is a natural number.
Obviously we have 
\[\Nm(1-bT^n)=1+\sum_{i\ge 1}a_iT^{in},\quad a_i\in A.\]
This shows, that, if we expand $\prod(1-\alpha_j T^j)$ in a power series, the coefficient of $T^r$ has to be zero, that is
\[\biggl(\sum_{l=2}^r\sum_{\substack{1\le j_1<\ldots<j_l \\ j_1+\ldots+j_l=r}}(-1)^l\alpha_{j_1}\ldots\alpha_{j_l}\biggr)-\alpha_r=0.\]
But by the definition of $r$ the term in brackets is zero, thus $\alpha_r=0$, a contradiction.

All in all the norm map induces, by \ref{0.10}, a map 
\[\Tr:\BW{S}{B}\longrightarrow \BW{S}{A},\]
which is additive and satisfies (iii), (iv), (v) and (iv).
Next we want to show (i).  For this we may assume $S=\N$. Recall from \ref{0.10} that the ghost map
is given by the composition
\[\gh: \BW{}{B}=(1+TB[[T]])^\times\stackrel{-T\frac{d }{d T}\log}{\longrightarrow}TB[[T]]\stackrel{\simeq}{\longrightarrow} B^\N\]
and that, under the latter isomorphism, the trace on $B^\N$ corresponds to the trace on $TB[[T]]$. Thus it remains to show
\[\frac{\Nm(P)'}{\Nm(P)}=\Tr\left(\frac{P'}{P}\right),\quad P\in (B[[T]])^\times,\]
where $'$ means the derivation relative to $T$. Here we follow an argument of \cite{Hu89}. Let $M=(a_{ij})$ be the matrix of $P$ in the Basis $e$
and $\tilde{M}=(\tilde{a}_{ij})$, with $\tilde{a}_{ij}=(-1)^{i+j}\det M_{ji}$, $M_{ji}$ arising from $M$ by canceling the $j$-th row and 
the $i$-th column. Then $M'$ is the matrix of $P'$ and
\[(\det M)'=\sum_{i,j}\tilde{a}_{ji}a_{ij}'.\]
(Since, if we view the $a_{ij}$'s as indeterminants, the Laplace expansion shows 
 $\partial M/\partial a_{ij}=(-1)^{j+i}\det M_{ij}$ and thus $d( \det M)=\sum_{i,j} (-1)^{j+i}\det M_{i,j}d a_{ij}$.)
So we get
\begin{eqnarray*}
\Tr\left(\frac{P'}{P}\right) & = & \Tr(M^{-1}M')=\frac{1}{\det M}\Tr(\tilde{M}M')= \frac{1}{\det M}\sum_{i,j}\tilde{a}_{ij}a_{ji}'\\
                             & = & \frac{(\det M)'}{\det M}=\frac{\Nm(P)'}{\Nm(P)}.
\end{eqnarray*}
The last thing we have to show is (ii), i.e. the linearity and that the trace commutes with Frobenius and Verschiebung. 
Denote $\tilde{A}=\Z[y_a|a\in A]$. We get a map $\tilde{A}\longrightarrow A,\quad y_a\mapsto a$. Now define $\tilde{B}:=\tilde{A}^n$.
There is a ring structure on  $\tilde{B}$ such that the natural map $\tilde{B}\longrightarrow B=\oplus_{i=1}^n A e_i$
is a surjective ring homomorphism and we have $\tilde{B}\otimes_{\tilde{A}}A=B$. 
(Take $a_{ijk}\in A$ such that $e_ie_j=\sum a_{ijk}e_k$, for $1\le i\le j$ and write $f_i$ for the standard basis
 vector in $\tilde{B}$ , i.e. $f_1=(1,0\ldots,0)$, etc. . Then define the ring structure on $\tilde{B}$ via, 
$f_if_j=f_jf_i=\sum y_{a_{ijk}}f_k$. Thus the map $\tilde{B}\longrightarrow B,\quad y_af_i\mapsto ae_i$ becomes a ring homomorphism.
Notice, that $\tilde{B}$ is in general not unital, but it does not need to be by remark \ref{0.11}, (iii).) 
Thus we get a commutative diagram as in (v). So we may assume that
$A$ is $\Z$-torsion free and since $\gh$ is now injective, the linearity follows from (i) and the linearity of the usual trace map.
The statement for the Verschiebung follows from property \ref{0.5}, (i) and  for the Frobenius from \ref{0.4}. 
\end{proof}

\begin{ex}\label{0.13}
Let $L\supset k$ be a finite separable field extension of degree n. Denote by $\sigma_1,\ldots,\sigma_n$ the $n$ distinguished embeddings from
$L$ into $\bar{k}$, the algebraic closure of $k$. Then the $\sigma_i$'s define maps from the Witt ring over $L$ in the Witt ring over $\bar{k}$
(again denoted by $\sigma_i$) and we get
\[\Tr w=\sum_i\sigma_i(w), \text{ for all } w\in \BW{S}{L}.\]
(This follows immediately from the corresponding relation for the norm.)

If $L\supset k$ is a purely inseparable field extension of degree $p^e$ and $a\in L$, we have
\[\Tr\V_n([a])=\V_n(\V_{p^e}\F_{p^e}([a]))=\V_{np^e}([a^{p^e}]),\]
by remark \ref{0.11}, (i) this is just multiplication with $p^e$. (This is because $[L((T)):k((T))]=p^e$ and thus
$\Nm(1-aT)=(1-aT)^{p^e}=(1-a^{p^e}T^{p^e})$.)  
\end{ex}

Let $p$ be a prime. We denote by $\Z_{(p)}$ the localization of $\Z$ with respect to the prime ideal $(p)\subset\Z$.
If $A$ is a $\Z_{(p)}$-algebra, the Witt ring may be described in terms of the $p$-typical one:

\begin{prop}\label{0.9}
Let $p$ be a prime, $A$ a $\Z_{(p)}$-algebra, $S\in J$, $\sP=\{1,p,p^2,\ldots\}$ and 
$I_p=\{n\in\N |(n,p)=1\}$. For $n\in I_p$ we define
\[\epsilon_n=
\prod_{\substack{q \text{ prime }\in S \\q\neq p}}\left(\tfrac{1}{n}\V_n(1)-\tfrac{1}{nq}\V_{nq}(1)\right)\in \BW{S}{A}.\]
Then  there is an isomorphism of rings 
\[\varphi: \BW{S}{A}\stackrel{\simeq}{\longrightarrow}\prod_{n\in I_p}\BW{\sP\cap S/n}{A},\quad 
                         w\mapsto (\F_1(w)_{|\sP\cap S},\ldots,F_n(w)_{|\sP\cap S/n},\ldots),\]
with inverse map 
\[\psi: \prod_{n\in I_p}\BW{\sP\cap S/n}{A}\longrightarrow \BW{S}{A},\quad
     (w_1,\ldots,w_n,\ldots)\mapsto \sum_{n\in \text{I}_p}\tfrac{1}{n}\epsilon_n\V_n(\tilde{w}_n),\]
where $\tilde{w}_n\in \BW{S/n}{A}$ is any lifting of $w_n$ under the restriction map. 
\end{prop}
\begin{proof}
The $\epsilon_n$'s enjoy the following properties
\begin{enumerate}
\item[a)] $\gh_s(\epsilon_n)=\begin{cases}
                              1 & \text{if } s\in n\sP\cap S\\
                              0 & \text{else}
                             \end{cases}$, for all $s\in S$.
\item[b)]For $m,n\in I_p$ and $m\neq n$, $\F_m(\epsilon_n)=0$.
\item[c)] $\F_n(\epsilon_n)_{|\sP\cap S/n}=1\in\BW{\sP\cap S/n}{A}$.
\end{enumerate}
Indeed, a) follows from the property \ref{0.5}, (i). For b) and c) it is enough to consider $A=\Z_{(p)}[x_1,\ldots]$.
Then b) follows from a) and $\gh_s\circ \F_m=\gh_{sm}$. For a prime $q$ with $q\neq p$, we remark that
$\V_q(w)_{|\sP\cap S/n}=0$, $w\in \BW{S}{A}$. By property \ref{0.7}, (iii)
\[\F_n(\epsilon_n)_{|\sP\cap S/n}=\prod_{p\neq q}\left(1-\tfrac{1}{q}\V_q(1)\right)_{|\sP\cap S/n}=1
                                                                                           \in \BW{\sP\cap S/n}{A},\]
thus c).

Now $\varphi$ is obviously a unital ring homomorphism, so it remains to check that $\psi$ is the inverse map. 
$\varphi\circ\psi=\id$ follows immediately form b), c) and \ref{0.7}, iii).
For $w\in \BW{S}{A}$, we want to show, $\psi(\varphi(w))=w$. We may assume $A=\Z_{(p)}[x_1,\ldots,]$. By a) we have, for
any $s\in S$, 
\[\gh_s\left(\frac{1}{n}\epsilon_n\V_n(\widetilde{\F_n(w)_{|\sP\cap S/n}})\right)=\begin{cases}
                                                                                \gh_s(w) & \text{if } s\in n\sP\cap S\\
                                                                                   0      & \text{else}
                                                                               \end{cases}\]
and, since $S=\bigsqcup_{n\in I_p\cap S}n\sP\cap S$, we see that $\gh_s(\psi(\varphi(w))=\gh_s(w)$, for all $s\in S$. 
\end{proof}
\newpage

%
%
%
%
%
%

\section{The generalized de Rham-Witt Complex}
We give the definition of Hesselholt's and Madsen's \cite{HeMa01} generalization of the de Rham-Witt complex of Bloch \cite{Bl78}, Deligne and Illusie
\cite{Il79}. It  is the initial object in a category whose objects are functors from the category of truncation sets into the category of differential
graded algebras, which are equipped with Frobenius and Verschiebung maps, satisfying certain relations. Hesselholt and Madsen show the existence 
in \cite{HeMa01} using category theory. For $\Z_{(p)}$-algebras, $p$ odd, we construct the generalized de Rham-Witt complex as a quotient of the absolute
de Rham complex over the Witt vectors, following the method of Illusie. Therefore, we first construct a Witt complex without Frobenius, which works 
completely analogous to the p-typical case and then we show that there is already a Frobenius on it, 
using the results of Illusie or the more general results of Hesselholt-Madsen in the $p$-typical case. Afterwards we list some properties, which all
follow easily from corresponding  results in the $p$-typical situation. \\
\\
We denote by DGA the category of differential graded $\Z$-algebras.

\begin{defn}\label{1.1}
Let $A$ be a ring. A \emph{V-complex over A} is  a contravariant functor
\[E:J\longrightarrow \text{DGA},\]
which transforms direct limits into inverse ones and posses the following additional structure:
We denote the elements of degree $i$ in the dga $E(S)$ by $E(S)^i$, then there is a natural transformation of rings
\[\lambda:\BW{S}{A}\longrightarrow E(S)^0,\quad \text{for all } S\in J\]
and for all $n\in \N$,  natural transformations of graded $\Z$-modules
\[\V_n:E(S/n)\longrightarrow E(S), \quad \text{for all } S\in J,\]
satisfying
\begin{enumerate}
\item $\V_1=\id,\quad \V_n\circ\V_m=\V_{nm},\quad \text{and }\V_n\circ\lambda=\lambda\circ\V_n$, where the latter $\V_n$ is 
      the Verschiebung on $\BW{S}{A}$.
\item $\V_n(xd y)=\V_n(x)d\V_n(y),\quad$ for all $x\in E(S/n)^i$ and $y\in E(S/n)^j$, all $i,j$.
\item $\V_n(x)d\lambda([a])=\V_n(x\lambda([a])^{n-1})d\lambda(\V_n([a])),\quad$ for all $a\in A$ and $x\in E(S/n)^i$, all $i$. 
\end{enumerate}
A morphism of V-complexes, $f:E\longrightarrow E'$, is a natural transformation of differential graded algebras in degree 0, which is 
compatible with the $\V_n$'s and $\lambda$'s. 

If $p$ is a fixed prime, $\sP=\{1,p,p^2,\ldots\}$ and we take the restricted category $J_\sP$ instead of $J$, we call such a functor 
\emph{p-typical $\V$-complex}.
\end{defn}

\begin{rmk}\label{1.2}
In \cite{Il79}, chapter I, Illusie defined a ($p$-typical) V-pro-complex over $A$ to be the following thing (see also \cite{HeMa04}, section 5):
A projective system of differential graded $\Z$-algebras, $E=((E_n)_{n\in \N}, R: E_n\to E_{n-1})$, together with a map of projective systems of rings,
$\lambda:\W{\kaydot}{A}\longrightarrow E^0_{\kaydot}$ and a map of projective systems of groups $\V: E_{\kaydot-1}\longrightarrow E_\kaydot$, satisfying
$\lambda\circ\V=\V\circ \lambda$, $\V(xdy)=\V(x)d\V(y)$ and $\V(x)d\lambda([a])=\V(x\lambda([a])^{p-1})d\V(\lambda([a]))$, 
for all $x,y\in E$ and $a\in A$. A $p$-typical-V-complex, in our sense, is obviously one in the sense of Illusie. On the other hand, if $E$ is a 
V-pro-complex in the sense of Illusie and $E_n^\kaydot$ is as a dga generated by $E_n^0$ and $\lambda$ is surjective, then $E$ is also a 
$p$-typical-V-complex in our sense. One only has to show $\V^s(x)d\lambda([a])=V^s(x\lambda([a])^{p^s-1})d\V^s(\lambda[a])$,
for all $s\ge 1$, $a\in A$ and $x\in E$.
But from the property \ref{0.7}, (vi), together with the extra assumption, one sees easily $\lambda([a])\V(x)=\V(\lambda([a]^{p})x)$ and 
the statement follows by induction. 
\end{rmk}

The following proposition is a straightforward generalization of \cite{Il79}, chapter I, theorem 1.3 (see also \cite{HeMa04}, proposition 5.1.1).

\begin{prop}\label{1.3}
Let $A$ be a ring. Then there is an initial object in the category of V-complexes,
\[S\mapsto \bb{W}_S\Omega^\kaydot_A,\]
which we call the \emph{de Rham-Witt complex}. If $S$ is a finite truncation set, then there is an epimorphism of dga's
\[\Omega_{\BW{S}{A}/\Z}^\kaydot\longrightarrow\bb{W}_S\Omega^\kaydot_A.\]
Furthermore,
\[\BWC{\{1\}}{\kaydot}{A}=\Omega_{A/\Z}^\kaydot,\quad \text{and} \quad \bb{W}_S\Omega^0_A=\BW{S}{A}.\] 
\end{prop}
\begin{proof}
We begin with the construction for finite truncation sets. Set
\[\bb{W}_{\{1\}}\Omega^\kaydot_A=\Omega^\kaydot_{A/\Z}.\]
Now let $S$ be a finite truncation set. 
Assume we have constructed a functor from the category $\{S_0\varsubsetneq S\}\subset J_S$ to DGA, $S_0\mapsto \bb{W}_{S_0}\Omega^\kaydot_A$, 
together with a natural transformations of  graded groups,
$\V_n: \bb{W}_{S_0/n}\Omega_A^\kaydot\longrightarrow \bb{W}_{S_0}\Omega_A^\kaydot$, for all $n$, satisfying the following properties, for all 
$S_0\varsubsetneq S$
\begin{enumerate}
\item[$(a)_{S_0}$] $\bb{W}_{S_0}\Omega^0_A=\BW{S_0}{A}$ and the $\V_n$'s and the restriction maps on either side coincide.
\item[$(b)_{S_0}$] $\V_1=id, \quad V_n\circ\V_m=\V_{nm}$.
\item[$(c)_{S_0}$] $\V_n(xdy)=\V_n(x)d\V_n(y), \quad$ for all $x\in \bb{W}_{S_0/n}\Omega^i_A$ and $y\in \bb{W}_{S_0/n}\Omega^j_A$.
\item[$(d)_{S_0}$] $\V_n(x)d[a]=\V_n(x[a]^{n-1})d\V_n([a]),\quad$ for all $x\in \bb{W}_{S_0/n}\Omega^i_A$ and $a\in A$.
\item[$(e)_{S_0}$] We have an epimorphism of dga's $\pi:\Omega^\kaydot_{\BW{S_0}{A}/\Z}\longrightarrow \bb{W}_{S_0}\Omega^\kaydot_A$.
\end{enumerate}
Now we define $N_S^\kaydot\subset\Omega^\kaydot_{\BW{S}{A}}$ to be the differential graded ideal, generated by the following elements, for all $n>1$
\begin{enumerate}
\item[$(I)_n$] \[\sum_j\V_n(x_j)d\V_n(y_{1j})\ldots d\V_n(y_{ij}),\]
               for all $x_j$, $y_{lj}\in \bb{W}_{S/n}\Omega^0_A=\BW{S/n}{A}$ with $\sum_j x_jdy_{1j}\ldots dy_{ij}=0$ in $\bb{W}_{S/n}\Omega^i_A$.
\item[$(II)_n$] \[\V_n(x)d[a]-\V_n(x[a]^{n-1})d\V_n([a]),\]
                    for all $a\in A$ and $x\in\bb{W}_{S/n}\Omega^0_A$.
\end{enumerate} 
We define
\[\bb{W}_S\Omega^\kaydot_A=\frac{\Omega^\kaydot_{\BW{S}{A}/\Z}}{N^\kaydot_S}.\]
Clearly, $\bb{W}_S\Omega^0_A=\BW{S}{A}$. Let $S_0\subset S$ be a truncation set. We have a restriction map $\BW{S}{A}\longrightarrow \BW{S_0}{A}$ and
composing the induced map on the absolute de Rham complex with the map from $(e)_{S_0}$, we get a map of dga's
$\Omega^\kaydot_{\BW{S}{A}}\longrightarrow \bb{W}_{S_0}\Omega^\kaydot_A$. Using $(c)_{S_0}$ and $(d)_{S_0}$, we see that this map factors to give a
map of dga's
\[R^S_{S_0}:\bb{W}_S\Omega^\kaydot_A\longrightarrow \bb{W}_{S_0}\Omega^\kaydot_A.\]
We obtain a functor
\[J_S\longrightarrow \text{DGA},\quad S_0\mapsto \bb{W}_{S_0}\Omega^\kaydot_A.\]   
For each $n\in \N$ we define a map of graded groups $\V_n:\bb{W}_{S/n}\Omega^\kaydot_A\longrightarrow \bb{W}_S\Omega^\kaydot_A$ by
$\V_n(xdy_1\ldots dy_i):=\V_n(x)d\V_n(y_1)\ldots d\V_n(y_i)$, for all $x,y_i\in \BW{S/n}{A}$, where the $\V_n$ on the right hand side is the Verschiebung
on the Witt vectors. By definition the $\V_n$'s are well defined and fulfill (i), (ii) and (iii) from the definition \ref{1.1}.

Now take any $S\in J$. We write $J^f_S$ for the category whose objects are finite truncation sets contained in $S$. Then we define in DGA
\[\bb{W}_S\Omega^\kaydot_A=\varprojlim_{S_0\in J^f_S}\bb{W}_{S_0}\Omega^\kaydot_A.\]
For $S'\subset S$, the restriction maps on the finite truncation sets defined above induce a map of dga's
\[\bb{W}_{S}\Omega^\kaydot_A\longrightarrow \bb{W}_{S'}\Omega^\kaydot_A.\]
We get a functor
\[J\longrightarrow \text{DGA},\quad S\mapsto\bb{W}_S\Omega^\kaydot_A,\]
which transforms direct limits into inverse ones and coincides in degree $0$ with $S\mapsto \BW{S}{A}$. Notice,
$J^f_{S/n}=\{S_0/n | S_0\in J_S^f\}$, thus we can write
\[\bb{W}_{S/n}\Omega^\kaydot_A=\varprojlim_{S_0\in J^f_S}\bb{W}_{S_0/n}\Omega^\kaydot_A\]
and we see that the $\V_n$'s, which we defined for finite truncation sets, induce maps of graded groups
\[\V_n:\bb{W}_{S/n}\Omega^\kaydot_A\longrightarrow \bb{W}_S\Omega^\kaydot_A.\]
Thus we get a V-complex $S\mapsto \bb{W}_S\Omega^\kaydot_A$. 

It remains to show that this is the initial object in the category of V-complexes. So take a V-complex $E$ and a finite truncation set $S$.
Then the map $\lambda:\BW{S}{A}\longrightarrow E(S)^0$ induces a homomorphism of dga's, which factors by definition \ref{1.1}, (i), (ii) and (iii) 
to give a map
\[\bb{W}_S\Omega^\kaydot_A\longrightarrow E(S)^\kaydot,\quad xdy_1\ldots dy_i\mapsto \lambda(x)d\lambda(y_1)\ldots d\lambda(y_i),
                                                                                                                    \quad x,y_i\in \BW{S}{A},\]
which commutes with restriction and Verschiebung. Furthermore this arrow is unique with the property that it is $\lambda$ in degree $0$ and since
$E$ transforms direct limits into inverse ones, we get a unique homomorphism of V-complexes $\bb{W}_{-}\Omega^\kaydot_A\longrightarrow E(-)^\kaydot$
and we are done. 
\end{proof}

\begin{rmk}\label{1.4}
We write $\BWC{}{\kaydot}{A}:=\BWC{\N}{\kaydot}{A}$. If $p$ is a fixed prime and $\sP=\{1,p,\ldots\}$, then we will say 
\emph{$p$-typical de Rham-Witt complex} instead of ``de Rham-Witt complex restricted to $J_\sP$''. 
By remark \ref{1.2} this coincides with the de Rham-Witt pro-complex of Deligne-Illusie.
 We write $\WC{n}{\kaydot}{A}$ for the $p$-typical de Rham-Witt complex evaluated in $\{1,\ldots,p^{n-1}\}$ 
and $\WC{}{\kaydot}{A}$, if we evaluate it in $\sP$ and $\V:=\V_p$. 
\end{rmk}

The following definition is  \cite{HeMa01}, definition 1.1.1., it generalizes the definition of a Witt pro-complex (with Frobenius) of Deligne-Illusie.

\begin{defn}\label{1.5}
Let $A$ be a ring.  A \emph{Witt complex over $A$} is a contravariant functor
\[E:J\longrightarrow \text{DGA},\] 
which transforms direct limits into inverse ones, together with natural transformations of graded rings
\[\F_n: E(S)\longrightarrow E(S/n),\quad \text{for all }n\in \N\]
and  natural transformations of graded groups
\[\V_n: E(S/n)\longrightarrow E(S),\quad \text{for all }n\in \N,\]
satisfying the following relations, for all $n,m \in\N$
\begin{enumerate}
\item $\F_1=\V_1=\id,\quad \F_m\F_n=\F_{mn},\quad \V_m\V_n=\V_{nm}.$
\item $\F_n\V_n=n$, and if $(m,n)=1$, then  $\F_m\V_n=\V_n\F_m$. 
\item $ \V_n(\F_n(x)y)=x\V_n(y)$, for all $x\in E(S)$, $y\in E(S/n)$ and all $n\in\N$.
\item  $\F_m d\V_n = k d\F_{m/c}\V_{n/c}+l \F_{m/c}\V_{n/c}d$,\\
       where $c=(m,n)$ and $l,k\in \Z$ are arbitrary with  $km+ln=(m,n)$.
\end{enumerate}
Furthermore, there is a natural transformation of rings
\[\lambda: \BW{S}{A}\longrightarrow E(S)^0,\]
which  commutes with $\F_n$ and $\V_n$ and satisfies
\begin{enumerate}
\item[(v)] $\F_n d\lambda([a])=\lambda([a]^{n-1})d\lambda([a])$, for all $a\in A$ and $n\in \N$.
\end{enumerate}

A morphism of Witt complexes, is a natural transformation of differential graded algebras, compatible with $\F_n$, $\V_n$ and $\lambda$.

If $p$ is a fixed prime, we say \emph{$p$-typical Witt complex} for a Witt complex on $J_\sP$ and write $\F:=\F_p$ and
$\V_p:=\V$.
\end{defn}

In \cite{HeMa01}, proposition 1.1.5., Hesselholt and Madsen proved, using the Freyd adjoint functor theorem that there  always exists an initial object 
in the category of Witt complexes. For $\Z_{(p)}$-algebras, we will show, using the corresponding results from the $p$-typical situation that
the de Rham-Witt complex is the initial object in the category of Witt complexes. 

\begin{properties}\label{1.7}
Let $E$ be a Witt complex over a ring $A$. Then 
the following equalities hold, for all $n\in \N$, $S\in J$, $x,y\in E(S/n)$ and all $a\in A$
\begin{enumerate}
\item[(a)] $\F_nd\V_n=d$.
\item[(b)] $d\F_n=n\F_nd,\quad \V_nd=nd\V_n$.
\item[(c)] $\V_n(xdy)=\V_n(x)d\V_n(y),\quad \V_n(x)d\lambda([a])=\V_n(x\lambda([a]^{n-1}))d\lambda\V_n([a]).$
\item[(d)] $F_md\V_n(\lambda([a]))=kd\V_{n/c}(\lambda([a])^{m/c})+
                                                             l\V_{n/c}(\lambda([a])^{\frac{m}{c}-1})d\V_{n/c}(\lambda([a]))$,\\
           where $c=(m,n)$ and $l,k\in \Z$ are arbitrary with  $km+ln=(m,n)$.
\end{enumerate}
\end{properties}
\begin{proof}
(a) follows immediately from (iv) of the definition. The first equation of (c) follows from (a) and (iii) of the definition, the second equation
 from the first, (iii) and (v). Now (b),
\begin{eqnarray*} \V_n(dx) & = & \V_n(1)d\V_n(x)= \V_n(1)d\V_n(x)+\V_n(x)d\V_n(1) = d(\V_n(x)\V_n(1))\\
                           & = & d\V_n(\F_n\V_n(x))=nd\V_n(x)
\end{eqnarray*}
and similarly
\begin{eqnarray*}
d\F_n(x) & = & \F_nd\V_n\F_n(x) = \F_nd(\V_n(1)x)= \F_n(d(\V_n(1))x)+\F_n(\V_n(1)dx)\\
         & = & \F_nd(V_n(1))\F_n(x)+ \F_n\V_n(1)\F_n(dx)=n\F_n(dx).
\end{eqnarray*}
It remains to show (d). First we consider $(m,n)=1$ and take $l,k\in\Z$ with $km+ln=1$. Then by (iv) and (ii)
\[\F_md\V_n([a])=kd\F_m\V_n([a])+l\F_m\V_nd [a]= kd\V_n([a^m])+l\V_n([a]^{m-1}d[a])\]
and by (c) we are done in this case. The general statement follows now from $\F_md\V_n=\F_{m/c}d\V_{n/c}$.
\end{proof}

\begin{rmk}\label{1.8}
Property (c) shows that we have a forgetful functor from the category of Witt complexes into the category of V-complexes.
\end{rmk}

The following theorem was first proved for $\bb{F}_p$-algebras by Illusie, see \cite{Il79}, chapter I, theorem 2.17, and then generalized
by Hesselholt-Madsen for $\Z_{(p)}$-algebras.

\begin{thm}[\cite{HeMa04}, Theorem D]\label{1.9}
Let $p$ be an odd prime and $A$ a $\Z_{(p)}$-algebra. Then there is a map of projective systems of graded rings
\[F: \WC{\kaydot}{\kaydot}{A}\longrightarrow \WC{\kaydot-1}{\kaydot}{A},\]
which is in degree zero the Frobenius of the Witt ring over $A$ and satisfies
\[\F d\V=d,\quad \F\V=p,\quad \V(\F(x)y)=x\V(y), \text{ for all } x\in\WC{n}{\kaydot}{A},y\in\WC{n-1}{\kaydot}{A}\]
and
\[\F d[a]=[a]^{p-1}d[a], \text{ for all }a\in A.\] 
\end{thm}

\begin{cor}\label{1.10}
$\WC{\kaydot}{\kaydot}{A}$ is the initial object in the category of $p$-typical Witt complexes.
\end{cor}
\begin{proof}
First we have to show that $\WC{\kaydot}{\kaydot}{A}$ is a $p$-typical Witt complex, i.e. we have to show that 
(iv) and (v) of definition \ref{1.5} hold true in our $p$-typical situation. (v) follows by induction and the fact that
$\F$ coincides with the Frobenius on the Witt vectors. For (iv) we must show
\[\F^s d\V^r=kd\F^{s-t}\V^{r-t}+l\F^{s-t}\V^{r-t}d,\]
with $t=\min(r,s)$ and $kp^s+lp^r=p^t$. But we have $\V d=pd\V$ and $d\F=p\F d$ 
(proved in the same way as in \ref{1.7}, b)). Thus if $t=s$, we have
\[kd\V^{r-s}+l\V^{r-s} d=d\V^{r-s}=\F^sd\V^r\]
and if $t=r$,
\[kd\F^{s-r}+l\F^{s-r} d=\F^{s-r}d=\F^s d\V^r.\]
So we see that $\F$ makes $\WC{\kaydot}{\kaydot}{A}$ a $p$-typical Witt complex. Now, if $E$ is any $p$-typical 
Witt complex over $A$, we have a unique map of $p$-typical V-complexes $\WC{\kaydot}{\kaydot}{A}\longrightarrow E$, 
which becomes automatically (by property \ref{1.7}, (d)) a map of Witt complexes. Thus $\WC{\kaydot}{\kaydot}{A}$
is the initial object in the category of $p$-typical Witt complexes.
\end{proof}

Now, for a $\Z_{(p)}$-algebra $A$, we want to decompose the de Rham-Witt complex into $p$-typical parts 
(analogous to the decomposition of the  Witt vectors into the $p$-typical ones, as in proposition \ref{0.9}).
Then the F from the theorem \ref{1.9} will define $\F_n$'s on the de Rham-Witt complex and we see that 
$\BWC{-}{\kaydot}{A}$ is the initial object in the category of Witt complexes. These calculations were already done
in \cite{HeMa01}, section 1.2., but since this construction helps us to reduce almost all statements concerning
the de Rham-Witt complex to the $p$-typical one, we will include it here.\\
\\
If $(\Omega,d)$ is a dga, we denote $\Omega(\frac{1}{j})=(\Omega,\frac{1}{j}d)$.

\begin{constr}[see \cite{HeMa01}, Section 1.2.]\label{1.11}
Let $p$ be an odd prime and $A$ a $\Z_{(p)}$-algebra. Then we define a contravariant functor
$M: J\longrightarrow \text{DGA}$, by 
\[M(S)^\kaydot=\prod_{j\in I_p}\BWC{\sP\cap S/j}{\kaydot}{A}(\tfrac{1}{j}),\]
with the product ring structure and the differential, $\partial$, given by
\[\partial(\alpha_j)_{j\in I_p}=(\tfrac{1}{j}d\alpha_j)_{j\in I_p}.\]
$M$ transforms direct limits into inverse ones and by proposition \ref{0.9}, we have an isomorphism in degree zero
\[\lambda: \BW{S}{A}\longrightarrow M(S)^0,\quad w\mapsto (\F_j(w)_{|\sP\cap S/j})_{j\in I_p},\]
in particular $\lambda([a])=([a]^j)_{j\in I_p}$, $a\in A$ . For $n\in\N$ we define maps
\[\V_n: M(S/n)\longrightarrow M(S),\quad \F_n:M(S)\longrightarrow M(S/n)\]
by
\[\V_n(\alpha)_j=\begin{cases}
                h\V^s(\alpha_{j/h}) & \text{if } h|j\\
                         0         & \text{else}
                \end{cases},\]
where $n=p^sh$ with $(h,p)=1$ and $\alpha=(\alpha_j)_{j\in I_p}\in M(S/n)$, and for $\beta=(\beta_j)_{j\in I_p}\in M(S)$
\[\F_n(\beta)_j=F^s(\beta_{jh}).\] 
Here $\V$ and $\F$ are the maps from proposition \ref{1.3} and theorem \ref{1.9}. (One easily checks that $\V_n$ and 
$\F_n$ map the elements, where they should be mapped to.) Obviously, $\F_n$ is a natural transformation of graded rings
and $\V_n$ is one of graded groups.

\begin{prop}\label{1.12}
The functor $M$ constructed above is a Witt complex.
\end{prop}
\begin{proof}
We have to check the relations of definition \ref{1.5}. This is all straightforward, so we will just present 
the most confusing part, i.e. that $\V_n$ and $\F_n$ commute with $\lambda$ and relation (iv).

First we show, $\V_n$ commutes with $\lambda$. Take $w\in\BW{S/n}{A}$, $j\in I_p\cap S$ and 
write $n=p^sh$, with $(h,p)=1$, then by definition
\[\V_n(\lambda(w))_j= \begin{cases}
                      h\V^s(\F_{j/h}(w)_{|\sP\cap S/(nj/h)}) & \text{if }h|j\\
                                        0                   & \text{else}
                     \end{cases}.\]
On the other hand, if we write $c=(j,n)$,
\[\lambda(\V_n(w))_j= \F_j(\V_n(w))_{|\sP\cap S/j}=c\F_{j/c}(\V_{n/c}(w))_{|\sP\cap S/j}= 
                                                                     c\V_{n/c}(\F_{j/c}(w))_{|\sP\cap S/j}.\]
Now if the last term is not zero, then $n/c\in \sP\cap S/j$, which is equivalent to $ c=h, n/c=p^s$. Thus we have
\[\lambda(\V_n(w))_j=\V_n(\lambda(w))_j.\]
Next, $\F_n$ commutes with $\lambda$. For $w\in \BW{S}{A}$ and $n$ as above,
\[\F_n(\lambda(w))_j=\F^s(\lambda(w)_{jh})=\F_{p^s}(\F_{jh}(w)_{|\sP\cap S/(jh)}).\]
But, this last term lives on $\sP\cap (S/jh)/p^s=\sP\cap S/(jn)$, so
\[\F_n(\lambda(w))_j=(\F_{p^s}\F_{jh}(w))_{|\sP\cap S/(jn)}=(\F_{j}\F_{n}(w))_{|\sP\cap S/(jn)}=\lambda(\F_n(w))_j.\] 
Finally relation (iv). After checking relation (i) (easy), this is equivalent to (by \ref{1.7} (a) and (b))
\[\F_n\partial\V_n=\partial,\, n\partial\V_n=\V_n\partial \text{ and } n\F_n\partial=\partial \F_n.\]
Now again this is easy, but we include the calculations as an example. Write $n=p^sh$ with $(h,p)=1$ and take
$\alpha\in M(S/n)$, $j\in I_p$. Then 
\[(\F_n\partial\V_n(\alpha))_j=\F^s(\partial\V_n(\alpha))_{jh}=\frac{1}{j}\F^sd\V^s(\alpha_j)=(\partial\alpha)_j\]
and if $h$ divides $j$
\[(n\partial\V_n(\alpha))_j= \frac{nh}{j}d\V^s(\alpha_{j/h})= \frac{h}{j}h\V^sd(\alpha_{j/h})=(\V_n\partial(\alpha))_j,\]
which obviously also holds if $h$ does not divide $j$.
For $\alpha\in M(S)$
\[(n\F_n\partial\alpha)_j=n\F^s(\frac{1}{jh}d\alpha_{jh})=\frac{1}{j}d\F^s(\alpha_{jh})=(\partial\F_n\alpha )_j\]
and we are done.
\end{proof}
\end{constr}

\begin{thm}\label{1.13}
Let $p$ be an odd prime and $A$ a $\Z_{(p)}$-algebra. Then the de Rham-Witt complex over $A$, $\BWC{-}{\kaydot}{A}$,
is the initial object in the category of Witt complexes. Furthermore, the isomorphism of proposition \ref{0.9} induces
an isomorphism of Witt complexes
\eq{1.13.1}{\varphi: \BWC{-}{\kaydot}{A}\stackrel{\simeq}{\longrightarrow} M(-),}
with
\[ \varphi(xdy_1\ldots dy_i)=
               \bigl(\F_j(x)\partial_j\F_j(y_1)\ldots\partial_j\F_j(y_i)_{|\sP\cap S/j}\bigr)_{j\in I_p},\quad x,y_l\in \BW{S}{A}.\]
\end{thm}
\begin{proof}
If we equip the de Rham-Witt complex with a structure of a Witt complex, then it is automatically the initial object.
Since, if we take a Witt complex $E$, we get by remark \ref{1.8} a unique map of V-complexes, 
which by property \ref{1.7}, (d), commutes with the $\F_n$'s. Thus it is enough to show that we have an isomorphism
of V-complexes $\varphi: \BWC{-}{\kaydot}{A}\stackrel{\simeq}{\longrightarrow} M(-)$, since then the $\F_n$'s on $M$
induce $\F_n$'s on the de Rham-Witt complex. As both functors take direct limits to inverse ones, we are reduced to
prove the isomorphism for finite truncation sets. So take a finite truncation set $S\in J$. By proposition \ref{0.9}
we have an isomorphism
\[\BW{S}{A}\longrightarrow \prod_{j\in I_p}\BW{\sP\cap S/j}{A},\quad w\mapsto (\F_j(w)_{|\sP\cap S/j})_{j\in I_p},\]
with inverse map
\[\psi: \prod_{j\in I_p}\BW{\sP\cap S/j}{A}\longrightarrow \BW{S}{A},\quad
     (w_j)_{j\in I_p}\mapsto \sum_{n\in \text{I}_p}\tfrac{1}{j}\epsilon_j\V_j(\tilde{w}_j),\]
where $\tilde{w}_j\in \BW{S/j}{A}$ is any lifting of $w_j$ under the restriction map and the $\epsilon_j$'s are
as in proposition \ref{0.9}. 
If we denote the differential of $\Omega_{\BW{\sP\cap S/j}{A}}^\kaydot(\frac{1}{j})$ by $\partial_j=\frac{1}{j}d$,
this map induces an isomorphism of dga's
\[\varphi : \Omega^\kaydot_{\BW{S}{A}}\stackrel{\simeq}{\longrightarrow} 
                                                   \prod_{j\in I_p}\Omega^\kaydot_{\BW{\sP\cap S/j}{A}}(\tfrac{1}{j}),\]
with
\[\varphi(xdy_1\ldots dy_i)=
               \bigl(\F_j(x)\partial_j\F_j(y_1)\ldots\partial_j\F_j(y_i)_{|\sP\cap S/j}\bigr)_{j\in I_p},\quad x,y_l\in \BW{S}{A}.\]
We observe that $\epsilon_j^2=\epsilon_j$ (this follows from (a) in the proof of proposition \ref{0.9}) and thus the
inverse map is given by
\eq{1.13.1.1}{\varphi^{-1}\bigl((x_j\partial_j y_{1j}\ldots\partial_j y_{ij})_{j\in I_p}\bigr)=
  \sum_{j\in I_p}\tfrac{1}{j^{i+1}}\epsilon_j\V_j(\tilde{x}_j)d\V_j(\tilde{y}_{1j})\ldots d\V_j(\tilde{y}_{ij}).}
Now we have to show that this gives a well defined isomorphism as in (\ref{1.13.1}). For $S=\{1\}$ nothing is to
prove. So we take any finite truncation set $S$ and assume that the statement is known for all $S_0\varsubsetneq S$.
Let $N^\kaydot_S$ be the differential graded ideal from the proof of proposition \ref{1.3}, generated by $(I)_n $ and $(II)_n$,
for $1<n\in S$. Then we have to show
\eq{1.13.1.5}{\varphi(N^i_S)= \prod_{j\in I_p} N^i_{\sP\cap S/j}.}
That the left hand side is contained in the right, follows if we prove 
\eq{1.13.2}{\pi\varphi(\V_n(x)d\V_n(y_1)\ldots d\V_n(y_i))=\V_n(\pi\varphi(x dy_1\ldots dy_i))\text{ in } M(S),}
where $\pi :\prod_{j\in I_p}\Omega^\kaydot_{\BW{S}{A}}(\tfrac{1}{j})\longrightarrow M(S)$ is the natural map, the
$\V_n$ on the right hand side is the one from construction \ref{1.11} and
$x,y_i\in\BW{S/n}{A}$. But if we write $\varphi=(\varphi_j)_j$, then
\[\pi\varphi_j(\V_n(x)d\V_n(y_1)\ldots d\V_n(y_i))= 
              \lambda(\V_n(x))_jd\lambda(\V_n(y_1))_j\ldots d\lambda(\V_n(y_i))_j\]
and thus (\ref{1.13.2}) follows from the fact that $M$ is a Witt complex. 
It remains to show that the right hand side of (\ref{1.13.1.5}) is contained in the left.
Denote  $e_j=(0,\ldots,0,1,0\ldots)$, with $1$ sitting in the $j$-th place. 
Now take $x_n, y_{ln}\in \BW{(\sP\cap S/j)/p^s}{A}$ with 
\[\sum_n x_n d y_{1n}\ldots dy_{in}=0\in\BWC{(\sP\cap S/j)/p^s}{i}{A}\]
 and write 
$\alpha=e_j\sum_n\V_{p^s}(x_n)\partial_j\V_{p^s}(y_{1n})\ldots\partial_j\V_{p^s}(y_{in})$. Then
\[\varphi^{-1}(\alpha)=
 \sum_n\V_{jp^s}(\tfrac{1}{j^{i+1}}\F_{jp^s}(\epsilon_j )\tilde{x}_n)
                                            d\V_{jp^s}(\tilde{y}_{1n})\ldots d\V_{jp^s}(\tilde{y}_{in}),\]
but $\F_{jp^s}(\epsilon_j)=\F_{jp^s}(\prod(\frac{1}{j}\V_j(1)-\frac{1}{jq}\V_{jq}(1))=\epsilon_1\in \BW{S/(jp^s)}{A}$
and 
\[\sum\epsilon_1\tilde{x}_nd\tilde{y}_{1n}\ldots\tilde{y}_{in}=\varphi^{-1}(e_1\sum_n x_n d y_{1n}\ldots dy_{in})=0.\]
This shows that $\varphi^{-1}{\alpha}$ is of type $(I)_j$ and thus in $N^i_S$. Finally we take the element
\[\beta=e_j\V_{p^s}(x)\partial_j[a]-e_j\V_{p^s}(x[a]^{p^s-1})\partial_j\V_{p^s}([a])\in \BWC{S/j}{1}{A}e_j\]
and get
\[\varphi^{-1}(\beta)=\V_j(\epsilon_1\V_{p^s}(\tfrac{1}{j^2}\tilde{x}))d\V_j([a])-
                                \V_j(\epsilon_1\V_{p^s}(\tfrac{1}{j^2}\tilde{x}[a]^{p^s-1}))d\V_j(\V_{p^s}([a])),\]
which is of type $(II)_{p^s}$, if $j=1$ and of type $(I)_j$, if $j>1$. 
This yields (\ref{1.13.1.5}) and thus an isomorphism (\ref{1.13.1}), which is
compatible with $\V_n$ by (\ref{1.13.2}).
\end{proof}

\begin{rmk}\label{1.14}
If $A$ is a $\Q$-algebra, we may choose any prime $p$. In particular, if $S$ is a finite truncation set, we can choose
a prime $p$, which is not contained in $S$. Then $\sP\cap S=\{1\}$ and we obtain by theorem \ref{1.13} an isomorphism 
\[\BWC{S}{\kaydot}{A}\stackrel{\simeq}{\longrightarrow}\prod_{s\in S}\Omega^\kaydot_A(\tfrac{1}{s}).\]
Therefore almost all statements for the de Rham-Witt complex follow immediately from the corresponding statements
for the de Rham complex.
\end{rmk}
\begin{rmk}\label{1.14.1}
Similar to the definition and construction of the $p$-typical de Rham-Witt complex of Deligne-Illusie 
(or Hesselholt-Madsen in the case of $\Z_{(p)}$-algebras), Langer-Zink give in \cite[Chapter 1]{LaZi04} a definition
and construction of a ($p$-typical) relative de Rham-Witt complex for a $R$-algebra $A$ (with $R$ a $\Z_{(p)}$-algebra),
denoted by $\WC{\kaydot}{\kaydot}{A/R}$. The differential in the $n$-th level is required to be $\W{n}{R}$-linear.
The universality of the $p$-typical de Rham-Witt complex gives a natural map 
$\WC{\kaydot}{\kaydot}{A}\to\WC{\kaydot}{\kaydot}{A/R}$.
\end{rmk}

\begin{lem}\label{1.14.1.1}
If $R$ is a perfect $\bb{F}_p$-algebra (i.e. $R^p=R$), with $p$ odd, or $R=\Q$ and $A$ is a $R$-algebra, then
the natural map 
\[\WC{\kaydot}{\kaydot}{A}\to\WC{\kaydot}{\kaydot}{A/R}\]
is an isomorphism 
\end{lem}
\begin{proof}
We have to check that the differential of a $p$-typical Witt complex is $\W{\kaydot}{R}$-linear. This follows from remark \ref{1.14},
if $R=\Q$. Else let $a\in R$ be an arbitrary element and take $b\in R$ such that $b^{p^n}=a$. Then
\[d\V^r([a])=d\V^r\F^r[b^{p^{n-r}}]=p^rd[b^{p^{n-r}}]=p^n[b^{p^{n-r}-1}]d[b]=0 \text{ in } \WC{n}{1}{A}.\] 
\end{proof}

There are two important examples of truncation sets, sets of the shape $\{1,p,\ldots,p^r\}$, these give the $p$-typical Witt complexes, and
the sets $\{1,2,\ldots, m\}$, $m\in \N$. The following definition says how a Witt complex on these sets look like.

\begin{defn}\label{1.14.2}
Let $A$ be a ring.  A \emph{restricted Witt complex over $A$} is a projective system of differential graded $\Z$-algebras
\[((E_m)_{m\in\N},\, \R : E_{m+1}\to E_m)\]
together with families of homomorphisms of graded rings
\[(\F_n: E_{nm+n-1}\longrightarrow E_{m})_{m,n\in \N}\]
and  homomorphisms of graded groups
\[(\V_n: E_m\longrightarrow E_{nm+n-1})_{m,n\in \N},\]
satisfying the following relations, for all $n,r \in\N$
\begin{enumerate}
\item $\R\F_n=\F_n\R^n,\quad \R^n\V_n=\V_n\R,\quad \F_1=\V_1=\id,\quad \F_n\F_r=\F_{nr},\quad \V_n\V_r=\V_{nr}.$
\item $\F_n\V_n=n$, and if $(n,r)=1$, then  $\F_r\V_n=\V_n\F_r$ on $E_{rm+r-1}$. 
\item $ \V_n(\F_n(x)y)=x\V_n(y)$, for all $x\in E_{nm+n-1}$, $y\in E_m$ and all $n\in\N$.
\item $\F_nd\V_n=d,\quad d\F_n=n\F_nd, \quad \V_nd=nd\V_n. $ 
\end{enumerate}
Furthermore, there is a homomorphism of projective systems of rings
\[(\lambda: \BW{m}{A}\longrightarrow E_m^0)_{m\in\N},\]
which  commutes with $\F_n$ and $\V_n$ and satisfies
\begin{enumerate}
\item[(v)] $\F_n d\lambda([a])=\lambda([a]^{n-1})d\lambda([a])$, for all $a\in A$ and $n\in \N$.
\end{enumerate}

A morphism of restricted Witt complexes over $A$ is a morphism of projective systems of dga's, compatible with $\F_n$, $\V_n$ and $\lambda$.
\end{defn}

Notice, that we have a forgetful functor from the category of Witt complexes over $A$ to the category of restricted Witt complexes over $A$,
$E(-)\mapsto (E(\{1,\ldots,m\}))_{m\in\N}$.

\begin{cor}\label{1.14.3}
Let $A$ be a $\Z_{(p)}$-algebra with $p$ an odd prime. Then $(\BWC{m}{\kaydot}{A})_{m\in\N}$ is the initial object in the category of restricted
Witt complexes.
\end{cor}
\begin{proof}
Let $E$ be a restricted Witt complex over $A$. Using the construction  of $\BWC{m}{\kaydot}{A}$ in the proof of proposition \ref{1.3} one checks 
immediately that we have a unique map of projective systems of dga's $\BWC{\kaydot}{\kaydot}{A}\to E$, commuting with $\V_n$ and $\lambda$.
It follows from \ref{1.7}, (d), that it automatically becomes a map of restricted Witt complexes.  
\end{proof}

We give some properties of the de Rham-Witt complex we need later on.

\begin{prop}\label{1.14.5}
Let $(A_i)_{i\in I}$ be a directed system, $A=\varinjlim A_i$.  For a finite truncation set $S$ define  
$E(S)=\varinjlim\BWC{S}{\kaydot}{A_i}\in\text{DGA}$ and for an arbitrary truncation set $S$ define 
$E(S)=\varprojlim_{S_0}E(S_0)\in\text{DGA}$, where the limit is taken over all finite truncation sets contained in $S$. Then
\[E: J\longrightarrow DGA,\quad S\mapsto E(S)\]
can be equipped in a natural way with the structure of a Witt complex and the natural map 
\[\BWC{-}{\kaydot}{A}\longrightarrow E(-)^{\kaydot}\]
is an isomorphism. 
\end{prop}

\begin{proof}
First of all notice, that direct and inverse limits exist in DGA (see \cite[2.11 Theorem]{Ku86}), thus the statement makes sense. Next we observe,
\[E(S)^0=\varinjlim\BW{S}{A_i}\cong\BW{S}{A}\]
(where the latter isomorphism is induced by the natural maps $\BW{S}{A_i}\to\BW{S}{A}$). With the same reasoning as in the proof of proposition \ref{1.3}
it is enough to show that $S\mapsto E(S)$ is a Witt complex on the category of finite truncation sets. For finite truncation sets $T\subset S$, 
$n,q\in\N$ and $\omega\in \varinjlim \BWC{S}{q}{A_j}$ define
\eq{1.14.5.1}{\omega_{|T}=\varphi_{j,T}({\omega_j}_{|T}),}
and
\eq{1.14.5.2}{\F_n(\omega)=\varphi_{j,S/n}( \F_n(\omega_j)),   }
where $\omega_j\in\BWC{S}{q}{A_j}$ is a representative of $\omega$ and $\varphi_{j,S}: \BWC{S}{q}{A_j}\to \varinjlim\BWC{S}{q}{A_i}$ is the natural map.
Furthermore define for $\omega\in \varinjlim\BWC{S/n}{q}{A_i}$
\eq{1.14.5.3}{\V_n(\omega)=\varphi_{j,S}(\V_n(\omega_j)). }
One immediately verifies that these maps are well defined and make $S\mapsto E(S)$ into a Witt complex over $A$. Thus, by the universal property of the
de Rham-Witt complex, we get a map
\[\BWC{-}{\kaydot}{A}\longrightarrow E(-)^\kaydot,\]
which is easily checked to be the inverse to the natural map
\[E(-)^\kaydot\longrightarrow \BWC{-}{\kaydot}{A},\]
coming from the maps
\[\BWC{S}{\kaydot}{A_i}\longrightarrow\BWC{S}{\kaydot}{A}, \text{ for } S \text{ finite},\]
hence the assertion. 
\end{proof}

\begin{prop}\label{1.15}
Let $k$ be a field of characteristic $p\neq 2$ and $A$ a $k$-algebra, $U\subset A$ a multiplicative system and $S$ a finite truncation set. 
Then, for all i, the natural map
\[\BW{S}{U^{-1}A}\otimes_{\BW{S}{A}}\BWC{S}{i}{A}\stackrel{\simeq}{\longrightarrow} \BWC{S}{i}{U^{-1}A}\]
is an isomorphism of $\BW{S}{U^{-1}A}$-modules.
\end{prop}
\begin{proof}
If $p=0$, this follows from the remark \ref{1.14}. If $p>2$, then by theorem \ref{1.13} and proposition \ref{0.9} we are reduced
to the $p$-typical situation and the statement becomes proposition 1.11. in chapter I of \cite{Il79}.
\end{proof}

\begin{prop}\label{1.16}
Let $k$ be a field of characteristic $p\neq 2$ and $A\longrightarrow B$ an \a'etale homomorphism of $k$-algebras and $S$ a finite truncation set.
Then, for all $i$,  the natural map
\[\BW{S}{B}\otimes_{\BW{S}{A}}\BWC{S}{i}{A}\stackrel{\simeq}{\longrightarrow} \BWC{S}{i}{B}\]
is an isomorphism of $\BW{S}{B}$-modules.
\end{prop}
\begin{proof}
Again in characteristic zero this follows from remark \ref{1.14} and the corresponding statement for K\"ahler differentials 
(see \cite{EGA IV}, corollaire (17.2.4)). If $p>2$, it is by theorem  \ref{1.13} and proposition \ref{0.9} enough to consider the $p$-typical
de Rham-Witt complex and the statement becomes proposition 1.14., in chapter I of \cite{Il79}.  
\end{proof}

The following Lemma was proved by Hesselholt in the $p$-typical case.

\begin{lem}[\cite{He04}, Lemma 1.2.2.]\label{1.17}
Let $p$ be an odd prime, $A$ a $\Z_{(p)}$-algebra and $I\subset A$ an ideal. Let $S$ be a finite truncation set and denote by $\sI_S$ the differential
graded ideal in $\BWC{S}{\kaydot}{A}$ generated by $\BW{S}{I}$. Then, for all $q$,
we have an isomorphism of $\BW{S}{A}$-modules
\[\BWC{S}{q}{A}/\sI^q_S\stackrel{\simeq}{\longrightarrow}\BWC{S}{q}{A/I}.\]
\end{lem}
\begin{proof}
Obviously we have a contravariant functor
\[\Gamma : J_S\longrightarrow \text{DGA},\quad S_0\mapsto \BWC{S_0}{\kaydot}{A}/\sI^\kaydot_{S_0}.\]
We want to show that this is a Witt complex over $A/I$ on the restricted category $J_S$. By lemma \ref{0.8.1} we have a natural isomorphism
\[\lambda : \BW{S_0}{A/I}\stackrel{\simeq}{\longrightarrow}\Gamma(S_0)^0=\BW{S_0}{A}/\BW{S_0}{I}.\]
Furthermore, it follows from the properties \ref{0.7}, (i),(iv) and \ref{1.7}, (c) and (d), that for $S_0\subset S$ and all $n,m$
\[\F_m(\sI_{S_0})\subset \sI_{S_0/m},\quad \V_n(\sI_{S_0/n})\subset \sI_{S_0}.\]
Thus the $\V_n$'s and $\F_m$'s on the de Rham-Witt complex over $A$, induce such maps on $\Gamma$. It follows that $\Gamma$ is a Witt-complex
over $A/I$ on $J_S$. Therefore we get a unique map $\bb{W}^S_{-}\Omega^\kaydot_{A/I}\longrightarrow \Gamma$ (where $\bb{W}^S_{-}\Omega^\kaydot_{A/I}$
denotes the de Rham-Witt complex restricted to $J_S$). One easily checks that this map and the natural map 
$\Gamma\longrightarrow \bb{W}^S_{-}\Omega^\kaydot_{A/I}$ are inverse to each other and we are done.  
\end{proof}

Finally, we want to introduce the great brother of the ghost map.

\begin{defn-prop}[cf. \cite{LaZi04}, 2.4]\label{1.18}
Let $p$ be an odd prime number, $S$ a finite truncation set and $A$ a $\Z_{(p)}$-algebra. Then there is a map
of graded $\BW{S}{A}$-algebras
\[\Gh:\BWC{S}{\kaydot}{A}\longrightarrow \prod_{s\in S}{\gh_s}_*\Omega_A^\kaydot,\]
which sends elements $\omega=wd\V_{n_1}([a_1])\ldots d\V_{n_r}([a_r])$, with $w\in \BW{S}{A}$ and $a_i\in A$, to
$(\Gh_s(\omega))_{s\in S}$ with
\eq{1.18.1}{\Gh_s(\omega)= \begin{cases}
               \gh_s(w)a_1^{\frac{s}{n_1}-1}\ldots a_r^{\frac{s}{n_r}-1}da_1\ldots da_r & \text{if } n_i|s,\text{ all }i\\
                                          0                                          & \text{else}.
                          \end{cases}}
(A sloppy way to say this is ``$\Gh_s(xdy_1\ldots dy_r)=\gh_s(x)\frac{1}{s}d\gh_s(y_1)\ldots\frac{1}{s}d\gh_s(y_r)$''.)
Gh is a natural transformation with respect to $A$ and if $A$ is a $\Q$-algebra, it is an isomorphism. 
We call $\Gh$ the \emph{ghost map on $\BWC{S}{\kaydot}{A}$}.
\end{defn-prop}
\begin{proof}
We define
\[\partial=(\partial_s)_s: \BW{S}{A}\longrightarrow \prod_{s\in S}{\gh_s}_*\Omega_A^1,
           \quad w=(w_s)_s\mapsto\bigl(\sum_{e|s}w_e^{\frac{s}{e}-1}dw_e\bigr)_s.\]
We observe that $\partial$ is a derivation. Indeed, it is enough to check this on $A=\Z[x_1,\ldots]$, but then
it follows immediately from $s\partial_s=d\gh_s$. We obtain a map of graded $\BW{S}{A}$-algebras
\[\Gh :\Omega^\kaydot_{\BW{S}{A}}\longrightarrow \prod_{s\in S}{\gh_s}_*\Omega^\kaydot_A,\]
which is clearly an isomorphism, if $A$ is a $\Q$-algebra. It remains to show that $\Gh$ factors through
$\BWC{S}{\kaydot}{A}$, i.e. the elements $(I)_n$, $(II)_n$ from the proof of proposition \ref{1.3} go to zero under $\Gh$.
Again we can assume $A=\Z[x_1\ldots]$ and then the assertion follows easily from 
$s^r\Gh_s(xdy_1\ldots dy_r)=\gh_s(x)d\gh_s(y_1)\ldots d\gh_s(y_r)$.

\end{proof}
\pagebreak[2]
\begin{properties}\label{1.19}
We have for all $s\in S$
\begin{enumerate}
\item $\Gh_s(\V_n(\omega))=\begin{cases}
                            n\Gh_{s/n}(\omega) & \text{if } n|s\\
                               0               & \text{else,}
                           \end{cases}$
       for $\omega\in \BWC{S/n}{\kaydot}{A}$.
\item $\Gh_s(\F_m(\omega))=\Gh_{sm}(\omega)$, for $\omega\in\BWC{S}{\kaydot}{A}$. 
\item $\Gh_s(\frac{d[a]}{[a]})=\frac{da}{a}$, for $a\in A^\times$.
\item $s\Gh_s\circ d=d \circ\Gh_s.$
\item If $\varphi:\BWC{S}{q}{A}\stackrel{\simeq}{\to}M(S)^q$ is the map from theorem \ref{1.13}, then for $j\in I_p$, $r\in\N_0$
      \[\Gh_{p^r}(\varphi(\omega)_j)=\Gh_{jp^r}(\omega).\]
\end{enumerate}
\end{properties}
\begin{proof}
(i), (iii) and (iv) follow immediately from (\ref{1.18.1}). It is enough to check (ii) for $\omega=d\V_r([a])$, $a \in A$. 
Write $c=(m,r)$ and take $l,k\in \Z$ with $km+lr=c $. Then, by property \ref{1.7}, (d), 
$\F_md\V_r([a])=kd\V_{r/c}([a]^{m/c})+l\V_{r/c}([a]^{\frac{m}{c}-1})d\V_{r/c}([a])$. Using (\ref{1.18.1}), we see
that $\Gh_s(\F_md\V_r([a]))$ is zero if $\frac{r}{c}\not|s$ or equivalently $r\not|sm$, and if $r|sm$
\[\Gh_s(\F_md\V_r([a]))=a^{\frac{(ms)}{r}-1}da=\Gh_{ms}(d\V_r([a])).\]
This proves (ii). Finally, (v) follows from 
\[\varphi(xdy_1 \ldots dy_q)_j=\F_j(x)\partial_j\F_j(y_1)\ldots\partial_j\F_j(y_q)_{|\sP\cap S/j}=\F_j(xd y_1\ldots dy_q)_{|\sP\cap S/j}\]
and (ii). 
\end{proof}
\begin{rmk}\label{1.20}
Let $R$ be a $\Z_{(p)}$-algebra and $A$ a $R$-algebra. Then Langer and Zink construct in \cite[2.4.]{LaZi04} in the same way as above a homomorphism
of projective systems of graded algebras
\[\Gh=(Gh_1,\Gh_p,\ldots,\Gh_{p^{n-1}}):\WC{n}{\kaydot}{A/R}\to \prod_{i=0}^{n-1}{\gh_{p^i}}_*\Omega^\kaydot_{A/R},\]
satisfying the relations in (\ref{1.18.1}) and \ref{1.19}, (i)-(iv) formulated for the $p$-typical situation. And they prove:
\end{rmk}

\begin{lem}[\cite{LaZi04}, Cor 2.18.]\label{1.21}
If $R$ is a $\Z_{(p)}$-algebra without $p$-torsion, then
\[\Gh:\WC{n}{\kaydot}{R[x_1,\ldots,x_d]/R}\longrightarrow \prod_{i=0}^{n-1}{\gh_{p^i}}_*\Omega^\kaydot_{R[x_1,\ldots,x_d]/R}\]
is injective.
\end{lem}
\begin{cor}\label{1.22}
If $A=\Z_{(p)}[x_1,\ldots,x_d]$ and $f\in A\setminus{\{0\}}$, then the natural map
\[\WC{n}{\kaydot}{A_f/\Z_{(p)}}\longrightarrow \WC{n}{\kaydot}{\Q\otimes_{\Z_{(p)}} A_f}\]
is injective.
\end{cor}
\begin{proof}
We have a commutative square
\[\xymatrix{\WC{n}{\kaydot}{\Z_{(p)}[x_1,\ldots,x_d]/\Z_{(p)}}\ar[r]\ar[d]^\Gh & \WC{n}{\kaydot}{\Q[x_1,\ldots,x_d]}\ar[d]^\Gh\\
             \prod_{i=0}^{n-1}{\gh_{p^i}}_*\Omega^\kaydot_{\Z_{(p)}[x_1,\ldots,x_d]/\Z_{(p)}}\ar[r] & 
                                                                   \prod_{i=0}^{n-1}{\gh_{p^i}}_*\Omega^\kaydot_{\Q[x_1,\ldots,x_d]}. }\]
Now the right map is an isomorphism, the left map and the lower map are injective (since $\Omega^\kaydot_{\Z_{(p)}[x_1,\ldots,x_d]/\Z_{(p)}}$ is a free 
$\Z_{(p)}$-module). Hence the upper map is injective too. Furthermore 
$\WC{n}{\kaydot}{A_f/\Z_{(p)}}=\W{n}{A}_{[f]}\otimes_{\W{n}{A}}\WC{n}{\kaydot}{A/\Z_{(p)}}$ (see \cite[(1.35)]{LaZi04}) and since
localization is flat, the statement follows. 
                                                                   
\end{proof}

\begin{rmk}\label{1.23}
Notice, that the above injectivity would be false, if we wrote the absolute de Rham-Witt complex at the left instead of the relative one.
For example $\WC{n}{1}{\Q}$ is zero, but $\WC{n}{1}{\Z_{(p)}}$ is not (see \cite[Example 1.2.4.]{HeMa04}).
\end{rmk}
\newpage

%
%
%
%
%
%
\section{A Residue Theorem}
In this section we will define a residue symbol on the de Rham-Witt complex over the function field $K$ of a curve, which is defined over an arbitrary
field $k$ of characteristic $\neq 2$.
This residue will generalize the one on the K\"ahler differentials on the one hand and the one defined by Witt on $\BW{S}{K}\times K^\times$ 
on the other hand (see \cite[2. Der Residuenvektor]{Wi36},  see also \cite[4.3.]{AnRo04}). To be able
to define a residue for curves defined over a field $k$, we first have to construct a trace on the de Rham-Witt complex for arbitrary finite field 
extensions. Here we proceed analogously to the construction of the trace for K\"ahler differentials 
(see \cite[\S 2. Der Begriff der Spur]{Ku64}, see also \cite[\S 16.]{Ku86}). The main goal of this section
is to prove a reciprocity law for the residue symbol. The proof will be done as in the classical case, by reduction to $\P^1$ via a trace formula
and a direct proof in this case. (See \cite[II, 12.]{Se88}, when the ground field is algebraically closed, and 
\cite[17.7 corollary]{Ku86}, for the general case.)\\
\\
Any field in this section has characteristic $\neq 2$.\\
\\
Before we can start with the construction of the trace, there are some technical points to settle. 
\begin{thm}[\cite{HeMa04}, Theorem 4.2.8.]\label{2.1}
Let $p$ be an odd prime, $A$ a $\Z_{(p)}$-algebra and $A[x]$ the polynomial ring over $A$. 
Then, for all $n, q\in \N$, the group $\WC{n}{q}{A[x]}$ is freely generated by elements of the following type
\eq{2.1.1}{a[x]^j, \text{ for } a\in \WC{n}{q}{A}, j\in\N_0,}
\eq{2.1.2}{b[x]^{j-1}d[x],\text{ for } b\in \WC{n}{q-1}{A}, j\in \N,}
\eq{2.1.3}{\V^s(a[x]^j),\text{ for } a\in \WC{n-s}{q}{A}, j\in I_p, s=1,\ldots,n-1,}
\eq{2.1.4}{d\V^s(b[x]^j),\text{ for } b\in\WC{n-s}{q-1}{A}, j\in I_p, s=1,\ldots,n-1.}
\end{thm}
\begin{rmk}\label{2.2}
In \cite[4.2.]{HeMa04} Hesselholt and Madsen say also how $d$, $\V$ and $\F$ act on these elements and they give rules for multiplication.
In fact, it is this result which enables Hesselholt and Madsen to define the Frobenius map form theorem \ref{1.9}. A more general statement,
for a polynomial ring with a finite number of variables, may be found in \cite[Theorem 1.2.1.]{He04}. (See also \cite[Theorem 2.8]{LaZi04} for the
corresponding result in the case of the relative version of the $p$-typical de Rham-Witt complex.) 

\end{rmk}
\begin{lem}\label{2.3}
\begin{enumerate}
\item Let $A\longrightarrow B$ be  an \'etale morphism of $\bb{F}_p$-algebras. Then for all $n\in \N$, the map $\W{n}{A}\longrightarrow\W{n}{B}$ is 
      \'etale and  
                \[\W{n+1}{B}\otimes_{\W{n+1}{A}}\F_*\W{n}{A}\stackrel{\simeq}{\longrightarrow} \W{n}{B},\quad b\otimes a \mapsto \F(b)a\]
      is an isomorphism of $\W{n}{A}$-algebras.
\item Let $A$ be a $\bb{F}_p$-algebra and $U\subset A$ a multiplicative system. Then the map
       $\W{n}{A}\to\W{n}{U^{-1}A}$ is flat and
         \[\varphi:\W{n+1}{U^{-1}A}\otimes_{\W{n+1}{A}}\F_*\W{n}{A}\stackrel{\simeq}{\longrightarrow} \W{n}{U^{-1}A},\quad q\otimes a\mapsto \F(q)a\]
        is an isomorphism of $\W{n}{A}$-algebras.
\end{enumerate}
\end{lem}
\begin{proof}
(i) may be found for example in \cite[A.8, A.11]{LaZi04}. To prove (ii), we define a map 
$\psi : \W{n}{U^{-1}A}\longrightarrow \W{n+1}{U^{-1}A}\otimes_{\W{n+1}{A}}\F_*\W{n}{A}$, by 
\eq{2.3.1}{\psi\bigl(\V^j([\tfrac{a}{u}])\bigr)= \tfrac{1}{[u]}\otimes\V^j([a][u]^{p^{j+1}-1}),\text{ and } 
                                                                 \psi\bigl(\sum_j\V^j([q_j])\bigr)=\sum_j\psi(\V^j([q_j])).}
We have to check that this is well defined. So take $a,b\in A$ and $t,u,v\in U$  with $t(av-ub)=0$, then
\[\psi(\V^j([\tfrac{a}{u}]))=\tfrac{1}{[tuv]}\otimes\V^j([atv][utv]^{p^{j+1}-1})=\psi(\V^j([\tfrac{b}{v}])).\] 
Obviously, $\varphi\circ\psi=\id$. Thus it remains to show $\psi\circ\varphi=\id$ and for this it is enough to consider elements of the form
$\alpha=\V^j([\frac{a}{u}])\otimes\V^i([b])$. Since we are in a $\bb{F}_p$-algebra, we have $\V\circ\F=\F\circ\V=p$, thus, if $t=\min(i,j)$, we
obtain with the help of property \ref{0.7}, (v)
\[\varphi(\alpha) =  \V^{j}(\F([\tfrac{a}{u}]))\V^i([b]) = p^{t}\V^{j+i-t}([\tfrac{a}{u}]^{p^{i+1-t}}[b]^{p^{j-t}})
                                                                        =  \V^{i+j}(\tfrac{[a]^{p^{i+1}}[b]^{p^{j}}}{[u]^{p^{i+1}}}).\]
So we get
\begin{eqnarray*}\psi(\varphi(\alpha)) & = & \tfrac{1}{[u]^{p^{i+1}}}\otimes \V^{i+j}([a]^{p^{i+1}}[b]^{p^j}([u]^{p^{i+1}})^{p^{j+i+1}-1})\\
                                       & =  & \tfrac{1}{[u]^{p^{i+1}}}\otimes \F(\V^j([a][u]^{p^{i+j+1}-1}))\V^i([b]) =  \alpha.             
\end{eqnarray*}
Furthermore, the flatness of $\W{n}{A}\to\W{n}{U^{-1}A}$, follows from
$\W{n}{U^{-1}A}=U_n^{-1}\W{n}{A}$, with $U_n=\{(a_0,\ldots,a_{n-1})\in \W{n}{A}|a_0\in U\}$ (see the proof of proposition 1.11. in \cite{Il79})
and the fact that localization is flat.
\end{proof}

The following lemma is a slight adjustment of \cite[I, 3.2., 3.4., 3.21.1.]{Il79} to our situation.  

\begin{lem}\label{2.4}
Let $k$ be a field of characteristic $p\ge 3$ and $i\ge 0$, $n\ge 1$. Then
\eq{2.4.1}{\ker( R:\WC{n+1}{i}{k}\longrightarrow\WC{n}{i}{k})=\ker(p:\WC{n+1}{i}{k}\longrightarrow\WC{n+1}{i}{k})=\V^n\WC{1}{i}{k}+d\V^n\WC{1}{i-1}{k},}
where $p$ means multiplication by $p$, and
\eq{2.4.2}{\ker(\V:\WC{n}{i}{k}\longrightarrow \WC{n+1}{i}{k})= d\V^{n-1}\WC{1}{i-1}{k}.}
\end{lem}
\begin{proof}
Illusie already proved the statement for smooth $\bb{F}_p$-schemes, in particular for $\A^r_{\bb{F}_p}$. We write $A=\bb{F}_p[x_1,\ldots,x_r]$
and $Q$ for its quotient field. From (\ref{2.4.1}) for $A$ instead of $k$ it follows that we have exact sequences of $\W{n+1}{A}$-modules
\[\xymatrix{0\ar[r]  & \V^n\WC{1}{i}{A}+d\V^n\WC{1}{i-1}{A}\ar[r] & \WC{n+1}{i}{A}\ar[r]^\R & \R_*\WC{n}{i}{A}\ar[r]  & 0   }\]
and
\[\xymatrix{0\ar[r]  & \V^n\WC{1}{i}{A}+d\V^n\WC{1}{i-1}{A}\ar[r] & \WC{n+1}{i}{A}\ar[r]^p & p\WC{n+1}{i}{A}\ar[r]  & 0   }.\]
Tensoring with $\W{n+1}{Q}$ together with $\W{n+1}{Q}\otimes\WC{n+1}{i}{A}=\WC{n+1}{i}{Q}$ (proposition \ref{1.15}), shows that 
the abelian groups $\ker\, R_Q$ and $\ker\,p_Q$ are generated by
\[a\V^n(\alpha)+bd\V^n(\beta)=\V^n(\F^n(a)\alpha)+d\V^n(\F^n(b)\beta)-\V^n(\F^n(db)\beta),\]
with $a,b\in \W{n+1}{Q}$, $\alpha\in \WC{1}{i}{A}$ and $\beta\in\WC{1}{i-1}{A}$. This shows (\ref{2.4.1}) for $Q$.
Furthermore, by combining proposition \ref{1.15} with lemma \ref{2.3}, (ii), we obtain an isomorphism
\[\W{n+1}{Q}\otimes_{\W{n+1}{A}}\F_*\WC{n}{i}{A}\stackrel{\simeq}{\longrightarrow}\WC{n}{i}{Q}, \quad q\otimes \alpha\mapsto F(q)\alpha.\]
This yields the following commutative  diagram with exact columns
\[\xymatrix{                                                                                                          &  & 0\ar[d]\\
            \W{n+1}{Q}\otimes_{\W{n+1}{A}} \ker\,\V_A\ar[rr]\ar[d]                                                    &  & \ker\,\V_Q\ar[d]\\
            \W{n+1}{Q}\otimes_{\W{n+1}{A}}\F_*\WC{n}{i}{A}\ar[rr]^(.7){\simeq\,(\ref{2.3}, (ii))}\ar[d]_{\id\otimes \V} &  & \WC{n}{i}{Q}\ar[d]^\V\\
            \W{n+1}{Q}\otimes_{\W{n+1}{A}}\V\WC{n}{i}{A}\ar[rr]^(.7){\simeq\,(\ref{1.15})}\ar[d]                          &  &  \V\WC{n}{i}{Q}\ar[d]\\
                                               0                                                                      &  &        0.
            }\]
{From} (\ref{2.4.2}) for $A$ it follows that we have a surjection
\[\W{n+1}{Q}\otimes_{\W{n+1}{A}}\F_*d\V^{n-1}\WC{1}{i-1}{A}\surj \ker\, \V_Q.\]
Thus the abelian group $\ker\, V_Q$ is generated by
\[\F(q)d\V^{n-1}(\alpha)=d\V^{n-1}(\F^n(q)\alpha)-\V^{n-1}(\F^{n-1}d\F(q)\alpha),\text{ for } q\in \W{n+1}{Q},\alpha\in\WC{1}{i-1}{A}.\]
But $\V^{n-1}(\F^{n-1}d\F(q)\alpha)=p\V^{n-1}(\F^n(dq)\alpha)=0$ in $\WC{n}{i}{Q}$, by what we saw above and this proves (\ref{2.4.2}) for $Q$. 
Now assume $k\supset Q$ is a finite separable field extension, then replacing $A$ by $Q$ and $Q$ by $k$ in the discussion above and
quoting proposition \ref{1.16} instead of \ref{1.15} and the part (i) of lemma \ref{2.3} instead of (ii), we see that the lemma holds in this case.

We observe: \emph{let $k\supset Q$ be an algebraic separable (resp. $k\supset \bb{F}_p$ any) field extension and
$L$ a finite separable (resp. finitely generated) intermediate field and $\alpha\in \WC{n}{i}{L}$ with $\alpha\mapsto 0$ under 
$\WC{n}{i}{L}\longrightarrow\WC{n}{i}{k}$. Then there is a finite separable (resp. finitely generated) intermediate field $L\subset L'\subset k$ with
$\alpha\mapsto 0$ under $\WC{n}{i}{L}\longrightarrow\WC{n}{i}{L'}$. }(Indeed, we may write $k=\varinjlim L$, where the limit is over all
finite separable field extensions $L\supset Q$ contained in $k$ (resp. over all finitely generated field extensions $L\supset \bb{F}_p$ contained in $k$)
and then the observation follows from proposition \ref{1.14.5}.)

Now let $k$ be  finitely generated over $\bb{F}_p$. Then there exists a transcendence basis $\{x_1,\ldots,x_r\}$, such that 
$k\supset \bb{F}_p(x_1,\ldots,x_r)=Q$ is an algebraic separable field extension. Take $\alpha\in \ker(R:\WC{n+1}{i}{k}\longrightarrow \WC{n}{i}{k})$.
Clearly, there is a finite separable intermediate field $Q\subset L\subset k$ and an element $\alpha'\in \WC{n+1}{i}{L}$, which maps to $\alpha$.
Thus $\R(\alpha')\mapsto 0$ under  $\WC{n}{i}{L}\longrightarrow \WC{n}{i}{k}$ and the assertion in this case follows from the observation
and by what we saw above. 
The statements for $p$ and $\V$ follow similarly. Finally, let $\bb{F}_p\subset k$ be an arbitrary field of characteristic $p$, 
then we replace in the discussion above  $Q$ by $\bb{F}_p$ and ``finite separable'' by ``finitely generated'' and the statement is proved. 
\end{proof}

The following proposition is essential to prove that the trace we will define in \ref{2.7} is well defined in the case of purely inseparable field
extensions of degree $p$.

\begin{prop}\label{2.5}
Let $k$ be a field of characteristic $p\ge 3$ and $L=k[x]/(x^p-b)$, for some $b\in k\setminus k^p$. Then
\ml{2.5.1}{\ker(\WC{n+1}{i}{k}\longrightarrow \WC{n+1}{i}{L})\\
          =\sum_{r=1}^n(\V^r([b])d\V^n\WC{1}{i-1}{k}+d\V^r([b])d\V^n\WC{1}{i-2}{k})+ 
                                                                    (\V^n\WC{1}{i-1}{k}+d\V^n\WC{1}{i-2}{k})d[b].}
\end{prop}
\begin{proof}
We denote by $\iota: \WC{n+1}{i}{k}\longrightarrow \WC{n+1}{i}{L}$ the map induced by the inclusion $k\subset L$. First we observe that 
by lemma \ref{2.4} the right hand side of the equation is contained in the left hand side, for example 
\[\iota(\V^r([b])d\V^n(\omega))=\V^r(\F([x]))d\V^n(\iota(\omega))= p\V^{r-1}([x])d\V^n(\iota(\omega))=0\]
or
\[\iota(\V^n(\omega)d[b])=p\V^n(\iota(\omega))[x]^{p-1}dx=0.\]
Thus our task is to show the other inclusion. Let $I$ be the kernel of $k[x]\to L$ and $\sI^i$ the kernel 
of the map $\WC{n+1}{i}{k[x]}\to\WC{n+1}{i}{L}$.
Then, since we have an inclusion $\WC{n+1}{i}{k}\subset\WC{n+1}{i}{k[x]}$ (by theorem \ref{2.1}), $\ker\,\iota$ equals $\sI^i\cap \WC{n+1}{i}{k}$.
Now by lemma \ref{1.17}, $\sI^i$ is generated by $\W{n+1}{I}$ and $d\W{n+1}{I}$. We claim that every element in $\W{n+1}{I}$ may be
written as a sum of elements
\[v_r(a,j):=V^r(a[x]^j([x]^p-[b])), \text{ for } a\in\W{n+1-r}{k}, r=0,\ldots,n, j\ge 0.\]
(To see this, we define an operator $\varphi$ in the following way: take $w=(w_0,\ldots,w_n)\in\W{n+1}{I}$, write $m=\min\{l|w_l\neq 0\}$ and
$w_m=\sum_ja_jx^j(x^p-b)$, for $a_j\in k$. Then  $u(w)=\sum_j\V^m([a_j][x]^j([x]^p-[b]))$ is in $\W{n+1}{I}$ (by lemma \ref{0.8.1}) and
$\varphi(w):=w-u(w)\in V^{m+1}\W{n}{I}$. It follows that there is a $N$ such that $\varphi^N(w)=0$, hence the assertion.)
Thus we can write every Element of $\sI^i$ as
\[\sum_j\sum_{r=0}^n\omega_{jr}v_r(a_{jr},j)+\eta_{ir}dv_r(b_{jr},j),\]
for $\omega_{ir}\in\WC{n+1}{i}{k[x]}$, $\eta_{ir}\in \WC{n+1}{i-1}{k[x]}$ and $a_{ir},b_{ir}\in \W{n+1-r}{k}$.
Now $k[x]\longrightarrow k$, $x\mapsto 0$, defines a map $s: \WC{n+1}{i}{k[x]}\longrightarrow \WC{n+1}{i}{k}$ such that  
the inclusion followed by $s$ is the identity. We obtain $\ker\, \iota=\sI^i\cap\WC{n+1}{i}{k}\subset s(\sI^i)$. 
Thus every element in $\ker\, \iota$ may be written as sum of elements of the following type (we write them directly in the form of theorem \ref{2.1})
\[\omega\V^r(a([x]^p-[b]))=\begin{cases}
                           a\omega[x]^p-a \omega[b] & \text{for } r=0\\
                            \V^{r-1}(V(F^r(\omega)a)[x])-\omega \V^r(a[b]) & \text{for } r\ge 1
                           \end{cases}\]
with $\omega\in\WC{n+1}{i}{k}$, $a\in \W{n+1-r}{k}$, and
\[\eta dV^r(c([x]^p-[b]))=\]
                          \[\begin{cases}
                          \eta dc[x]^p+\eta cp[x]^{p-1}d[x]-\eta d(c[b]) & \text{for } r=0\\
                          \eta d\V(c)[x]+\eta\V(c)d[x]-\eta d\V(c[b]) & \text{for } r=1\\
                          (-1)^{i-1}\Bigl(d\V^{r-1}(\F^{r-1}(\eta)\V(c)[x])-\V^{r-1}(\F^{r-1}(d\eta)\V(c)[x])\Bigr)-\eta d\V^r(c[b]) & \text{for } r>1
                          \end{cases}\] 
with $\eta \in\WC{n+1}{i-1}{k}$, $c\in\W{n+1-r}{k}$. Thus we can write every element $\alpha\in\ker\,\iota$ in the following way \pagebreak[3]
\[\alpha=\sum_{i_0}a_{0i_0}\omega_{0i_0}[x]^p-\sum_{i_0}a_{0i_0}\omega_{0i_0}[b] + 
    \sum_{i_1}\V(\F(\omega_{1i_1})a_{1i_1})[x]-\sum_{i_1}\omega_{1i_1}\V(a_{1i_1}[b])\]
\[+\sum_{r=2}^n\biggl(\sum_{i_r}\V^{r-1}\bigl(\V(\F^r(\omega_{ri_r})a_{ri_r}))[x]\bigr)-\sum_{i_r}\omega_{ri_r}\V^r(a_{ri_r}[b])\biggr)\]
\[+\sum_{j_0}\eta_{0j_0}dc_{0j_0}[x]^p+\eta_{0j_0}c_{0j_0}p[x]^{p-1}d[x]-\eta_{0j_0}d(c_{0j_0}[b])\]
\[+ \sum_{j_1}\eta_{1j_1}d\V(c_{1j_1})[x]+\eta_{1j_1}\V(c_{1j_1})d[x]-\eta_{1j_1}d\V(c_{1j_1}[b])\]
\[+\sum_{r=2}^n\biggl(\sum_{j_r}(-1)^{i-1}d\V^{r-1}\bigl(\V(c_{rj_r})\F^{r-1}(\eta_{rj_r})[x]\bigr)-\]
\[(-1)^{i-1}\V^{r-1}\bigl(\V(c_{rj_r})\F^{r-1}(d\eta_{rj_r})[x]\bigr)-\eta_{rj_r}d\V^r(c_{rj_r}[b])\biggr).\]
By theorem \ref{2.1} such an element is in $\sI^i\cap \WC{n+1}{i}{k}$ if and only if the coefficients of the $x$-terms cancel out each other, that
is
\begin{enumerate}
\item[(a)] $\sum_{i_1}\V(\F(\omega_{1i_1})a_{1i_1})+ \sum_{j_1}\eta_{1j_1}d\V(c_{1j_1})=0$ in $\WC{n+1}{i}{k}$,
\item[(b)] $\sum_{i_0}a_{0i_0}\omega_{0i_0}+ \sum_{j_0}\eta_{0j_0}dc_{0j_0}=0$ in $\WC{n+1}{i}{k}$,
\item[(c)] $\sum_{j_1}\eta_{1j_1}\V(c_{1j_1})=0$ in $\WC{n+1}{i-1}{k}$,
\item[(d)] $p\sum_{j_0}\eta_{0j_0}c_{0j_0}=0$ in $\WC{n+1}{i-1}{k}$,
\item[(e)] $\sum_{i_r}\V(\F^r(\omega_{ri_r})a_{ri_r})-(-1)^{i-1}\sum_{j_r}\V(c_{rj_r})\F^{r-1}(d\eta_{rj_r})=0$ in $\WC{n+2-r}{i}{k}$, $r\ge 2$,
\item[(f)] $\sum_{j_r}\V(c_{rj_r})\F^{r-1}(\eta_{rj_r})=0$ in $\WC{n+2-r}{i-1}{k}$, $r\ge 2$.
\end{enumerate}
If we write 
\[\gamma_r=\sum_{i_r}\omega_{ri_r}\V^r(a_{ri_r}[b])+\sum_{j_r}\eta_{rj_r}dV^r(c_{rj_r}[b])\]
we obtain 
\[\alpha=-\sum_{r=0}^n\gamma_r\]
and thus it is enough to show that $\gamma_r$  is contained in the right hand side of (\ref{2.5.1}), for $r\ge 0$. 
We start with $\gamma_0$. By (b)
\[\gamma_0=\sum_{j_0}\eta_{0j_0}c_{0j_0}d[b]\]
and from (d) and  (\ref{2.4.1}) we obtain $\gamma_0\in \V^n\WC{1}{i}{k}d[b]+d\V^n\WC{1}{i-1}{k}d[b]$. Next $\gamma_1$. Writing (c) as
$\V\bigl(\sum_{j_1}\F(\eta_{1j_1})c_{1j_1}\bigr)=0$, we see from (\ref{2.4.2}) that there is a $\beta\in \WC{1}{i-1}{k}$ such that
\eq{2.5.2}{\sum_{j_1}\F(\eta_{1j_1})c_{1j_1}=d\V^{n-1}(\beta).}
Furthermore, (c) tells us $\sum_{j_1}\eta_{1j_1}d\V(c_{1j_1})=(-1)^i\sum_{j_1}\V(c_{1j_1})d\eta_{1j_1}$, thus (a) becomes
\[\V\bigl(\sum_{i_1}\F(\omega_{1i_1})a_{1i_1}+ (-1)^i \sum_{j_1}c_{1j_1}\F(d\eta_{1j_1})\bigr)=0\]
and by (\ref{2.4.2}) there exists an $\beta'\in\WC{1}{i-1}{k}$ such that
\eq{2.5.3}{\sum_{i_1}\F(\omega_{1i_1})a_{1i_1}+(-1)^i\sum_{j_1}c_{1j_1}\F(d\eta_{1j_1})=d\V^{n-1}(\beta').} 
All in all we get by (\ref{2.5.2}) and (\ref{2.5.3})
\[\gamma_1=\sum_{i_1}\omega_{1i_1}\V(a_{1i_1}[b])+(-1)^{i-1}\sum_{j_1}d\V(\F(\eta_{1j_1})c_{1j_1}[b])
                                                 +(-1)^i\sum_{j_1}\V(c_{1j_1}\F(d\eta_{1j_1})[b])\]
\[=\V(d\V^{n-1}(\beta')[b])+(-1)^{i-1}d\V(d\V^{n-1}(\beta)[b])=\V([b])d\V^n(\beta')+d\V([b])d\V^n((-1)^{i-1}\beta).\]
Finally $\gamma_r$, for $r\ge 2$. We see from (e) and (\ref{2.4.2}) that there is a $\beta\in \WC{1}{i-1}{k}$ with
\eq{2.5.4}{\sum_{i_r}\F^r(\omega_{ri_r})a_{ri_r}+(-1)^i \sum_{j_r}c_{rj_r}\F^r(d\eta_{rj_r})=d\V^{n-r}(\beta).}
By (f) and (\ref{2.4.2}) there is a $\beta'\in \WC{1}{i-2}{k}$ such that
\eq{2.5.5}{\sum_{j_r}c_{rj_r}\F^r(\eta_{rj_r})=d\V^{n-r}(\beta').}
Using  (\ref{2.5.4}) and (\ref{2.5.5}), we obtain
\[ \gamma_r= \V^r\bigl(\sum_{i_r}\F^r(\omega_{ri_r})a_{ri_r}[b]\bigr)+
  \sum_{j_r}(-1)^i\V^r(c_{rj_r}[b]\F^r(d\eta_{rj_r})) +(-1)^{i-1}\sum_{j_r}d\V^r(\F^r(\eta_{rj_r})c_{rj_r}[b])\]
\[=\V^r\bigl(d\V^{n-r}(\beta)[b]\bigr)+(-1)^{i-1}d\V^r\bigl(d\V^{n-r}(\beta')[b]\bigr)=
\V^r([b])d\V^n(\beta)+d\V^r([b])d\V^n((-1)^{i-1}\beta')\]
and we are done.
\end{proof}

\begin{defn}[cf. \cite{Il79}, Proposition I.3.4.]\label{2.6}
Let $k$ be a field of characteristic $p\ge 3$ and $S$ a finite truncation set. We just write $pS$ instead of $S\cup pS$
(which is again a truncation set). Then we define a map of dga's
\[\underline{p}: \BWC{S}{\kaydot}{k}\longrightarrow \BWC{pS}{\kaydot}{k},\quad \omega\mapsto p\widetilde{\omega},\]
where $\widetilde{\omega}$ is any lifting of $\omega$ to $\BWC{pS}{\kaydot}{k}$. By lemma \ref{2.4} $\underline{p}$
is well defined and injective in the $p$-typical case, $\underline{p}: \WC{n}{\kaydot}{k}\to\WC{n+1}{\kaydot}{k}$,
and thus by theorem \ref{1.13} also in general.
\end{defn}
The following construction of the trace in the case of a finite field extension, by defining it separately for separable field extensions and for
purely inseparable ones of degree $p$ and then in the general case by combining these two, is the way it was done for K\"ahler differentials
in \cite[\S2.]{Ku64}. We recall the definition of the trace on K\"ahler differentials for separable and purely inseparable field extensions, to point out
that our construction harmonizes with the one of Kunz. If $L\supset k$ is a separable extension, the trace is given by 
$\Tr=\Tr_{L/k}\otimes 1: L\otimes_k\Omega^n_k=\Omega^n_L\to\Omega^n_k$. If $L\supset k$ is purely inseparable of degree $p$, we can write
$L=k[x]/(x^p-b)$, $b\in k\setminus k^p$. Then each $\omega\in \Omega^n_L$ may be written in the form 
$\sum_{i=0}^{p-1}\alpha_ix^i+\sum_{i=1}^{p-1}\beta_ix^idx$, with $\alpha_i\in\Omega^n_k$ and $\beta_i\in\Omega^{n-1}_k$ and then, by definition,
$\Tr(\omega)=\beta_{p-1}db$.

\begin{thm}\label{2.7}
Let $k\subset L$ be a finite field extension of characteristic $p\neq 2$ and $S$ a finite truncation set. 
Then there is a map of differential graded $\BWC{S}{\kaydot}{k}$-modules
\[\Tr=\Tr_{L/k}: \BWC{S}{\kaydot}{L}\longrightarrow \BWC{S}{\kaydot}{k},\]
satisfying the following properties (usually we just write $\Tr$ but depending on the situation, we also use
$\Tr_{L/k}$, $\Tr_S$  or $\Tr^i$ to indicate that we are in degree $i$, 
combinations of these notations may also occur, but we try to omit $Tr_{L/k,S}^q$).
\begin{enumerate}
\item $\Tr^0$ is the trace from proposition \ref{0.12}. For $S=\{1\}$, $\Tr_S$ is the old trace on K\"ahler differentials.
\item If $k\subset L$ is a separable field extension, then we may identify (by proposition \ref{1.16})
        $\BW{S}{L}\otimes\BWC{S}{i}{k}=\BWC{S}{i}{L}$ and then $\Tr$ is given by
           \[\Tr=\Tr^0\otimes \id :\BW{S}{L}\otimes\BWC{S}{i}{k}\longrightarrow \BWC{S}{i}{k}.\]
\item If $k\subset L$ is purely inseparable of degree $p$ and $\omega\in\BWC{S}{i}{L}$, then 
      $\underline{p}(\omega)$ is in the image of the natural map $\BWC{pS}{i}{k}\to\BWC{pS}{i}{L}$ and if we take an
     $\alpha\in\BWC{pS}{i}{k}$ which maps to $\underline{p}(\omega)$, we have
         \[\Tr(\omega)=\alpha_{|S}.\] 
\item If $k\subset E\subset L$ are finite field extensions, then
                  \[\Tr_{L/k}=\Tr_{E/k}\circ\Tr_{L/E}.\]
\item $\Tr$ commutes with $\V_n$, $\F_n$ and restriction maps.
       
\end{enumerate}
\end{thm}
\begin{proof}
First assume $k\subset L$ separable. Denote by $\sigma_j$, $j=1,\ldots, d$, the embeddings of $L$ into its algebraic 
closure, which leave $k$ invariant. These induce maps on the level of the de Rham-Witt complex, again denoted by 
$\sigma_j$. For $\omega\in\BWC{S}{i}{L}$, we define 
\[\Tr(\omega)=\sum_j\sigma_j(\omega).\]
By example \ref{0.13} this gives (ii) and (i). Obviously it is a  map of differential graded 
$\BWC{S}{\kaydot}{k}$-modules and fulfills (iv) and (v).  

Now we consider the case of a purely inseparable field extension $k\subset L$ of degree p. We want to define the trace 
in this case by (iii). To be able to do this, we have to show 
\eq{2.7.1}{\im(\underline{p}_L)\subset \im(\BWC{pS}{\kaydot}{k}\to\BWC{pS}{\kaydot}{L})}
and check that the trace is independent of the choices. 
By theorem \ref{1.13} it is enough to consider the $p$-typical case. We can find an element $b\in k\setminus k^p$
such that we may identify
\[L=k[x]/(x^p-b).\]
If we denote by $y$ the image of $x$ in $L$, we see by theorem \ref{2.1} that every $\omega\in \WC{n}{i}{L}$ 
may be written in the following way
\[\omega=\sum_{j=0}^{p-1} a_{0j}[y]^j + \sum_{j=1}^{p-1}b_{0j}d[y]^j+b_{0p}[y]^{p-1}d[y]
+\sum_{s=1}^{n-1}\sum_{j=1}^{p-1}\V^s(a_{sj}[y]^j)+d\V^s(b_{sj}[y]^j),\]
with $a_{sj}\in \WC{n-s}{i}{k}$ and $b_{sj}\in\WC{n-s}{i-1}{k}$. Recalling from remark \ref{0.11}, (i), that 
$\V\circ F=p$, we obtain
\ml{2.7.2}{\underline{p}(\omega)= \sum_{j=0}^{p-1} \widetilde{a}_{0j}\V([b]^j)+
                                         \sum_{j=1}^{p-1}\widetilde{b}_{0j}d\V([b]^j)+\widetilde{b}_{0p}d[b]\\
                                   + \sum_{s=1}^{n-1}\sum_{j=1}^{p-1}\V^s\bigl(\widetilde{a}_{sj}\V([b]^j)\bigr)+
                                      d\V^s\bigl(\widetilde{b}_{sj}\V([b]^j)\bigr), }
where $\widetilde{}$ denotes any lifting from the level $n-s$ to the level $n-s+1$. This shows (\ref{2.7.1}). So take an 
$\alpha\in \WC{n+1}{i}{k}$, which maps to $\omega$ and define $\Tr(\omega)=\alpha_{|n}$. By proposition \ref{2.5} this
is independent of the choice of $\alpha$. Example \ref{0.13} shows that this coincides in degree $0$ 
with the trace on Witt vectors and for $n=1$ it is by (\ref{2.7.2})  the old trace for K\"ahler differentials, hence (i).
No doubt, $\Tr$ is a map of differential graded modules and satisfies (v).

Finally, the general case. Here we copy the proof of \cite[2.3.]{Ku64} almost word by word.
So let $k\subset L$ be any finite field extension of characteristic $p\ge 3$. Then we find a tower of field extensions
\eq{2.7.2.1}{k=k_0\subset k_1\subset \ldots \subset k_m=L,}
with $k_{j+1}/k_j$ either separable or purely inseparable of degree $p$. Then we define
\[\Tr_{L/k}=\Tr_{k_1/k_0}\circ\Tr_{k_2/k_1}\circ\ldots\circ\Tr_{k_m/k_{m-1}}.\]
Of course, we have to show that this is independent of the choice of the $k_j$. (Notice that the theorem follows, once this 
is proved.) Again we consider first some special cases.

\emph{First case. Assume there is an $a\in L$, with $a^p\in k$ and $a\not\in k$ such that 
$L\supset L'=k(a)$ is finitely separable. Then there is a primitive element $b\in L$ with $L=L'(b)$.
Since $L'(b)=L'(b^p)$, we may assume $b$ is separable over $k$. Define $L''=k(b)$. Then}
\eq{2.7.3}{\Tr_{L'/k}\circ\Tr_{L/L'}=\Tr_{L''/k}\circ\Tr_{L/L''}.}
Indeed, take $\omega\in \BWC{S}{i}{L}$, since $L''\subset L$ is purely inseparable of degree $p$ we find an 
$\alpha\in \BWC{pS}{i}{L''}$ with $\alpha\mapsto \underline{p}(\omega)$ under $\BWC{pS}{i}{L''}\to \BWC{pS}{i}{L}$.
Now denote by $\overline{L}$ the algebraic closure of $L$, then 
\[\Hom_{L'}(L,\overline{L})=
              \{\sigma_1,\ldots,\sigma_d\}\stackrel{\simeq}{\longrightarrow}\Hom_k(L'',\overline{L}),\quad
                                                                                        \sigma\mapsto \sigma_{|L''}\]
defines a bijection (since $L'\subset L$ and $k\subset L''$ are separable and $[L:L']=[L'':k]$).
Therefore
\[\Tr_{L''/k}\Tr_{L/L''}(\omega)=\sum_j\sigma_j(\alpha_{|S}).\]
On the other hand $\sum_j\sigma_j(\alpha)\mapsto \underline{p}\bigl(\sum_j\sigma_j(\omega)\bigr)$ under 
$\BWC{pS}{i}{k}\to\BWC{pS}{i}{L'}$, hence (\ref{2.7.3}).

\emph{Second case. Assume $L=k(a,b)$ with $a^p,b^p\in k$, $a,b\not\in k$ and $[L:k]=p^2$. 
      Write $L'=k(a)$ and $L''=k(b)$. Then}
\eq{2.7.4}{\Tr_{L''/k}\circ\Tr_{L/L''}=\Tr_{L'/k}\circ\Tr_{L/L'}.}
This follows from the definition.

\emph{General case.} If we have a tower of field extensions as in (\ref{2.7.2.1}), we can assume by the first case
that there is a natural number $r$ such that the $k_{j+1}/k_j$ are separable extensions for $j=0,\ldots, r$ and
$k_{j+1}/k_j$ are purely inseparable of degree $p$, for $j=r+1,\ldots,m-1$. Since (iv) holds for separable extensions,
we are reduced to $k\subset L$ purely inseparable of degree $p^m$. We make induction over $m$. For $m=1$ there is
nothing to prove. Now assume we have two towers of field extensions
\[k=k_0\subset k_1\subset \ldots \subset k_m=L,\]
\[k=\overline{k}_0\subset \overline{k}_1\subset \ldots \subset \overline{k}_m=L,\]
with $k_{i+1}/k_i$ (resp. $\overline{k}_{i+1}/\overline{k}_i$) purely inseparable of degree $p$.
Write $k_1=k(a)$ and $\overline{k}_1=k(b)$. We may assume $[k(a,b):k]=p^2$, else $k_1=\overline{k}_1$ and we are done by
induction. Now take any tower
\[ k(a,b)\subset k_3'\subset\ldots\subset k_m'=L\]
and by induction
\[\Tr_{k(a,b)/k_1}\circ\Tr_{k_3'/k(a,b)}\circ\ldots\circ\Tr_{k_m'/k_{m-1}'}=
                                                                    \Tr_{k_2/k_1}\circ\ldots\circ\Tr_{k_m/k_{m-1}}\]
and
\[\Tr_{k(a,b)/\overline{k}_1}\circ\Tr_{k_3'/k(a,b)}\circ\ldots\circ\Tr_{k_m'/k_{m-1}'}=
                      \Tr_{\overline{k}_2/\overline{k}_1}\circ\ldots\circ\Tr_{\overline{k}_m/\overline{k}_{m-1}}.\]
Thus from the second case it follows
\[\Tr_{k_1/k_0}\circ\ldots\circ\Tr_{k_m/k_{m-1}}=
                           \Tr_{\overline{k}_1/\overline{k}_0}\circ\ldots\circ\Tr_{\overline{k}_m/\overline{k}_{m-1}}\]
and we are done.
\end{proof}

The following proposition is the formulation of proposition \ref{0.12}, (v) for the de Rham-Witt complex, whose technical appearance
is due to the fact that we defined the trace only for field extensions.

\begin{prop}\label{2.7.5} 
Let $S$ be finite truncation set and $E, L\supset k$ finite field extensions. 
Write $E\otimes_k L=\oplus_{i=1}^n A_i$, where the $A_i$'s are local artinian rings with
residue field $(A_i)_{\text{red}}=L_i$ and length $l_i=l_{A_i}(A_i)$. We denote by $\sigma_i : E\inj L_i $ the natural inclusion.
Then the following square is commutative
\[\xymatrix{ \BWC{S}{q}{E}\ar[rr]^{\prod_i l_i\sigma_i}\ar[d]_{\Tr_{E/k}} &  & \prod_i\BWC{S}{q}{L_i}\ar[d]^{\sum_i\Tr_{L_i/L}}\\
             \BWC{S}{q}{k}\ar[rr]                                         &   & \BWC{S}{q}{L}.
}\]
\end{prop}
The following lemma will be useful in the proof of the proposition.

\begin{lem}\label{2.7.5.1}
Let $k$ be a field of characteristic $p$, $L\supset k$ a field extension.
\begin{enumerate}
\item If $E\supset k$ is a finite separable field extension, then there are finite field extensions $L_i\supset L$ such that we have an 
      isomorphism of $L$-algebras
      \[E\otimes_k L \cong \prod_i L_i.\]
\item If $E\supset k$ is purely inseparable, then $E\otimes_k L$ is local artinian.
\item Let $A$ and $B$ be two local artinian $k$-algebras such that $A\otimes_k B$ is again local artinian. Then 
$(A\otimes_k B_{\text{red}})_{\text{red}}=(A\otimes_k B)_\text{red}$ and
\[l_{A\otimes B}(A\otimes B)=l_{A\otimes B_\text{red}}(A\otimes B_\text{red})l_B(B).\] 
\end{enumerate}
\end{lem}
\begin{proof}
(i) follows from \cite[V,\S 6, No. 7, Prop. 5, Thm. 4]{Bo90} and (ii) from \cite[V, Exercises \S 14, 11)]{Bo90}.
The first statement in (iii) is a special case of \cite[I,\S 4, Cor.(4.5.12)]{EGA I}.
For the second statement we notice that $B\to A\otimes_k B$ is a flat and local homomorphism and then \cite[Lemma A.4.1.]{Fu84}
gives the assertion.
\end{proof}
\begin{proof}[Proof of Proposition \ref{2.7.5}]
We have to consider several cases.

{\em First case: $E\supset k$ is separable.} We have $E\otimes L\cong \prod_i L_i$ and $l_i=1$ (by \ref{2.7.5.1},(i)). Hence the statement follows
                                           from the definition of the trace and proposition \ref{0.12}, (vi) and (v).
                                           
{\em Second case: $E\supset k$ is purely inseparable of degree $p$.} Now $E\otimes_k L$ is local artinian (by \ref{2.7.5.1}, (ii)) and 
                                                       $(E\otimes_k L)_\text{red}=L_1$. Hence we have to show that
\[\xymatrix{\BWC{S}{q}{E}\ar[r]^{l_1\sigma_1}\ar[d]_{\Tr_{E/k}} & \BWC{S}{q}{L_1}\ar[d]^{\Tr_{L_1/L}}\\
             \BWC{S}{q}{k}\ar[r]                        & \BWC{S}{q}{L}                           }\]
commutes. Write $E=k[x]/(x^p-a)$, $a\in k\setminus k^p$. If $a\in L^p$, then $L_1=L$ and $l_1=p$. Thus the statement follows from theorem \ref{2.7}, (iii).
If $a \not\in L^p$, then $E\otimes_k L=L_1=L[x]/(x^p-a)$ and $l_1=1$ and we may argue again with \ref{2.7}, (iii).

{\em Third case: $E\supset k$ is purely inseparable of degree $p^s$.} By induction over $s$. The induction step remains to be checked. 
                            Write $ E\supset E'\supset k$ with $E\supset E'$ (resp. $E'\supset k$) purely inseparable of degree $p$ (resp. $p^{s-1}$).
                            Now  by \ref{2.7.5.1}, (ii) and (iii) 
\[l_1=l_{E\otimes_{E'}(E'\otimes_k L)_\text{red}}(E\otimes_{E'}(E'\otimes_k L)_\text{red})l_{E'\otimes_k L}(E'\otimes_k L).\] 
Then the assertion follows from the induction hypothesis and the transitivity of the trace.

{\em General case: $E\supset k$ finite field extension.} Write $E\supset E'\supset k$ with $E\supset E'$ purely inseparable and $E'\supset k$ separable
                                                         and apply the second and the third case, together with the transitivity of the trace.
\end{proof}

Since we want to prove a residue theorem, we need a residue. Now we define a residue on the de Rham-Witt complex 
over the ring of Laurent series over a $\Z_{(p)}$-algebra, generalizing the residue on K\"ahler differentials.
 
\begin{prop}\label{2.8}
Let $p$ be an odd prime, $S$ a finite truncation set and $A$ a $\Z_{(p)}$-algebra. Denote by $A((t))$ the ring
of Laurent series over $A$. Then there is a
$\BW{S}{A}$-linear map
\[\Res_S=\Res : \BWC{S}{1}{A((t))}\longrightarrow \BW{S}{A},\]
satisfying the following properties
\begin{enumerate}
\item $\gh_s\circ\Res_S=\Res\circ \Gh_s$, for all $s\in S$, where the latter $\Res$ is the usual 
         residue on $\Omega_{A((t))/A}^1$ (of course here we mean by $\Gh_s$ the composition of $\Gh_s$ with the natural map
        $\Omega^1_{A((t))}\to \Omega^1_{A((t))/A}$.)
\item $\Res\bigl(\V_n([at^j])d\V_m([bt^i])\bigr)=\begin{cases}
                                                       \sgn(i)(i,j)\V_{(mn)/c}([a]^{m/c}[b]^{n/c}) & \text{if } jm+in=0\\
                                                                           0                       & \text{else},
                                                     \end{cases}$
      where $a, b\in A$, $i,j\in\Z$, $n,m\in \N$, $c=(m,n)$ and $\sgn(i)=i/|i|$, for $i\neq 0$, $\sgn(0)=0$.
\item $\Res$ is a natural transformation in $A$.
\item If $u\in A[[t]]^\times$ and $\tau=tu$, then $A((t))=A((\tau))$ and $\Res_t=\Res_\tau$ .
\item If $\omega\in \BWC{S}{1}{A[[t]]}$ or $\in\BWC{S}{1}{\frac{1}{t}A[[\frac{1}{t}]]}$, then $\Res(\omega)=0$.
\item $\Res$ commutes with $\V_n$, $\F_n$ and restriction.
\end{enumerate}
\end{prop}
Before we start with the proof of this proposition, we need the following Lemma.
\begin{lem}\label{2.9}
In the situation of \ref{2.8}, the infinite sums  $\sum_{j=r}^\infty [a_j][t^j]$, with $r\in\Z$ and $a_j\in A$,
are well defined elements in $\BW{S}{A((t))}$ and every $w\in \BW{S}{A((t))}$ can be uniquely written
\[w=\sum_{n\in S}\V_n\bigl(\sum_{j>>-\infty}^\infty [a_{nj}][t]^j\bigr),\quad a_{nj}\in A,\]
where $j>>-\infty$ means, only finitely many negative $j$'s.
\end{lem}
\begin{proof}
It is enough to check the first statement for $A=\Z_{(p)}[x_1,\ldots]$, but if it is true for its quotient field, its true for $A$, thus 
we may assume $A$ to be a $\Q$-algebra. So we are reduced to show that $\gh_s(\sum_{j=r}^n [a_j][t^j])$ converges for $n\to \infty$ and all $s\in S$.
But this is obvious (since $\gh_s([a])=a^s$). Now take $w=(w_s)_{s\in S}\in \BW{S}{A((t))}$. Write $s_0=\min\{s\in S\,|\, w_s\neq 0\}$ and
$w_{s_0}=\sum_{j\ge r} a_jt^j\in A((t))$. Then define an operator $\varphi$ by $\varphi(w)=w-\V_{s_0}(\sum_{j\ge r}[a_j][t^j])$. 
Now the lemma follows from the observation that there is a natural number $N$ with $\varphi^N(w)=0$. 
(Indeed, we may assume $A$ is without $\Z$-torsion and then by the definition
 of $s_0$ and $\varphi$ we have $\gh_s(\varphi(w))=0$ for all $s<s_0+1$. Thus such a $N$ exists, since $S$ is finite .)  
\end{proof}
\begin{proof}[Proof of Proposition \ref{2.8}]
We define a map 
\[\Res: \BW{S}{A((t))}\times\BW{S}{A((t))}\longrightarrow \BW{S}{A}\]
by
\ml{2.9.1}{\Res\biggl(\sum_{n\in S}\V_n\bigl(\sum_{j>>-\infty}[a_{nj}][t]^j\bigr),\sum_{m\in S}\V_m\bigl(\sum_{i>>-\infty}[b_{mi}][t]^i\bigr)\biggr)=\\
           \sum_{in+jm=0}\sgn(i)(i,j) \V_{(mn)/(m,n)}([a_{nj}]^{m/(m,n)}[b_{mi}]^{n/(m,n)}).}
Notice that the sum is finite. Now we claim
\eq{2.9.2}{\gh_s(\Res(\alpha,\beta))=\Res(\Gh_s(\alpha d\beta)), \text{all }s\in S}
where $\alpha,\beta\in \BW{S}{A((t))}$, $\alpha d \beta\in\BWC{S}{1}{A((t))}$ and the $\Res$ on the right is the residue on $\Omega^1_{A((t))/A}$. 
We may assume $\alpha= \V_n([a][t]^j)$ and $\beta=\V_m([b][t]^i)$ and we write $c=(m,n)$. If $mj+ni\neq 0$, then $\Res(\alpha,\beta)=0$, 
if $mj+in=0$, then by property \ref{0.5}, (i)
\[gh_s(\Res(\alpha,\beta))=\gh_s\bigl(\sgn(i)(i,j)\V_{mn/c}([a]^{m/c}[b]^{n/c})\bigr)=\begin{cases}
                                                                                     \sgn(i)(i,j)\frac{mn}{c}(a^{s/n}b^{s/m}) & \text{if }\frac{mn}{c}|s\\
                                                                                                    0                          & \text{else}.
                                                                                     \end{cases} \]
If $\frac{mn}{c}\not|s$, then $\Gh_s(\alpha d \beta )$ is zero by definition, else 
\[\Gh_s(\alpha d \beta)= nia^{s/n} b^{s/m}t^{(js/n)+(is/m)}\tfrac{dt}{t}\in\Omega^1_{A((t))/A}.\]
Therefore
\[\Res(\Gh_s(\alpha d\beta))=\begin{cases}
                              ni a^{s/n}b^{s/m}  & \text{if } mj+ni=0\\
                                  0              & \text{else.}
                             \end{cases}\]
Now the claim follows, since $in=\sgn(i)(i,j)\frac{mn}{c}$, if $jm+in=0$. It follows that $\Res$ induces a well defined map on 
$\BWC{S}{1}{A((t))}$, denoted by $\Res$ again. Indeed, $\BWC{S}{1}{A((t))}$ is a quotient of $\BW{S}{A((t))}\times\BW{S}{A((t))}$ 
(recall that $B\otimes_\Z B/(a\otimes bc-ab\otimes c-ac\otimes b) \to \Omega^1_B$ via $a\otimes b\mapsto a db$ is an isomorphism for any ring $B$),
thus we have to show that $\Res$ vanishes on the kernel of the quotient map. Again it is enough to consider the case where $A$ has no $\Z$-torsion, 
but then the vanishing follows immediately from the claim, since $\Gh$ is well defined on $\BWC{S}{1}{A((t))}$. The other statements 
of proposition \ref{2.8} follow with the same reasoning and the corresponding properties of the usual residue on $\Omega^1_{A((t))}$.  

\end{proof}

\begin{rmk}\label{2.10}
\ref{2.8}, (i) shows that the composition $\BW{S}{A((t))}\times A((t))^\times\to \BWC{S}{1}{A((t))}\stackrel{\Res}{\to} \BW{S}{A}$, where the 
first map is given by $(w,a)\mapsto wd[a]/[a]$, is the residue symbol defined by Witt in \cite[2. Der Residuenvektor $(\alpha,\beta)$]{Wi36}.
See also \cite[4.3.4.]{AnRo04} for a definition of Witts residue symbol using the Contou-Carr\`ere symbol. 
 
\end{rmk}

\begin{rmk}\label{2.10.1}
It follows from (\ref{2.9.1}) that the map $\Res:\WC{n}{1}{A((t))}\to \W{n}{A}$ factors through the relative $p$-typical
de Rham-Witt complex $\WC{n}{1}{A((t))/A}$.
\end{rmk}
The following lemma will be essential for the residue theorem. It is a generalization of lemma 5 in chapter II of \cite{Se88}. See also
the proof of \cite[17.6. Theorem]{Ku86}.
\begin{lem}\label{2.11}
Let $p$ be an odd prime, $S$ a finite truncation set, $A$ a $\Z_{(p)}$-algebra and $B$ an $A$-algebra, which is a free $A$-module of finite rank $r$.
Take $z=t^e+\sum_{i>e}a_it^i\in A((t))$ with $e\ge 1$. Then $B((t))$ is a free
$A((z))$-module of rank $re$ and for $\beta\in \BW{S}{B((t))}$ and $\alpha\in \BW{S}{A((z))}$ we have
\eq{2.11.1}{\Tr_{B/A}(\Res_{B((t))}(\beta d\alpha))=\Res_{A((z))}(\Tr_{B((t))/A((z))}(\beta)d \alpha)\in \BW{S}{A}.} 
\end{lem}
\begin{proof}
Only the second statement must be proved. We observe that by lemma \ref{2.9} it is enough to prove (\ref{2.11.1}) 
for 
\[\beta=\V_n([bt^i]) \text{ and } \alpha=\V_m([az^j]),\]
with  $b\in B$, $a\in A$ and $i,j\in \Z$. 
By the same argument as in the proof of (ii) of proposition \ref{0.12} and using that $\Res$ is 
a natural transformation, we may assume $A$ to be $\Z$-torsion free. But then $A$ sits injective in $A\otimes_{\Z}\Q$, which extends to the 
level of Witt vectors and by base change for the trace (see proposition \ref{0.12}, (v)) we are reduced to the case, 
where $A$ and $B$ are $\Q$-algebras. Now write $z=t^eu$ with $u=1+a_{e+1}t+\ldots\in 1+tA[[t]]$. Then there is a $v\in 1+tA[[t]]$
with $v^e=u$ (namely $v=\exp(1/e\log(u))$). By proposition \ref{2.8}, (iv), we may replace $t$ by $tv$, and thus assume $t^e=z$.
Now it is time for calculations. We have 
\eq{2.11.2}{\gh_s\bigl(\Tr_{A((t))/A((z)),S}([t]^i)\bigr)=\Tr_{A((t))/A((z))}(t^{is})=\begin{cases}
                                                                               e z^{is/e} & \text{if } e|is\\
                                                                                   0      & \text{else.}
                                                                           \end{cases} }
Now write $c=(m,n)$, if $\frac{mn}{c}\not|s$, then $\gh_s \bigl(\Res_{S,z}(\Tr(\beta)d\alpha)\bigr)$ is zero, else
\[\gh_s \bigl(\Res_{S,z}(\Tr_S(\beta)d\alpha)\bigr)=
 \Res_z\Bigl(n\gh_{s/n}\bigl(\Tr_{A((t))/A((z)),S}([t]^i)\Tr_{B/A,S}([b])\bigr)a^{s/m-1}(z^j)^{s/m-1}d(az^j)\Bigr).\]
By (\ref{2.11.2}) this is zero, if $e\not|\frac{is}{n}$, else this is equal to
\[\Res_z(nej \Tr_{B/A}(b^{s/n})a^{s/m}z^{(is/(ne))+(js/m)}\tfrac{d z}{z}).\]
All in all we obtain
\eq{2.11.3}{\gh_s \bigl(\Res_{S,z}(\Tr_{B((t))/A((z)),S}(\beta)d\alpha)\bigr)=\begin{cases}
                                                               nej\Tr_{B/A}(b^{s/n})a^{s/m} & \text{if } \frac{mn}{c}|s \text{ and } im+jne=0\\
                                                                           0                & \text{else}
                                                              \end{cases} }
 and it is easy to check that this coincides with $\gh_s\Bigl(\Tr_{B/A,S}\bigl(\Res_{S,t}(\beta d\alpha)\bigr)\Bigr)$. Hence the assertion.

\end{proof}

We recall the notion of absolute and relative Frobenius, see for example \cite[\S 9]{EsVi92}.

\begin{defn}\label{2.12}
Let $k$ be a field of characteristic exponent $p$ (i.e if the characteristic of $k$ is $0$, then $p=1$ else $p$ equals the characteristic of $k$)
and $X$ a $k$-scheme of finite type. Then the \emph{absolute Frobenius} $\F_X$ on $X$ is defined to be the identity on the topological space of $X$ and
the $p$-th power on the structure sheaf, that is
\[F_X^*:\sO_X\longrightarrow \sO_X,\quad a\mapsto a^p.\]
We write $X^{(p)}=X\times_{\F_k}\Spec k$, then we have a commutative diagram
\[\xymatrix{X\ar[r]^\F\ar[dr] &  X^{(p)}\ar[r]^{\text{pr}_1}\ar[d] & X\ar[d]\\
                             &         \Spec k\ar[r]^{\F_k}           & \Spec k,}\]
where $\text{pr}_1$ is the projection map and $\F$ is the unique map, which satisfies $\F\circ\text{pr}_1=\F_X$.
 We call $\F$ the \emph{relative Frobenius} or just Frobenius.
We define $X^{(p^n)}$ via $X^{(p^n)}=(X^{p^{n-1}})^{(p)}$ and the composition of the iterated relative Frobenius
\[X\stackrel{\F}{\longrightarrow}X^{(p)}\stackrel{\F}{\longrightarrow}X^{(p^2)}\longrightarrow\ldots \longrightarrow X^{(p^n)}\]
will be denoted by
\[\F^n:X\longrightarrow X^{(p^n)}.\]
(Notice, that this notation coincides with the $n$-th Frobenius on the $p$-typical de Rham-Witt complex, but it should be clear from the context
whom of these two we are talking about.)     
\end{defn}

\begin{lem}\label{2.13}
Let $k$ be a field of characteristic exponent $p\ge 1$ and $X$ a scheme of finite type over $k$. Let $P\in X$ be a closed point and denote by
$k(P)$ the residue field of $P$. Then there is a natural number $N$ such that for all $n\ge N$,
$k(\F^n(P))=k(\F^N(P))$ is the  relative separable closure of $k$ in $k(P)$.
\end{lem}
\begin{proof}
It follows from the definition of the relative Frobenius that $k(\F^n(P))$ is the compositum of $k$ with the $p^n$-th power of $k(P)$, i.e.
$k(\F^n(P))=k\bigl(k(P)^{p^n}\bigr)$. But the chain $k(P)\supset k(\F(P))\supset k(\F^2(P))\supset\ldots\supset k$ becomes stationary, since
$[k(P):k]$ is finite and hence the assertion.
\end{proof}

\begin{lem}\label{3.3}
Let $C$ be a smooth projective curve with function field K. For $j\ge 1$ let  $\F^j:C\to C^{(p^j)}$ be the relative Frobenius and $K_j$ the function field 
of $C^{(p^j)}$. Then for all closed points $P\in C$
\[ p^j=[K:K_j]=e_Pf_P,\]
where $e_p$ is the ramification index, i.e. $e_P=\va_P(z)$ with $z\in\sO_{C^{(p^j)},\F^j(P)}$ a local parameter in $\F^j(P)$, and $f_P=[k(P):k(\F^j(P))]$.
\end{lem}

\begin{proof}
The first equation, $p^j=[K:K_j]$, is \cite[Proposition 2.11, (c)]{Si92} and $[K:K_j]=e_Pf_P$ is a special case of \cite[I,\S 4, Proposition 10]{Se68}.
\end{proof}

\begin{defn-prop}[cf. \cite{Ku86}, 17.4.]\label{2.14}
Let $S$ be a finite truncation set, $k$ a field with characteristic exponent $p\ge 1$, $C$ a smooth projective curve over $k$ with function field
$K=k(C)$ and $P\in C$ a closed point. Let $n$ be a natural number such that $k(\F^n(P))$ is the relative separable  closure of $k$ in $k(P)$ (such an
$n$ exists by the lemma). We write $P_n:=\F^n(P)\in C^{(p^n)}$, $\kappa_n=k(P_n)$ and $K_n=k(C^{(p^n)})=k(K^{p^n})$. Finally we denote by
$\widehat{K}_n$ the completion of $K_n$ in $P_n$. Now the choice of a local parameter $t$ in $P_n$ identifies $\widehat{K}_n\cong\kappa_n((t))$ (this
is an isomorphism over $k$, since $\kappa_n\supset k$ is separable) and we have a natural inclusion $\iota : K_n\inj \kappa_n((t))$. 
Take $\omega\in \BWC{S}{1}{K}$, then we define the \emph{residue of $\omega$ in $P$} to be
\eq{2.14.1}{\Res_{S,P}(\omega)=\Res_P(\omega)=\Tr_{\kappa_n/k}\Bigl(\Res_{t,S}\bigl(\iota(\Tr_{K/K_n}(\omega))\bigr)\Bigr)\in \BW{S}{k}, }
where the $\Res_{t,S}$ on the right hand side, is the residue on $\BWC{S}{1}{\kappa_n((t))}$ from proposition \ref{2.8}.
The residue is well defined, i.e. independent from the choice of the local parameter $t$ and the number $n$.
\end{defn-prop}
\begin{proof}
The independence of the choice of the local parameter follows from proposition \ref{2.8}, (iv). Thus it remains to check the independence of the choice
of $n$. Therefore it is enough to show that nothing changes, if we replace $n$ by $n+1$. By the choice of $n$, 
we have $\kappa_n=\kappa_{n+1}=:\kappa$ and 
$K_{n+1}=k(K_n^p)\subset K_n$ is purely inseparable of degree $p$ (lemma \ref{3.3}). If $t$ is a local parameter in $P_n$, $z=t^p$ is one in $P_{n+1}$
(by lemma \ref{3.3}).
Thus our situation is
\[\xymatrix{K_n\ar@{^{(}->}[r]^{\iota_n} & \kappa((t))\\
            K_{n+1}\ar@{^{(}->}[u]\ar@{^{(}->}[r]^{\iota_{n+1}} & \kappa((z))\ar@{^{(}->}[u] }\]
and we want to show
\eq{2.14.1.5}{\Res_t\bigl(\iota_n(\Tr_{K/K_n}(\omega))\bigr)=\Res_z\bigl(\iota_{n+1}(\Tr_{K/K_{n+1}}(\omega))\bigr).}
It follows easily from theorem \ref{2.7}, (iii), that in $\BWC{S}{1}{K_n}$
\[\iota_{n+1}\circ\Tr_{K_n/K_{n+1}}= \Tr_{\kappa((t))/\kappa((z))}\circ\iota_n.\]  
Together with $\Tr_{K/K_{n+1}}= \Tr_{K_n/K_{n+1}}\circ\Tr_{K/K_n}$, we see that it is enough to prove
\eq{2.14.2}{\Res_t(\alpha)=\Res_z\bigl(\Tr_{\kappa((t))/\kappa((z))}(\alpha)\bigr), \text{ for  }\alpha \in \BWC{S}{1}{\kappa((t))}.}
Furthermore, we may assume by lemma \ref{2.9}
\[\alpha=\V_m([at^i])d\V_n([bt^j]),\]
with $n, m\in S$, $i,j\in\Z$ and $a,b\in \kappa$. Then we consider the following cases.

\emph{First case:} $n=m=1$. Thus 
\[\alpha=j[abt^{i+j-1}]d[t]+[at^{i+j}]d[b],\]
 we denote the first summand by $\alpha'$. If we write
$i+j-1=pl+r$ with $l\in\Z$ and $0\le r\le p-1$, we have
\[\alpha'=j[abz^l][t^r]d[t]\]
and calculating the trace as in (\ref{2.7.2}) we obtain
\[\Tr(\alpha')=\begin{cases}
                 \tfrac{j}{r+1}[abz^l]d\V_p([z]^{r+1}) & \text{if } r<p-1\\
                  j[abz^l]d[z]                          & \text{ if } r=p-1.
               \end{cases}\]
Thus by proposition \ref{2.8}, (ii)
\begin{eqnarray*}
\Res_z(\Tr(\alpha)) & = & \Res_z(\Tr(\alpha'))=\begin{cases}
                                            \tfrac{j}{r+1} (l,r+1)p[ab] & \text {if } pl+r+1=0, r<p-1\\
                                              j[ab]                      & \text{if } l+1=0, r=p-1\\
                                              0                          & else
                                           \end{cases}\\
              &       = &  \begin{cases}
                           j[ab] &\text{if } i+j=0\\
                            0    & \text{else}
                            \end{cases}\\                 
              &       = & \Res_t(\alpha).
\end{eqnarray*}

\emph{Second case: $n=1$, $m$ arbitrary.} Then $\alpha=\V_m(j[ab^m][t^{i+m(j-1)}])d[t]+\V_m([a] [t^{i+mj}])d[b]$. Again we denote the first summand by
$\alpha'$, thus
\[\alpha'=\V_m(j[ab^m][t^{i+mj-1}]d[t]).\]
Now both, $\Res_t$ and $\Res_z$, only see $\alpha'$, but since they commute with $\V_m$ we  are done by the first case.

\emph{Third case: $p\not|n$ and $m$ arbitrary. } Write $c=(m,n)$. Then
\[\alpha=\V_m([a][t^i])\tfrac{1}{n}\V_n(1)d\V_n([b][t^j])=\V_n\bigl(\V_{m/c}(\tfrac{c}{n}[a]^{n/c}[t^{in/c}])d([b][t^j])\bigr)\]
and we are done by the second case.

\emph{Fourth case: $n$ arbitrary and $p\not|m$.} Here we write
\[\alpha=d\bigl(\V_m([a][t^i])\V_n([b][t^j])\bigr)-\V_n([b][t^j])d\V_m([a][t^i])\]   
and since $\Tr$ commutes with $d$ and $\Res\circ d=0$, the result follows from the third case.

\emph{Fifth case: $m,n$ arbitrary.} Write $c=(m,n)$, $r=v_p(c)$ and $m'=m/p^r$, $n'=n/p^r$, then
\[\alpha=\V_{p^r}\bigl(\V_{m'}([a][t^i])d\V_{n'}([b][t^j])\bigr)\]
and the assertion follows either from the third or the fourth case. 
\end{proof}

\begin{rmk}\label{2.15}
In the situation of \ref{2.14} we have
\[\Res_P(\omega)=\Res_{P_n}(\Tr_{K/K_n}(\omega)), \text{ for all }n\in \N.\]
(This holds by definition for $n$ sufficiently large and  follows else from $\Tr_{K/K_m}=\Tr_{K_n/K_m}\circ\Tr_{K/K_n}$.) 

\end{rmk}

\begin{rmk}\label{2.15.5}
Still in the situation of \ref{2.14}. If $\omega\in \BWC{S}{1}{K}$ has no pole in $P\in C$ (i.e. $\omega$ is in the image of the 
natural map $\BWC{S}{1}{\sO_{C,P}}\to\BWC{S}{1}{K}$), then $\Res_P(\omega)=0$.
(This follows from \ref{2.8}, (v) and the fact that the trace map $\Tr_{K/K_n}$ restricts by construction to a map 
$\Tr:\BWC{S}{1}{\sO_{C,P}}\to\BWC{S}{1}{\sO_{C^{(p^n)},P_n}}$.)
\end{rmk}

The following proposition is a generalization of the trace formula for residues on K\"ahler differentials, see \cite[17.6.]{Ku86} or 
\cite[II, Lemma 4.]{Se88}.

\begin{prop}\label{2.16}
Let $S$ be a finite truncation set and $f: C\longrightarrow C'$ be a separable dominant morphism  of smooth projective curves over a field $k$. Denote
by $K$ (resp. $K'$) the function field of $C$ (resp. $C'$). Let $P\in C'$ be a closed point and $\omega\in \BWC{S}{1}{K}$.
Then
\[\sum_{Q\mapsto P}\Res_Q(\omega)=\Res_P(\Tr_{K/K'}(\omega)).\]
\end{prop}
\begin{proof}
First we assume that the residue fields of $P$ and the $Q$'s mapping to $P$ are separable over the ground field $k$. We denote by
$\widehat{K}_Q$ (resp. $\widehat{K'}_P$) the completion of $K$ (resp. $K'$) in $Q$ (resp. $P$). Since $K$ is separable over $K'$,
we see from \cite{Bo89}, chapter VI, $\S 8$, No.2, corollary 2, that we have an isomorphism
\[\widehat{K'}_P\otimes_{K'}K\stackrel{\simeq}{\longrightarrow} \prod_{Q\mapsto P}\widehat{K}_Q, \quad e\otimes f\mapsto (ef,ef,\ldots).\]
Thus
\[\BW{S}{\widehat{K'}_P\otimes_{K'}K}\cong\prod_{Q\mapsto P}\BW{S}{\widehat{K}_Q}.\]
Since $\BWC{S}{1}{K}=\BW{S}{K}\otimes_{\BW{S}{K'}}\BWC{S}{1}{K'}$, we may write $\omega=xdy$ with $x\in\BW{S}{K}$ and $y\in \BW{S}{K'}$.
Identifying the elements of $\BW{S}{K}$ (resp. $\BW{S}{K'}$) with their images in $\BW{S}{\widehat{K'}_P\otimes_{K'}K}$ (resp. $\BW{S}{\widehat{K'}_P}$)
and writing $x_Q$ for the image of $x$ in $\BW{S}{\widehat{K}_Q}$, we obtain in $\BW{S}{\widehat{K'}_P}$
\[\Tr_{K/K'}(x)=\Tr_{\widehat{K'}_P\otimes_{K'}K/\widehat{K'}_P}(x)=\sum_{Q\mapsto P}\Tr_{\widehat{K}_Q/\widehat{K'}_P}(x_Q).\]
Thus we have to show
\eq{2.16.1}{\Res_Q(x_Qdy_P)=\Res_P\bigl(\Tr_{\widehat{K}_Q/\widehat{K'}_P}(x_Q)dy_P\bigr).}
Now let $t$ be a local parameter in $Q$ and $z$ a parameter in $P$, since we assumed $k(Q)$ and $k(P)$ to be separable over k, we may
identify
\[\widehat{K'}_P=k(P)((z)) \text{ and } \widehat{K}_Q=k(Q)((t)).\]
Furthermore, we  may write $z=t^e+\sum_{i>e}a_it^i$, with $e=v_Q(z)$. But then (\ref{2.16.1}) follows from lemma \ref{2.11} 
and definition (\ref{2.14.1}).

Now in general we choose $n$ sufficiently large such that $k(\F^n(Q))$ and $k(\F^n(P))$ are separable over $k$, for all $Q\mapsto P$. Then
applying the case above to $C^{(p^n)}\to {C'}^{(p^n)}$ and by remark \ref{2.15}, we obtain (with the notation of \ref{2.14}) 
\begin{eqnarray*} \sum_{Q\mapsto P}\Res_Q(\omega) & = & \sum_{Q\mapsto P}\Res_{Q_n}\bigl(\Tr_{K/K_n}(\omega)\bigr)
                                                        =\Res_{P_n}\bigl(\Tr_{K_n/{K'}_n}\Tr_{K/K_n}(\omega)\bigr)\\
                                                   & =  & \Res_{P_n}\bigl(\Tr_{K'/{K'}_n}\Tr_{K/K'}(\omega)\bigr)=\Res_P\Tr_{K/K'}(\omega)
\end{eqnarray*}
and we are done.
\end{proof}

Now we can prove the residue theorem for the de Rham-Witt complex in degree one. This generalizes the classical residue theorem for 
K\"ahler differentials on curves. Our proof is the classical one (with some adjustments), see \cite[II, Proposition 6]{Se88} or 
\cite[17.7.]{Ku86}. Witt also stated such a theorem for his residue symbol on $\BW{S}{K}\times K^\times$, where $K$ is the function field
of a curve, see \cite[9. Analogon zum Residuensatz]{Wi36}. Anderson and Romo give in \cite[4.3.6.]{AnRo04} a different
proof of this,
by first defining the residue symbol via the Contou-Carr\`ere symbol and then using a general reciprocity law for the Contou-Carr\`ere symbol,
which they prove by methods, similar to the one of Tate in \cite{Ta68}.

\begin{thm}\label{2.17}
Let $S$ be a finite truncation set, $k$ a field and $C$ a smooth projective curve over $k$ with function field $K$. Then
\[\sum_{P\in C}\Res_P(\omega)=0, \text{ for all } \omega\in \BWC{S}{1}{K}. \]
(Notice, that $\Res_P(\omega)=0$, if $\omega\in \BWC{S}{1}{\sO_{C,P}}$, thus the sum is finite.) 
\end{thm}
\begin{proof}
Since $C$ is smooth over $k$, there exists an element $x\in K$ such that $k(x)\subset K$ is separable and finite. Thus $x$ defines a separable dominant
morphism $C\to \P^1_k$ and by proposition \ref{2.16} we are reduced to $C=\P^1_k$. Furthermore, we may assume that the points $P\in \P^1_k$ with
$\Res_P(\omega)\neq 0$ are \'etale over $\Spec k$, since by remark \ref{2.15} and with the notations of \ref{2.14}
\[\sum_{P\in C}\Res_P(\omega)=\sum_{P_n\in C^{(p^n)}}\Res_{P_n}(\Tr_{K/K_n}(\omega)), \text{ for all } n\in\N.\]
Denote by $\bar{k}$ the algebraic closure of $k$. Then the pull-back of a closed point $P$ which is \'etale over $k$ under $\P^1_{\bar{k}}\to \P^1_k$
decomposes into finitely many $\bar{k}$-rational points without multiplicities, i.e. 
$k(P)\otimes_k \bar{k}=\prod_i \bar{k}$. Thus we obtain a commutative diagram
\[\xymatrix{\BWC{S}{1}{k(P)((t))}\ar[r]^{\Res_t}\ar[d] & \BW{S}{k(P)}\ar[r]^{\Tr_{k(P)/k}}\ar[d] & \BW{S}{k}\ar@{^{(}->}[d]\\
            \prod_i\BWC{S}{1}{\bar{k}((t))}\ar[r]^{\prod_i \Res_t}  &  \prod_i\BW{S}{\bar{k}}\ar[r]^\sum & \BW{S}{\bar{k}}.}\]
Therefore, if $\omega'$ is the image of $\omega $ under $\BWC{S}{1}{k(x)}\to \BWC{S}{1}{\bar{k}(x)}$ and $\sum_i Q_i$ is the pull-back of $P$
under $\P^1_{\bar{k}}\to\P^1_k$, then
\[\Res_P(\omega)=\sum_i\Res_{Q_i}(\omega') \text{ in } \BW{S}{\bar{k}},\]
hence we may assume $k$ to be algebraically closed. Now we need a small lemma.

\begin{lem}\label{2.18}
Let k be an algebraically closed field, then any $w\in \BW{S}{k(x)}$ may be uniquely written in the following way
\[ w=\sum_{n\in S}\V_n\Bigl(\sum_{j\ge 0}[a_{jn}][x^j]+\sum_i[b_{in}][\tfrac{1}{(x-c_{in})^{r_{in}}}]\Bigr),\]
where all sums are finite and $a_{in}, b_{in},c_{in}\in k$ and $r_{in}\ge 1$.
\end{lem}
\begin{proof}[Proof of the lemma]
It is known that the elements
\[x^j, \text{ for } j\ge 0, \text{ and } \tfrac{1}{(x-a)^r}, \text{ for } a\in k, r\ge 1\]
form a $k$-basis of $k(x)$. Now the rest of the proof works exactly as in lemma \ref{2.9}.
\end{proof}

It follows from the lemma and proposition \ref{1.15}, that we have to prove the theorem only for the following two types of differentials,
\[\omega_1=\V_n([a][x]^j)d\V_m([b][x]^i), \text{ with } a,b\in k, i,j\ge 0\]
and
\[\omega_2=\V_n(\tfrac{[a]}{[x-c]^r})d\V_m([b][x]^i), \text{ with } a,b,c\in k, i\ge 0, r\ge 1.\]
Obviously, $\Res_P(\omega_1)=0$ for all $P\in \P^1_k\setminus\{\infty\}$ and in $P=\infty$ the element 
$t=1/x$ is a local parameter, thus by proposition \ref{2.8}, (ii)
\[\Res_\infty(\omega_1)=\Res_t\bigl(\V_n([a][t]^{-j})d\V_m([b][t]^{-i})\bigr)=0.\]
Now $\omega_2$.  Denote by $p$ the characteristic of $k$ and write $A=\Z_{(p)}[z_a,z_b,z_c]$. We have a map $\varphi:A\to k$ sending
$z_a$ (resp. $z_b$, $z_c$) to $a$ (resp. $b$, $c$). The only points with non vanishing residue are $P=c$ and $\infty$.
We have $t=x-c$ is a local parameter in $P=c$ and
\[\widetilde{\omega}_{2,P}=\V_n(\tfrac{[z_a]}{[t^r]})d\V_m([z_b][t+z_c]^i)\in \BWC{S}{1}{A((t))}\]
is a lifting of $\omega_{2,P}\in\BWC{S}{1}{k((t))}$, in particular $\Res_P(\omega_2)=\varphi(\Res_t(\widetilde{\omega}_{2,P}))$.
It follows from the definition of $\Res$ that $\gh_s(\Res(\widetilde{\omega}_{2,P}))=0$ if $m\not|s$ or $n\not|s$, else
\ml{2.17.1}{\gh_s(\Res(\widetilde{\omega}_{2,P}))=\Res(niz_a^{s/n}z_b^{s/m} t^{-sr/n}(t+z_c)^{si/m-1}dt)\\
=\begin{cases}
       ni\binom{\frac{is}{m}-1}{\frac{rs}{n}-1}z_a^{s/n}z_b^{s/m}z_c^{s(in-rm)/(nm)} & \text{if } n,m|s, \tfrac{rs}{n}\le \tfrac{is}{m}\\
                                   0                                                 &  \text{else}
 \end{cases}  }
In $P=\infty$ we have that $t=1/x$ is a local parameter and 
\[\widetilde{\omega}_{2,\infty}=\V_n(\tfrac{[z_a][t]^r}{[1-z_ct]^r})d\V_m([z_b][t]^{-i})\in\BWC{S}{1}{A((t))}\]
is a lifting of $\omega_{2,\infty}\in \BWC{S}{1}{k((t))}$, in particular $\Res_\infty(\omega_2)=\varphi(\Res_t(\widetilde{\omega}_{2,\infty}))$.
Again, $\gh_s(\Res(\widetilde{\omega}_{2,\infty}))=0$, if $m\not|s$ or $n\not|s$ and else
\[\gh_s(\Res(\widetilde{\omega}_{2,\infty}))=\Res(-inz_a^{s/n}z_b^{s/m}t^{s(rm-in)/(mn)}(1+z_ct+z_c^2t^2+\ldots)^{rs/n}\tfrac{dt}{t}).\]
But the coefficient of $t^j$ in $(1+z_ct+z_c^2t^2+\ldots)^{rs/n}$ is 
\[z_c^j\#\{I\in\N_0^{rs/n}\text{ with } |I|=j\}=z_c^j\binom{\frac{rs}{n}+j-1}{\frac{rs}{n}-1}.\]
Thus we obtain 
\eq{2.17.2}{\gh_s(\Res(\widetilde{\omega}_{2,\infty}))=\begin{cases}
                                 -in\binom{\frac{rs}{n}+\frac{is}{m}-\frac{rs}{n}-1}{\frac{rs}{n}-1} z_a^{s/n}z_b^{s/m}z_c^{s(in-rm)/(nm)} & 
                                                                    \text{if } m,n|s, \tfrac{rs}{n}\le \tfrac{is}{m}\\
                                                           0                                                & \text{else}.
                                                 \end{cases}}
All together we see
\[\gh_s(\Res(\widetilde{\omega}_{2,c}))+\gh_s(\Res(\widetilde{\omega}_{2,\infty}))=0,\text{ for all }s\in S\]
and since $A$ has no $\Z$-torsion, we obtain
\[0=\varphi(\Res(\widetilde{\omega}_{2,c})+\Res(\widetilde{\omega}_{2,\infty}))=\sum_{P\in \P^1_k}\Res_P(\omega_2)\]
and we are done.

\end{proof}
Next we want to define a residue in higher degrees, cf. \cite{BlEs03}, (6.14).

\begin{rmk}\label{2.19}
We recall the definition of the residue on $ \Omega^q_{A((t))/\bb{Z}}=\Omega^q_{A((t))}$, for any ring $A$. Take
$\omega\in\Omega^q_{A((t))}$. Then we can write $\omega$ uniquely in the following way
\[\omega=\sum_{i=1}^m\beta_i\tfrac{dt}{t^i}+\omega_0, \]
with $\beta_i\in\Omega^{q-1}_{A}$, 
$\omega_0\in\bigl(\Omega^q_{A[[t]]}+A((t))\otimes_A\Omega^q_{A}\bigr)$ and $m\in\N$. Then the 
residue in $\omega$ is defined to be
\[\Res^q(\omega)=\beta_1.\]
We obtain a map $\Res^q:\Omega^q_{A((t))}\to \Omega^{q-1}_{A}$, which satisfies
\begin{enumerate}
\item $\Res^q$ is $\Omega^{\kaydot}_{A}$-linear, i.e. $\Res^q(\alpha\omega)=\alpha\Res^{q-i}(\omega)$, for
        $\alpha\in\Omega^i_A$ and $\omega\in\Omega^{q-i}_{A((t))}$.
\item $\Res^q$ is a natural transformation in $A$.
\item $\Res^q\circ d=d\circ\Res^{q-1}$.
\item If $u\in (A[[t]])^\times$ and $\pi=tu$, then $A((t))=A((\pi))$ and $\Res^q_t=\Res^q_\pi$.
\item $\Res^q\omega=0$, for $\omega\in\Omega^q_{A[[t]]}$ or  $\omega\in\Omega^q_{\frac{1}{t}A[[\frac{1}{t}]]}$.
\item Take $\beta,\gamma\in A[[t]]$ and $\alpha_i=t^{m_i}u_i$, with $u_i\in (A[[t]])^\times$, $m_i\in\bb{Z}$, for $i=1,\ldots, q$. Then
        \[\Res^q(\beta\tfrac{d\alpha_1}{\alpha_1}\ldots\tfrac{d\alpha_{q-1}}{\alpha_{q-1}}\tfrac{dt}{t})=
          \beta(0)\tfrac{du_1(0)}{u_1(0)}\ldots\tfrac{du_{q-1}(0)}{u_{q-1}(0)}\]
      and
      \[\Res^q(\beta\tfrac{d(1+t\gamma)}{(1+t\gamma)}\tfrac{d\alpha_2}{\alpha_2}\ldots\tfrac{d\alpha_q}{\alpha_q})=0.\]
\end{enumerate}
(By linearity the proof of (iv) reduces to the known case $q=1$, the rest is easy.)
\end{rmk}

For the construction of a residue on the de Rham-Witt complex in higher degrees with respect to an arbitrary $\Z_{(p)}$-algebra  $A$, we need that
the ghost map $\Gh$ is injective, if $A$ has no $p$-torsion. But the injectivity only holds in the relative situation, $\WC{n}{q}{A/\Z_{(p)}}$
(cf. \ref{1.20}-\ref{1.23}). Unfortunately we have no  relative version of the generalized de Rham-Witt complex. Thus we first
construct the residue in the $p$-typical situation for arbitrary $\Z_{(p)}$-algebras ($p$ odd), using the relative de Rham-Witt complex of
Langer and Zink and as a corollary we will get a residue in the absolute situation for $A$ being any $k$-algebra, $k$ a field.  

\begin{prop}\label{2.19.5}
Let $p$ be an odd prime, $q,n\in\N$ and $A$ a $\Z_{(p)}$-algebra. Then there is a map
\[\Res^q_{n,t}=\Res^q:\WC{n}{q}{A((t))}\longrightarrow \WC{n}{q-1}{A/\Z_{(p)}},\]
which we will call the ($p$-typical) residue map, satisfying the following properties
\begin{enumerate}
\item For all $i\in \N_0$,
        \[\Gh_{p^i}\circ\Res^q_n=\Res^q\circ \Gh_{p^i},\]
     where the latter $\Res^q$ is the one from remark \ref{2.19}.
\item $\Res^q$ is $\WC{n}{\kaydot}{A}$-linear, i.e. $\Res^q(\alpha\omega)=\bar{\alpha}\Res^{q-i}(\omega)$, for
        $\alpha\in\WC{n}{i}{A}$, $\bar{\alpha}$ its image in $\WC{n}{i}{A/\Z_{(p)}}$ and $\omega\in\WC{n}{q-i}{A((t))}$.
\item $\Res^q$ is a natural transformation in $A$.
\item $\Res^q\circ d=d\circ\Res^{q-1}$ and $\Res^q$ commutes with $\V^i$, $F^i$ and restriction.
\item If $u\in (A[[t]])^\times$ and $\pi=tu$, then $A((t))=A((\pi))$ and $\Res^q_t=\Res^q_\pi$.
\item $\Res^q\omega=0$, for $\omega\in\WC{n}{q}{A[[t]]}$ or  $\omega\in\WC{n}{q}{\frac{1}{t}A[[\frac{1}{t}]]}$.
\item Take $w\in \W{n}{A[[t]]}$, $\gamma\in A[[t]]$ and $\alpha_i=t^{m_i}u_i$, with $u_i\in (A[[t]])^\times$, $m_i\in\bb{Z}$, for $i=1,\ldots, q$. 
      Then
        \[\Res^q(w\tfrac{d[\alpha_1]}{[\alpha_1]}\ldots\tfrac{d[\alpha_{q-1}]}{[\alpha_{q-1}]}\tfrac{d[t]}{[t]})=
          w(0)\tfrac{d[u_1(0)]}{[u_1(0)]}\ldots\tfrac{d[u_{q-1}(0)]}{[u_{q-1}(0)]}\]
      and
      \[\Res^q(w\tfrac{d[1+t\gamma]}{[1+t\gamma]}\tfrac{d[\alpha_2]}{[\alpha_2]}\ldots\tfrac{d[\alpha_q]}{[\alpha_q]})=0.\]
\end{enumerate}
\end{prop}
\begin{proof}
First assume $A$ is a $\Q$-algebra. Then $\Gh$ is an isomorphism and we define
\[\Res^q_n=\Gh^{-1}\circ\prod_{i=0}^n(\Res^q\circ\Gh_{p^i}).\]
It follows immediately from remark $\ref{2.19}$ and the properties $\ref{1.19}$ that $\Res^q_n$ fulfills (i)-(vii).
Now assume $A=\Z_{(p)}[x_1,\ldots,x_d]$ or $A=\Z_{(p)}[x_1,\ldots,x_d,\frac{1}{x_1},\ldots,\frac{1}{x_d}]$ and write $B=\Q\otimes_{\Z_{(p)}}A$.
By corollary \ref{1.22} we have an inclusion
\eq{2.20.1}{\WC{n}{q}{A/\Z_{(p)}}\inj\WC{n}{q}{B}=\WC{n}{q}{B/\Z_{(p)}}.}
We claim that the composition
\[\WC{n}{q}{A((t))}\stackrel{\iota}{\longrightarrow} 
                                                \WC{n}{q}{B((t))}\stackrel{\Res^q_B}{\longrightarrow}\WC{n}{q-1}{B}\]
factors to give a map
\[\Res^q_A:\WC{n}{q}{A((t))}\longrightarrow\WC{n}{q-1}{A/\Z_{(p)}},\]
which then automatically satisfies (i), (ii) and (iv)-(vii). By the injectivity of (\ref{2.20.1}) it is enough to show
that for all $\omega\in\WC{n}{q}{A((t))}$
\eq{2.20.2}{\Res^q_B(\iota(\omega))\in\im(\WC{n}{q-1}{A/\Z_{(p)}}\inj\WC{n}{q-1}{B}).}
By lemma \ref{2.9} we may assume
\eq{2.20.2.5}{ \omega=\V^{n_0}\bigr(\sum_{j_0>>-\infty}^\infty[a_{n_0j_0}t^{j_0}]\bigl)
                    d\V^{n_1}\bigr(\sum_{j_1>>-\infty}^\infty[a_{n_1j_1}t^{j_1}]\bigl)\ldots
                    d\V^{n_q}\bigr(\sum_{j_q>>-\infty}^\infty[a_{n_qj_q}t^{j_q}]\bigl),\quad a_{ij}\in A}
and since the calculation of $\Res^q$ on $\Omega^q_{A((t))}$ involves only finitely many terms, we may even assume
\[\omega=\V^{n_0}([a_{n_0j_0}t^{j_0}])d\V^{n_1}([a_{n_1j_1}t^{j_1}]\bigl)\ldots d\V^{n_q}([a_{n_qj_q}t^{j_q}]),\]
with $a_{lk}\in A$ and $j_i\in\Z$. That is, $\omega$ is in the image of the map
\[\W{n}{A(t)}\otimes_{\W{n}{A[t]}}\WC{n}{q}{A[t]}=\WC{n}{q}{A(t)}\longrightarrow \WC{n}{q}{A((t))}\]
and therefore $\omega$ is a sum of products of the elements (\ref{2.1.1})-(\ref{2.1.4}) with elements from $\W{n}{A(t)}$.
But for $w\in \W{n}{A(t)}$
\eq{2.20.3}{ \Res^q_B(wa[t]^j) =  0, \text{ for } a\in\WC{n}{q}{A},j\in \N_0,}
\eq{2.20.4}{ \Res^q_B(wb[t]^{j-1}d[t])  =  b\Res^1_B(w[t]^{j-1}d[t]), \text{ for } b\in\WC{n}{q-1}{A}, j\in\N,}
\eq{2.20.5}{ \Res^q_B(wV^s(a[t]^j))  =  0, \text{ for } a\in \WC{n-s}{q}{A}, j\in I_p,}
\ml{2.20.6}{\Res^q_B(wd\V^s(b[t]^j))  =   (-1)^{q-1}d\V^s(b_1)\ldots d\V^s(b_{q-1})\Res^1_B(wd\V^s(b_0[t]^j)),\\
                                        \text{ for } b=b_0db_1\ldots db_{q-1}\in\WC{n-s}{q-1}{A}, j\in I_p,} 
and since the claim is known for $\Res^1_B$, it follows for all $q$.

Now let  $A$ be any ring and take $\omega\in \WC{n}{q}{A((t))}$. Write $\omega$ as in (\ref{2.20.2.5}) and define
$\omega_m$ ($m\in\N$) by
\eq{2.20.6.5}{\omega=\V^{n_0}\bigr(\sum_{j_0>>-\infty}^m[a_{n_0j_0}t^{j_0}]\bigl)
                    d\V^{n_1}\bigr(\sum_{j_1>>-\infty}^m[a_{n_1j_1}t^{j_1}]\bigl)\ldots
                    d\V^{n_q}\bigr(\sum_{j_q>>-\infty}^m[a_{n_qj_q}t^{j_q}]\bigl).} 
Now for each $m$ there is a polynomial ring $\tilde{A}_m=\Z_{(p)}[x_1,\ldots,x_d]$ with a ring homomorphism $\varphi:\tilde{A}_m\to A$ such that
$\omega_m$ has a lifting $\tilde{\omega}_m\in\WC{n}{q}{\tilde{A}_m}$. Then 
$m_1=\min\{m|\varphi\bigl(\Res^q_{\tilde{A}_m}(\tilde{\omega}_m)\bigr)=\varphi\bigl(\Res^q_{\tilde{A}_{m+r}}(\tilde{\omega}_{m+r})\bigr), 
                                                                                                                                   \text{ all }r\in\N\}$ 
exists (by the construction of $\Res^q$ above) and we define
\[\Res^q_{n,A}(\omega):=\varphi(\Res^q_{n,\tilde{A}_{m_1}} ( \tilde{\omega}_{m_1}))\in\WC{n}{q-1}{A/\Z_{(p)}}.\]
It follows from (\ref{2.20.3})-(\ref{2.20.6}) and the naturalness of $\Res^1$ that this map is independent of the lifting.
Obviously $\Res^q_A$ fulfills (i)-(iv) and (vi). To prove (v) and (vii) we must be able to lift units $u\in A[[t]]^\times$. But by definition of
$\Res^q_A$ the statements depend only on a finite part of $u$. Thus we may assume $u=\sum_{i=0}^m a_it^i$, with $a_i\in A, a_0\in A^\times$.
If we now consider the map $\tilde{A}=\Z_{(p)}[x_0,\frac{1}{x_0},x_1,\ldots, x_m]\to A$, $x_i\mapsto a_i$, then 
$\tilde{u}=\sum_{i=0}^m x_it^i\in (\tilde{A}[[t]])^\times$ lifts $u$. Being able to lift $u$ to $(\tilde{A}[[t]])^\times$ the statements (v) and (vii)
follow from the first case and we are done.
\end{proof}

\begin{cor}\label{2.20}
Let $k$ be a field of characteristic $p\neq 2$, $S$ a finite truncation set, $q\in \N$ and $A$ a $k$-algebra. Then there is a map
\[\Res^q_{S,t}=\Res^q:\BWC{S}{q}{A((t))}\longrightarrow \BWC{S}{q-1}{A},\]
which we will call the residue map, satisfying the following properties
\begin{enumerate}
\item For all $s\in S$,
        \[\Gh_s\circ\Res^q_S=\Res^q\circ \Gh_s,\]
     where the latter $\Res^q$ is the one from remark \ref{2.19}.
\item $\Res^q$ is $\BWC{S}{\kaydot}{A}$-linear, i.e. $\Res^q(\alpha\omega)=\alpha\Res^{q-i}(\omega)$, for
        $\alpha\in\BWC{S}{i}{A}$ and $\omega\in\BWC{S}{q-i}{A((t))}$.
\item $\Res^q$ is a natural transformation in $A$.
\item $\Res^q\circ d=d\circ\Res^{q-1}$ and $\Res^q$ commutes with $\V_n$, $F_n$ and restriction.
\item If $u\in (A[[t]])^\times$ and $\pi=tu$, then $A((t))=A((\pi))$ and $\Res^q_t=\Res^q_\pi$.
\item $\Res^q\omega=0$, for $\omega\in\BWC{S}{n}{A[[t]]}$ or  $\omega\in\BWC{S}{q}{\frac{1}{t}A[[\frac{1}{t}]]}$.
\item Take $w\in \BW{S}{A[[t]]}$, $\gamma\in A[[t]]$ and $\alpha_i=t^{m_i}u_i$, with $u_i\in (A[[t]])^\times$, $m_i\in\bb{Z}$, for $i=1,\ldots, q$. 
      Then
        \[\Res^q(w\tfrac{d[\alpha_1]}{[\alpha_1]}\ldots\tfrac{d[\alpha_{q-1}]}{[\alpha_{q-1}]}\tfrac{d[t]}{[t]})=
          w(0)\tfrac{d[u_1(0)]}{[u_1(0)]}\ldots\tfrac{d[u_{q-1}(0)]}{[u_{q-1}(0)]}\]
      and
      \[\Res^q(w\frac{d[1+t\gamma]}{[1+t\gamma]}\tfrac{d[\alpha_2]}{[\alpha_2]}\ldots\tfrac{d[\alpha_q]}{[\alpha_q]})=0.\]
\end{enumerate}
\end{cor}

\begin{proof}
We define $\Res^q_{S,t}$ to be the composition
\begin{multline}\BWC{S}{q}{A((t))}\stackrel{\simeq(\ref{1.13})}{\longrightarrow}\prod_{j\in I_p}\BWC{\sP\cap S/j}{q}{A((t))}(\frac{1}{j})
\stackrel{\ref{2.19.5}}{\longrightarrow}\prod_{j\in I_p}\BWC{\sP\cap S/j}{q-1}{A/\Z_{(p)}}(\frac{1}{j})
\stackrel{\ref{1.14.1.1}}{=}\prod_{j\in I_p}\BWC{\sP\cap S/j}{q-1}{A}(\frac{1}{j})\\
\stackrel{\simeq(\ref{1.13})}{\longrightarrow}\BWC{S}{q-1}{A}\end{multline}
and the statements follow easily from proposition \ref{2.19.5}.

\end{proof}

\begin{defn}\label{2.21}
Let $S$ be a finite truncation set, $q\in\N$, $k$ a field with characteristic exponent $p\neq 2$, 
$C$ a smooth projective curve over $k$ with function field
$K=k(C)$ and $P\in C$ a closed point. Let $n$ be a natural number such that $k(\F^n(P))$ is separable over $k$. 
Write $P_n:=\F^n(P)\in C^{(p^n)}$, $\kappa_n=k(P_n)$ and $K_n=k(C^{(p^n)})=k(K^{p^n})$. Finally we denote by
$\widehat{K}_n$ the completion of $K_n$ in $P_n$. The choice of a local parameter $t$ in $P_n$ identifies $\widehat{K}_n=\kappa_n((t))$  
and we have a natural inclusion $\iota : K_n\inj \kappa_n((t))$. Take $\omega\in \BWC{S}{q}{K}$, then
we define the \emph{residue of $\omega$ in $P$} to be
\eq{2.21.1}{\Res^q_{S,P}(\omega)=\Res^q_P(\omega)=\Tr_{\kappa_n/k}\Bigl(\Res^q_{S,t}\bigl(\iota(\Tr_{K/K_n}(\omega))\bigr)\Bigr)\in \BWC{S}{q-1}{k}, }
where the $\Res^q_{S,t}$ on the right hand side, is the residue on $\BWC{S}{q}{\kappa_n((t))}$ from proposition \ref{2.20}.
By \ref{2.20}, (v), $\Res^q_{S,P}$ is independent of the choice of the local parameter $t$. To show that it is
independent of the choice of $n$, one reduces it, in the same way as it was done in the proof of \ref{2.14}, to 
\[\Res^q_t(\alpha)=\Res^q_z(\Tr_{\kappa((t))/\kappa((z))}(\alpha)), \text{ for } \alpha\in \BWC{S}{q}{\kappa((t))},\] 
with $\kappa=\kappa_n=\kappa_{n+1}$ and $z=t^p$. But this follows form (\ref{2.20.3})-(\ref{2.20.6}) and the corresponding property of $\Res^1$.

\end{defn}
\begin{rmk}\label{2.22}
In the situation of \ref{2.21} we have 
\[\Res_P^q(\omega)=\Res_{P_n}^q(\Tr_{K/K_n}(\omega)), \text{ for all } n\in\N. \]
\end{rmk}

\begin{rmk}\label{2.22.1}
In the situation of \ref{2.21}, we have $\Res^q_P(\omega)=0$, if $\omega\in\BWC{S}{q}{K}$ is regular in $P\in C$ (cf. remark \ref{2.15.5}). 
\end{rmk}

\begin{cor}[cf. \cite{BlEs03}, (6.14)]\label{2.23}
Let $S$ be a finite truncation set, $q\in \N$, $k$ a field of characteristic $\neq 2$ and 
$C$ a smooth projective curve over $k$ with function field $K$. Then
\[\sum_{P\in C}\Res^q_P(\omega)=0, \text{ for all }\omega\in \BWC{S}{q}{K}.\]

\end{cor}
\begin{proof} It is enough to consider the $p$-typical case.
Let $x$ be a separating transcendence basis of $K$ over $k$. Then we have  $\WC{n}{q}{K}=\W{n}{K}\otimes_{\W{n}{k[x]}}\WC{n}{q}{k[x]}$ and
the statement follows with the same reasoning as in (\ref{2.20.3})-(\ref{2.20.6}) and theorem \ref{2.17}.
\end{proof}
\pagebreak

%
%
%
%
%
%

\section{Additive cubical Chow groups with higher modulus}
\noindent
In this section we will define the  additive cubical higher Chow groups for fields on the level of zero cycles with modulus $(m+1)$, as it was done
by Bloch and Esnault in \cite[Section 6.]{BlEs03} (there for $m=1$). We will prove our main theorem that the additive Chow groups with modulus $(m+1)$
over an arbitrary field $k$ of characteristic $\neq 2$ are isomorphic to the de Rham-Witt complex of length $m$ over $k$.
This generalizes \cite[Theorem 6.4.]{BlEs03}, where the statement was proven for $m=1$. The proof is as follows: First we construct a map
from the additive Chow groups to the de Rham-Witt complex. That the map is well defined follows essentially
from the residue theorem of the previous section. This is just a generalization of the arguments in \cite{BlEs03}. 
Then we equip the additive Chow groups with the structure of a restricted Witt complex
to obtain an inverse map. To verify that the constructed maps satisfy the relations, which the multiplication, the differential, the Frobenius and the
Verschiebung should satisfy in a restricted Witt complex, we use a result of Nesterenko-Suslin and Totaro, which identifies Bloch's higher Chow groups 
of a field on the level of zero cycles with the Milnor $K$-theory, see \cite{NeSu89} and \cite{To92}.\\
\\
We will use the definitions and conventions from appendix A, without mentioning it explicitly.\\
\\
In the following $k$ is always assumed to be a field of characteristic $p\neq 2$. \\
\\
We write
\[X_{n,k}=X_n=\G_{m,k}\times_k\bigl(\P^1_k\setminus\{1\}\bigr)^n, \text{ with coordinates } (x,y_1,\ldots,y_n).\]
For $j=0,\infty$, we denote by $F^j_{n,i}=F^j_i\subset X_n$ the face given by $y_i=j$. Furthermore we define the divisors  
$F_{n,0}=F_0\subset \P^1\times(\P^1)^n=\overline{X}_n$ by $x=0$ and $F^1_{n,i}=F^1_i\subset \overline{X}_n$ by $y_i=1$.  
The {\em boundary maps}
\[\partial^j_i:\Zy_1(X_n)\longrightarrow \Zy_0(X_{n-1})\quad i=1,\ldots,n,\, j=0,\infty \]
are given by the composition of the intersection map $\Zy_1(X_n)\to \Zy_0(F^j_{n,i})$, $\alpha\mapsto F^j_i.\alpha$ 
with the canonical isomorphism $\Zy_0(F^j_{n,i})\cong\Zy_0(X_{n-1})$. Furthermore we denote
\[Y_n=\bigcup_{\substack{i=1,\ldots,n\\ j=0,\infty}} F^j_i=
      \bigcup_{i=1}^\infty \G_m\times\Bigl(\bigl(\P^1\setminus\{1\}\bigr)^{i-1}\times \{0,\infty\}\times\bigl(\P^1\setminus\{1\}\bigr)^{n-i}\Bigr).\]
Finally we define for any $m\in\N$ 
\eq{3.1.2}{\psi_{n,m}=\psi_n=\frac{1}{[x]}\frac{d[y_1]}{[y_1]}\ldots \frac{d[y_n]}{[y_n]}\in \BWC{m}{n}{\Gamma(X_n\setminus Y_n,\sO_{X_n})}.}
If $C\subset X_n$ is a curve with function field $K=k(C)$, we denote the image of $\psi_n$ in $\BWC{m}{n}{K}$ again by $\psi_n$
and if $P\in X_n\setminus Y_n$ is a point, we write $\psi_n(P)$ for the image of $\psi_n$ in $\BWC{m}{n}{k(P)}$.

\begin{defn}[\cite{BlEs03}, 6.]\label{3.1}
Let $m\ge 1$ be an integer. We define $\Zy_1(X_n;m)$ to be the subgroup of $\Zy_1(X_n),$ which is freely generated by 
1-dimensional subvarieties $C\subset X_n$, $C\not\subset Y_n$ satisfying the following properties
         \begin{enumerate}
          \item[(a)] ({\em Good position}) $\partial_i^j [C]\in\Zy_0(X_{n-1}\setminus Y_{n-1})$, for $i=1,\ldots,n$, $j=0,\infty$. 
          \item[(b)] ({\em Modulus $(m+1)$ condition}) If $\nu: \widetilde{C}\to \P^1\times \bigl(\P^1\bigr)^n$ is the normalization of the 
                     compactification of $C$, then
              \eq{3.1.4}{(m+1)[\nu^*F_0]\le\sum_{i=1}^n[\nu^*F^1_i]\quad\text{in } \Zy_0(\widetilde{C}).} 
         \end{enumerate}

We have a map $\partial=\sum_{i=1}^n (-1)^i(\partial^0_i-\partial^\infty_i):\Zy_1(X_n;m)\to\Zy_0(X_{n-1}\setminus Y_{n-1})$ 
and we define
\[\Th{n}{k}{m}=\frac{\Zy_0(X_{n-1}\setminus Y_{n-1})}{\partial\Zy_1(X_n;m)}.\]
We call these groups {\em additive cubical higher Chow groups of $k$ of level $n$, codimension $n$ and modulus $(m+1)$}.

For $m'\le m$ we have $\Zy_1(X_n;m)\subset \Zy_1(X_n;m')$, hence we naturally obtain  a projective system
\eq{3.1.1}{\bigl((\Th{n}{k}{m})_{m\in\N}, \, \R:\Th{n}{k}{m+1}\to \Th{n}{k}{m}\bigr).}
\end{defn}

\begin{rmk}\label{3.1.5}
The condition \ref{3.1.4} can be spelled out explicitly as
\[(m+1)\va_P(\nu^*x)\le\sum_{i=1}^n\va_P(\nu^*y_i-1), \text{ for all } P\in \nu^{-1}\bigl(\{0\}\times(\P^1)^n\bigr),\]
with $\va_P$ the discrete valuation of $\sO_{\widetilde{C},P}$.
\end{rmk}

\begin{prop}\label{3.1.6}
Let $C\subset X_n$ be a curve with $[C]\in\Zy_1(X_n;m)$, $K=k(C)$ its function field and 
$\nu :\widetilde{C}\to \P^1\times(\P^1)^n$ be the normalization of the compactification. 
Then $\nu^*\psi_n\in\BWC{m}{n}{K}$ has only poles in
\[\Sigma:=\bigcup_{i=1}^n\nu^{-1}\Bigl(\G_m\times\bigl(\P^1\setminus\{0,1,\infty\})^{i-1}\times\{0,\infty\}\times
          \bigl(\P^1\setminus\{0,1,\infty\}\bigr)^{n-i}\Bigr)\subset\widetilde{C}.\]
\end{prop}

Before we start with the proof of the proposition, we need a little lemma.

\begin{lem}\label{3.1.7}
Let $S$ be a finite truncation set. Then for all $n\in S$ there are $f_n\in\Z[x]$ such that
\[[1+x]=[1]+\sum_{n\in S}\V_n([xf_n(x)])\in\BW{S}{\Z[x]}.\] 
\end{lem}
\begin{proof}
We do induction over $S$. If $S=\{1\}$, we have $[1+x]=[1]+[x]$ and we are done. Now take any finite $S$ and write 
$r=\max S$ for the greatest number in $S$. Then, by induction hypothesis, there are $f_n\in\Z[x]$, 
for $n\in S\setminus\{r\}$, and  a polynomial $g\in\Z[x]$ such that 
\[[1+x]-([1]+\sum_{n\in S}\V_n ([xf_n(x)]))=\V_r([g]).\]
Applying $\gh_r$, we obtain
\[rg=(1+x)^r-(1+\sum_{\substack{n|r\\ n<r}} n(xf_n(x))^{r/n} ).\]
Now the right hand side is divisible by $x$, hence $g$ too and this gives the assertion.
\end{proof}

\begin{proof}[Proof of Proposition \ref{3.1.6}]
A priori $\nu^*\psi_n$ may have  poles only in $\nu^{-1}(x=0)$ and $\nu^{-1}(y_i=0,\infty)$, $i=1,\ldots,n$. 
Now we have to show
\begin{enumerate}
\item[(a)] $\nu^*\psi_n$ has no poles in $\nu^{-1}(x=0)$.
\item[(b)] If $\nu^*\psi_n$ has a pole in $P\in\nu^{-1}(y_i=0,\infty)$, for some $i$, then 
            $P\notin \cup_{j\neq i}\nu^{-1}(y_j=0,1,\infty)$.
\end{enumerate}
We will see that (a) holds because $C$ satisfies the modulus $(m+1)$ condition and (b) since $C$ is in good position.

We begin with the proof of (a). Take $P\in\nu^{-1}(x=0)$ and let $t\in\sO_{\widetilde{C},P}$ be a local parameter in $P$.
The modulus condition tells us in particular that there is at least one $i_0$ with $P\in\nu^{-1}({y_i}_0=1)$. 
Thus we may assume
\[P\in\bigcap_{i=1}^r\nu^{-1}(y_i=1) \text{ and } P\notin\bigcup_{i=r+1}^n\nu^{-1}(y_i=1),\quad 
  \text{for some } r \text{ with } 1\le r\le n.\]
Hence we can write
\[\nu^*x=t^av,\quad \nu^*y_i=1+t^{b_i}u_i,\text{ for } i=1,\ldots,r,\quad \nu^*y_i=t^{c_i}u_i,\text{ for } i=r+1,\ldots,n,\]
with $v,u_i\in\sO_{\widetilde{C},P}^\times$, $a,b_i\ge 1$ and $c_i\in\Z$ for $i=1,\ldots,n$, and the modulus $(m+1)$ condition
tells us 
\eq{3.1.6.1}{(m+1)a\le b_1+\ldots+ b_r.} 
Now we have in $\BWC{m}{n}{K}$
\[\nu^*\psi_n=\frac{1}{[t^av]}\frac{d[1+t^{b_1}u_1]}{[1+t^{b_1}u_1]}\cdots\frac{d[1+t^{b_r}u_r]}{[1+t^{b_r}u_r]}
                      \frac{d[t^{c_{r+1}}u_{r+1}]}{[t^{c_{r+1}}u_{r+1}]}\cdots\frac{d[t^{c_n}u_n]}{[t^{c_n}u_n]}\]
and to show that this is regular in $t=0$ amounts to show that the forms
\[\frac{1}{[t]^a}d[1+t^{b_1}u_1]\ldots d[1+t^{b_r}u_r]\frac{d[t]}{[t]}\]
and
\[\frac{1}{[t]^a}d[1+t^{b_1}u_1]\ldots d[1+t^{b_r}u_r]\]
are regular in $t=0$. We will use the shorthand $z=[t]$. Then, by lemma \ref{3.1.7}, we can write
\[[1+t^{b_i}u_i]=[1]+\sum_{j=1}^m\V_j(z^{b_i}w_{ij}),\quad 
                                                   \text{for some }w_{ij}\in\BW{m}{\sO_{\widetilde{C},P}}.\] 
Hence it is enough to show that the differentials
\[\frac{1}{z^a}d\V_{j_1}(z^{b_1}w_1)\ldots d\V_{j_r}(z^{b_r}w_r)\frac{dz}{z}\]
and
\[\frac{1}{z^a}d\V_{j_1}(z^{b_1}w_1)\ldots d\V_{j_r}(z^{b_r}w_r),\]
with $w_i\in\BW{m}{\sO_{\widetilde{C},P}}$ and $1\le j_1\le\ldots\le j_r\le m$, are regular in $z=0$.  
Now we claim that the differential $\omega=d\V_{j_1}(z^{b_1}w_1)\ldots d\V_{j_r}(z^{b_r}w_r)$ can be written as 
\eq{3.1.6.2}{d\V_{j_r}(z^{b_1+\ldots+b_r}\alpha+z^{b_1+\ldots+b_r-1}\beta dz),}
with $\beta\in\BWC{m}{r-2}{\sO_{\widetilde{C},P}}$ and $\alpha\in\BWC{m}{r-1}{\sO_{\widetilde{C},P}}$.

By induction on $r$ it is enough to show that 
\[\eta=d\V_{i}(z^ew)d\V_j(z^e\alpha'+z^b\tfrac{dz}{z}\beta'),\quad i\le j \text{ and }\alpha',\beta'\text{ appropriate},\]
can be written as $d\V_j(z^{b+e}\alpha+z^{b+e}\frac{dz}{z}\beta)$. We write $c=(i,j)$ and take $h,l$ with $hi+lj=c$. Then
by \ref{1.5}, (iv)
\[\F_id\V_j(z^b\alpha'+z^b\tfrac{dz}{z}\beta')=hd\F_{i/c}\V_{j/c}(z^b\alpha'+z^b\tfrac{dz}{z}\beta')+
                                                                              l\F_{i/c}\V_{j/c}d(z^b\alpha'+z^b\tfrac{dz}{z}\beta')\]
and since $\F_{i/c}\V_{j/c}=\V_{j/c}\F_{i/c}$ and $\F_i(\gamma)dz/z=\F_i(\gamma dz/z)$, 
there are appropriate $\alpha_0,\alpha_1,\beta_0,\beta_1$ such that this may be written as
\[d\V_{j/c}(z^{bi/c}\alpha_0+z^{bi/c}\tfrac{dz}{z}\beta_0)+\V_{j/c}(z^{bi/c}\alpha_1+z^{bi/c}\tfrac{dz}{z}\beta_1).\]
Therefore $\eta$ is the sum of the differentials
\[d\V_i(z^ewd\V_{j/c}(z^{bi/c}\alpha_0+z^{bi/c}\tfrac{dz}{z}\beta_0))=
                                               -d\V_i((ez^e\tfrac{dz}{z}w+z^edw)\V_{j/c}(z^{bi/c}\alpha_0+z^{bi/c}\tfrac{dz}{z}\beta_0))\]
and 
\[d\V_i(z^ew\V_{j/c}( z^{bi/c}\alpha_1+z^{bi/c}\tfrac{dz}{z}\beta_1)).\]
Both can be written in the form
\[d\V_j\V_{i/c}(z^{(bi+ej)/c}\alpha_2+z^{(bi+ej)/c}\tfrac{dz}{z}\beta_2)=
  d\V_j(z^{b+e}\V_{i/c}(z^{e(j-i)/c}\alpha_2)+z^{b+e}\tfrac{dz}{z}\V_{i/c}(z^{e(j-i)/c})\beta_2),\]
which is of the promised shape.

Thus we are reduced to show that 
\[\frac{1}{z^a}d\V_{j_r}(z^{b_1+\ldots+b_r}\alpha+z^{b_1+\ldots+b_r-1}\beta dz)\frac{dz}{z}\]
and
\[\frac{1}{z^a}d\V_{j_r}(z^{b_1+\ldots+b_r}\alpha+z^{b_1+\ldots+b_r-1}\beta dz),\]
with $\alpha,\beta$ as in (\ref{3.1.6.2}), are regular in $z=0$.
Since $\V_j(\gamma)dz/z=\V_j(\gamma dz/z)$, we have 
\begin{multline*}\frac{1}{z^a}d\V_{j_r}(z^{b_1+\ldots+b_r}\alpha+z^{b_1+\ldots+b_r-1}\beta dz)\frac{dz}{z}=
\frac{1}{z^a}d\V_{j_r}(z^{b_1+\ldots+b_r-1}\alpha dz)\\
 =d\V_{j_r}(z^{b_1+\ldots+ b_r-1-j_ra}\alpha dz)-
                                       d(\frac{1}{z^a})\V_{j_r}(z^{b_1+\ldots+ b_r-1}\alpha dz).
\end{multline*}
The second term is zero and the first is regular in $z=0$ because $j_ra+1\le (m+1)a\le b_1+\ldots+b_r$ by (\ref{3.1.6.1}).
Without the $\frac{dz}{z}$-term at the end it works the same. This proves (a).

Now (b). We know that $C$ is in good position, i.e.
\[C\cap \bigl(\G_m\times(\P^1\setminus\{1\})^{i-1}\times\{0,\infty\}\times(\P^1\setminus\{1\})^{n-i}\bigr)\subset
  \G_m\times\bigl(\P^1\setminus\{0,1,\infty\})^{i-1}\times\{0,\infty\}\times \bigl(\P^1\setminus\{0,1,\infty\}\bigr)^{n-i}.\]
Thus a point $P\in \widetilde{C}\setminus\nu^{-1}(x=0)$, which lies in $\nu^{-1}(y_i=0,\infty)$ for more than one $i$, 
must lie in $\nu^{-1}(y_j=1)$ for at least one j.
We may assume $j=1$ and then it is enough to show that $\nu^*\psi_n$ is regular in $\nu^{-1}(y_1=1)\cap\nu^{-1}(x\neq 0)$. 
Take $P\in\nu^{-1}(y_1=1)\cap\nu^{-1}(x\neq 0)$ and let $t$ be a local parameter in $P$. Since  $\nu^*x\in\sO_{\widetilde{C},P}^\times$, we see
that the question of regularity of $\nu^*\psi_n$ in $P$ reduces to show
\[d[1+t^bu]\frac{d[t]}{[t]},\quad b\ge 1, u\in\sO_{\widetilde{C},P}^\times\]
is regular in $t=0$. By lemma \ref{3.1.7} it is enough to prove the regularity in $z=0$ of
\[d\V_j(z^bw)\frac{dz}{z}=d\V_j(z^{b-1}wdz),\quad w\in\BW{m}{\sO_{\widetilde{C},P}},\]
which is obvious.

\end{proof}

\begin{thm}[cf. \cite{BlEs03}, Proposition 6.2.]\label{3.2}
Let $m$ be a positive number. Then the group homomorphism
\[\theta:\Zy_0(X_{n-1}\setminus Y_{n-1})\longrightarrow \BWC{m}{n-1}{k}, \quad [P]\mapsto \Tr_{k(P)/k}(\psi_n(P))\]
factors to give a map (which we will call $\theta$ again)
\[\theta: \Th{n}{k}{m}\longrightarrow \BWC{m}{n-1}{k}.\]
\end{thm}
\begin{proof}
Let $C\subset X_n$ be a curve with $[C]\in\Zy_1(X_n;m)$. We have to show $\theta(\partial[C])=0$ in $\BWC{m}{n-1}{k}$.
Denote by $\nu: \widetilde{C}\to \P^1\times(\P^1)^n$ the normalization of the compactification of $C$. It follows from
proposition \ref{3.1.6} and remark \ref{2.22.1} that $\Res_P(\nu^*\psi_n)=0$ for all $P\in\widetilde{C}\setminus\Sigma$.
Next we want to calculate $\Res_P(\nu^*\psi_n)$ for $P\in \Sigma$.

Take $P\in \nu^{-1}(y_n=0)\cap\Sigma$. 
Choose $j\ge 0$, such that $\kappa_j=k(\F^j(P))\supset k$ is 
separable and write $K$ (resp. $K_j$) for the function field of $\widetilde{C}$ (resp. $\widetilde{C}^{(p^j)}$). Furthermore write
\[ e_P=p^r,\quad f_P=p^s,\]
then by lemma \ref{3.3}
\eq{3.3.1}{j=r+s.}
Let $z$ be a local parameter in $\F^j(P)$ and $t$ one in $P$. Then we may write
\eq{3.3.2}{z=t^{e_P}u=t^{p^r}u, \text{ with } u\in \sO_{\widetilde{C},P}^\times}
and
\[\nu^*x=u_0,\quad \nu^*y_i=u_i,\, i=1,\ldots,n-1,\quad \nu^*y_n=t^au_n,\]
with $u_i\in\sO_{\widetilde{C},P}^\times$, $i=0,\ldots,n$ and $a=\va_P(\nu^*y_n)$. Therefore in $\BWC{p^{j(n-1)+r}m}{n}{k}$ 
\begin{multline*}
\underline{p}^{j(n-1)+r}\Res_P(\nu^*\psi_n)  =  
a\Tr_{\kappa_j/k}\Res_z\Tr_{K/K_j}
             \Bigl(\frac{1}{[u_0]}\frac{d[u_1^{p^j}]}{[u_1^{p^j}]}\cdots\frac{d[u_{n-1}^{p^j}]}{[u_{n-1}^{p^j}]}\frac{d[z]}{[z]}\Bigr)\\
 +  \Tr_{\kappa_j/k}\Res_z\Tr_{K/K_j}\Bigl(\frac{1}{[u_0]}\frac{d[u_1^{p^j}]}{[u_1^{p^j}]}\cdots\frac{d[u_{n-1}^{p^j}]}{[u_{n-1}^{p^j}]}
                                                                                    \frac{d[u_n^{p^r}/u^a]}{[u_n^{p^r}/u^a]}\Bigr),
\end{multline*}
here $\underline{p}$ is the map from definition \ref{2.6}, which is,  lifting to the level $pm$ and then multiplication with $p$.
By remark \ref{2.22.1}, the second term on the right hand side is zero and since $\Tr_{K/K_j}(\frac{1}{[u_0]})=\V_{p^j}(\frac{1}{[u_0^{p^j}]})$, we obtain
by the linearity of the trace and proposition \ref{2.20} (vii)
\[\underline{p}^{j(n-1)+r}\Res_P(\nu^*\psi_n)=a\Tr_{\kappa_j/k}\left(\V_{p^j}\Bigl(\frac{1}{[\nu^*x(P)^{p^j}]}\Bigr)
                                           \frac{d[\nu^*y_1(P)^{p^j}]}{[\nu^*y_1(P)^{p^j}]}\cdots
                                                    \frac{d[\nu^*y_{n-1}(P)^{p^j}]}{[\nu^*y_{n-1}(P)^{p^j}]}\right).\]
Denoting $k(P)=\kappa$, we have $f_P=[\kappa:\kappa_j]=p^s$ and thus
\[\underline{p}^{j(n-1)+r}\Res_P(\nu^*\psi_n)  =  a \Tr_{\kappa_j/k}\Tr_{\kappa/\kappa_j}\left(p^{j(n-1)+j-s}
    \frac{1}{[\nu^*x(P)]}\frac{d[\nu^*y_1(P)]}{[\nu^*y_1(P)]}\cdots \frac{d[\nu^*y_{n-1}(P)]}{[\nu^*y_{n-1}(P)]}\right)\]
\[= ap^{j(n-1)+r}\Tr_{\kappa/k}\left(
      \frac{1}{[x(\nu(P))]}\frac{d[y_1(\nu(P))]}{[y_1(\nu(P))]}\cdots \frac{d[y_{n-1}(\nu(P))]}{[y_{n-1}(\nu(P))]}\right)\]
\[= \underline{p}^{j(n-1)+s}a[\kappa:k(\nu(P))]\Tr_{k(\nu(P))/k}\Bigl(\psi_{n-1}\bigl(\nu(P)^0_n\bigr)\Bigr),\]
where $\nu(P)^0_n$ is the point $\nu(P)\in F^0_n$ viewed as a point in $X_{n-1}$ via the canonical isomorphism $F^0_n\cong X_{n-1}$.
Since $\underline{p}$ is injective, we obtain in $\BWC{m}{n-1}{k}$
\eq{1.3.4}{\Res_P(\nu^*\psi_n)=\va_P(\nu^*y_n) \,[k(P):k(\nu(P))]\, \theta\bigl([\nu(P)^0_n]\bigr).  }
Similarly, one shows
\[\Res_P(\nu^*\psi_n)= (-1)^{n-i}\va_P(\nu^*y_i)\,[k(P):k(\nu(P))]\, \theta\bigl([\nu(P)^j_i]\bigr),\]
for all $P\in \Sigma\cap \nu^{-1}(y_i=j)$ and $j=0,\infty$, with $\nu(P)^j_i$ being the point $P$ viewed as a point in $X_{n-1}$. 
Thus by corollary \ref{2.23}, proposition \ref{3.1.6} and example \ref{A3} we obtain
\begin{eqnarray*}
0 & = & \sum_{Q\in\Sigma}\Res_Q(\nu^*\psi_n)\\
  & = & \sum_{i=1}^{n}\sum_{Q\in\Sigma\cap\nu^{-1}(y_i=0)}(-1)^{n-i}\va_Q(\nu^*y_i)\,[k(Q):k(\nu(Q))]\, \theta\bigl([\nu(Q)^0_i]\bigr)\\
  &\phantom{=} & 
           -\sum_{Q\in\Sigma\cap\nu^{-1}(y_i=\infty)}-(-1)^{n-i}\va_Q(\nu^*y_i)\,[k(Q):k(\nu(Q))]\, \theta\bigl([\nu(Q)^\infty_i]\bigr)\\
  & = & (-1)^n\sum_{i=1}^n\sum_{P\in C\cap (y_i=0)}(-1)^i\ord_P(y_i)\theta([P^0_i]))
                                            -\sum_{P\in C\cap (y_i=\infty)}(-1)^i(-\ord_P(y_i))\theta([P^\infty_i])\\
 & =& (-1)^n \sum_{i=1}^n(-1)^i\theta({\partial_i^0}[C]-{\partial_i^\infty}[C]) = (-1)^n\theta(\partial [C])
\end{eqnarray*}
and this proves the theorem.
\end{proof}

\begin{cor}\label{3.4}
The map 
\eq{3.4.1}{\theta: \Th{1}{k}{m}\longrightarrow \BW{m}{k},\quad [P]\mapsto \Tr_{k(P)/k}(\tfrac{1}{[x(P)]})}
is an isomorphism of groups. (Notice that the brackets $[ \ ]$ on the left  denote a cycle class, whereas on the right the Teichm\"uller lift.) 
\end{cor}

\begin{proof}
First we observe that we can represent every cycle class $\xi\in \Th{1}{k}{m}$ by a principal divisor
\eq{3.4.2}{\xi=\Div(f)-\Div(g),\quad f,g\in (k^\times+xk[x]). }
(By abuse of notation we write $\Div(f)$ for the residue class of the divisor of f in $\Th{1}{k}{S}$.) Now we claim
\eq{3.4.3}{\theta(\Div(f)-\Div(g))= w(\frac{f(T)}{f(0)})-w(\frac{g(T)}{g(0)}),     }
where  $w(\frac{f(T)}{f(0)})$ is the Witt vector coming from $\frac{f(T)}{f(0)}\in \bigl(\frac{1+Tk[[T]]}{1+T^{m+1}k[[T]]}\bigr)^\times$
(see \ref{0.10}). Obviously it is enough to check
\[\theta(\Div(1-xh(x))= w(1-Th(T)),\]
for $h=\sum_{i=1}^na_ix^i\in k[x]$, such that $1-xh(x)\in k[x]$ is irreducible. 
Then $\Div(1-xh(x))=P\in \G_m$ is a point with residue field $k(P)=k[x,\frac{1}{x}]/(1-xh(x))$. Denoting the residue class of $x$ in $k(P)$ by $\alpha$,
we have
\[\theta([P])=\Tr_{k(P)/k}(\frac{1}{[\alpha]})=\Tr_{k(P)/k}\bigl([h(\alpha)]\bigr)= w\bigl(\Nm_{k(P)[[T]]/k[[T]]}(1-Th(\alpha))\bigr).\]
Now the matrix of $1-Th(\alpha)$ in the basis $1, \alpha,\ldots,\alpha^n$ is
\[M_n=\left(\begin{array}{cccccc}
          1-a_0T & -T & 0 & \cdots & 0& 0\\
          -a_1T  & 1  & -T & \ddots&  0& 0\\
          -a_2T & 0  & 1   & \ddots & 0& 0\\
           \vdots& \vdots& \ddots& \ddots  &\ddots &  \vdots\\
           -a_{n-1}T & 0 & 0 &\cdots & 1 &-T\\
           -a_nT     & 0 & 0 &\cdots & 0 & 1

  \end{array}\right).\]
By induction over $n$ (Laplace expansion with respect to the last column) we see
\[\det M_n= 1-Th(T).\]
This proves the claim (\ref{3.4.3}). It follows immediately that $\theta$ is surjective. Thus it remains to check the injectivity of $\theta$.
Take $\xi\in \Th{1}{k}{m}$ with $\theta(\xi)=0$. Write
\[\xi=\Div f-\Div g,\quad f,g\in (1+xk[x]), \text{ with } (f,g)=1.\]
Then by (\ref{3.4.3}) we have
\eq{3.4.4}{ \frac{f(T)}{g(T)}\in(1+T^{m+1}k[[T]])^\times.  }
Now we define a curve $C\subset \G_m\times\P^1\setminus\{1\}$ by
\[C \,:\, f(x)-yg(x)=0.\]
Since $(f,g)=1$ it follows from the Eisenstein criterion that $C$ is integral. Furthermore $\overline{C}$ is nonsingular in $P=(0,1)$, we have 
$C\not\subset \{y=0\}\cup\{y=\infty\}$ and $\overline{C}\cap(\{0\}\times\P^1)= P$. Additionally $x$ is a local parameter in $P$ and thus it follows
from (\ref{3.4.4}) and remark \ref{3.1.5} that $C$ fulfills the modulus condition. Hence $[C]\in\Zy_1(X_1;m)$ and $\partial[C]=\xi$.
\end{proof}

Next we want to equip the additive Chow groups with the structure of a restricted Witt complex. Then the universality of the de Rham-Witt complex
will yield an inverse map to the map $\theta$ of theorem \ref{3.2}.
\begin{smox}\label{3.5}
\hfill
\begin{enumerate}
\item Let $\mu_0: \G_m\times\G_m\to\G_m$ be the multiplication morphism on $\G_m$ (i.e. the morphism induced by $k[x,1/x]\to k[x_1,1/x_1,x_2,1/x_2]$, 
      $x\mapsto x_1x_2$). Then we define the morphism
      \eq{3.5.1}{\mu_{r,s}=\mu: X_r\times X_s\longrightarrow X_{r+s}}
      to be the composition
      \[X_r\times X_s\cong \G_m\times\G_m\times (\P^1\setminus\{1\})^r\times(\P^1\setminus\{1\})^s
                                                        \stackrel{\mu_0\times \id}{\longrightarrow}\G_m\times(\P^1\setminus\{1\})^{r+s}=X_{r+s}.\]
      If $P=(a,b_1,\ldots,b_r)\in X_r$ and $P'=(a',c_1,\ldots,c_s)\in X_s$ are $k$-rational points, then
      \[\mu(P\times P')=(aa',b_1,\ldots,b_r,c_1,\ldots,c_s)\in X_{r+s}.\]
      If $P\in X_r$ is any closed point, then $\mu_{|P\times X_s}$ is finite.
\item Let $\Delta: \G_m\to \G_m\times\G_m$ be the diagonal embedding. Then we define the morphism
      \eq{3.5.2}{\Delta_n: (\G_m\setminus\{1\})\times (\P^1\setminus\{1\})^{n-1}\longrightarrow X_n}
      to be the composition
      \[(\G_m\setminus\{1\})\times (\P^1\setminus\{1\})^{n-1}\stackrel{\Delta\times \id}{\longrightarrow}
                            (\G_m\setminus\{1\})\times(\G_m\setminus\{1\})\times (\P^1\setminus\{1\})^{n-1}\inj X_n,\]
      where the last map is the natural inclusion. If $P=(a,b_1,\ldots,b_{n-1})\in X_{n-1}$ is a $k$-rational point with $a\neq 1$, then
       \[\Delta_n(P)=(a,a,b_1,\ldots,b_{n-1}).\]
       $\Delta_n$ is a closed immersion.
\item Let $r\ge 1$ be a number and $\varphi_{r,0}:\G_m\to \G_m$ the morphism induced by $k[x,1/x]\to k[x,1/x]$, $x\mapsto x^r$. We define
      the morphism $\varphi_{r,n}=\varphi_r$ by 
      \eq{3.5.3}{\varphi_r=\varphi_{r,0}\times\id: X_n=\G_m\times(\P^1\setminus\{1\})^n\longrightarrow X_n.}
      If $P=(a,b_1,\ldots,b_n)\in X_n$ is a $k$-rational point, then 
      \[\varphi_r(P)=(a^r,b_1,\ldots,b_n).\]
      For all $r\ge 1$, $\varphi_r$ is flat and finite of degree $r$. 
\end{enumerate}
\end{smox}

\begin{defn-prop}\label{3.6}
Let $m,n,r$ and $s$ be natural numbers.
\begin{enumerate}
\item The map 
      \[k\longrightarrow \Zy_0(\G_m),\quad a\mapsto\Div(1-ax)\]
      induces a map
      \eq{3.6.1}{\{-\}:k\longrightarrow \Th{1}{k}{m},\quad a\mapsto \{a\}.}
\item The composition
      \[\Zy_0(X_{r-1})\otimes\Zy_0(X_{s-1})\stackrel{\times}{\longrightarrow}\Zy_0(X_{r-1}\times X_{s-1})
                                                                                           \stackrel{\mu_*}{\longrightarrow}\Zy_0(X_{r+s-2})\]
      induces a map
      \eq{3.6.2}{\st:\Th{r}{k}{m}\otimes\Th{s}{k}{m}\longrightarrow\Th{r+s-1}{k}{m},\quad \xi\otimes\eta\mapsto \xi\st\eta .}
\item The map
       \[-{\Delta_n}_*:\Zy_0(X_{n-1})\longrightarrow \Zy_0(X_n),\]
      where we define ${\Delta_n}_*([P]):=0$ for all closed points $P\in \{1\}\times(\P^1\setminus\{1\})^{n-1}$,
      induces a map
      \eq{3.6.3}{\sD_n=\sD=-{\Delta_n}_*: \Th{n}{k}{m}\longrightarrow\Th{n+1}{k}{m}.}
\item The map
      \[\varphi_r^*: \Zy_0(X_{n-1})\longrightarrow \Zy_0(X_{n-1})\]
      induces a map
      \eq{3.6.4}{\sV_r: \Th{n}{k}{m}\longrightarrow\Th{n}{k}{rm+r-1}. }
\item The map
      \[{\varphi_r}_*: \Zy_0(X_{n-1})\longrightarrow \Zy_0(X_{n-1})\]
      induces a map
      \eq{3.6.5}{\sF_r : \Th{n}{k}{rm+r-1}\longrightarrow\Th{n}{k}{m}.}
\end{enumerate}
\end{defn-prop}
\begin{proof}
For (i) nothing is to show. Since the composition in (ii) obviously restricts to 
$\Zy_0(X_{r-1}\setminus Y_{r-1})\otimes\Zy_0(X_{s-1}\setminus Y_{s-1})\to\Zy_0(X_{r+s-2}\setminus Y_{r+s-2})$, (ii) amounts to show that
for a curve $C\subset X_s$ with $[C]\in \Zy_1(X_s;m)$ and a point $P\in X_{r-1}\setminus Y_{r-1}$ there is a 1-cycle $\xi\in\Zy_1(X_{s+r-1};m)$ with
\eq{3.6.6}{\mu_*([P]\times\partial[C])=\partial \xi.}
We set
\[\xi:=\mu_*([P\times C])\in\Zy_1(X_{r+s-1}).\]
Then with the use of proposition \ref{A11} in the appendix we see
\[F^j_{r+s-1,i}.\xi=\mu_*(\mu^*(F^j_{r+s-1,i}).[P\times C])=
\begin{cases}
\mu_*( p_{r-1}^*(F^j_{r-1,i}).[P\times C]) &  1\le i\le r-1\\
\mu_*( p_s^*(F^j_{s,i-(r-1)}).[P\times C]) &  r\le i\ \le r+s-1,
\end{cases}\]
where $p_{r-1}$ and $p_s$ are the projections from $X_{r-1}\times X_s$ to $X_{r-1}$ and $X_s$ respectively.
Therefore $F^j_{r+s-1,i}.\xi$ is $0$ for $1\le i\le r-1$ and it equals $\mu_*([P]\times(F^j_{s,i-(r-1)}.[C]))$ for $r\le i\le s$.
This yields (\ref{3.6.6}) and that $\xi$ is in good position. It remains to show that $\xi$ satisfies the modulus $(m+1)$ condition.
Let $C''\subset \mu(P\times C)_{red}$ be an irreducible component and $C'\subset (P\times C)_\red$ an irreducible component mapping to $C''$
under the finite map $(P\times C)_\red\surj\mu(P\times C)_\red$. Denoting by $\pi$ the projection 
$\overline{X}_{r-1}\times\overline{X}_s\to\overline{X}_{s}$, with $\overline{X}_n=\P^1\times(\P^1)^n$, 
compactification and normalization give us the following diagram
\[\xymatrix{  \widetilde{C}\ar[d]^\nu & \widetilde{C'}\ar[d]^{\nu'}\ar[l]_{\tilde{\pi}}\ar[r]^{\tilde{\mu}} & \widetilde{C''}\ar[d]^{\nu''}\\
              \overline{X}_s      &  \overline{X}_{r-1}\times\overline{X}_s\ar[l]_(.6){\pi}         & \overline{X}_{r+s-1}.
            }\]
Now notice, that $\tilde{\pi}^*\nu^*(F_{s,0})=\tilde{\mu}^*{\nu''}^*(F_{r+s-1,0})$ (here we need that the x-coordinate of $P$ is not zero) and
$\tilde{\pi}^*\nu^*(F^1_{s,i})=\tilde{\mu}^*{\nu''}^*(F^1_{r+s-1, i+r-1})$. Thus, applying $\tilde{\mu}_*\tilde{\pi}^*$ to the modulus condition 
satisfied by $C$, we obtain
\[\text{deg}(\tilde{\mu})(m+1)[{\nu''}^*F_{r+s-1,0}]= (m+1)\tilde{\mu}_*\tilde{\pi}^*[\nu^*F_{s,0}]\le 
                 \sum_{i=1}^s\tilde{\mu}_*\tilde{\pi}^*[\nu^*F^1_{s,i}]=\text{deg}(\tilde{\mu})\sum_{i=1}^{r+s-1}[{\nu''}^*F^1_{r+s-1,i}].\]
(Here we need that $\tilde{\pi}$ is flat, which is the case by \cite[III, Proposition 9.7.]{Ha77}.)
This proves the modulus condition for $C''$ and thus $\xi\in\Zy_1(X_{r+s-1};m)$. 

Next (iii). Obviously ${\Delta_n}_*$ restricts to $\Zy_0(X_{n-1}\setminus Y_{n-1})\to\Zy_0(X_n\setminus Y_n)$. Now let $[C]\in \Zy_1(X_n;m)$ be a curve.
We have to show that there is a 1-cycle $\xi\in\Zy_1(X_{n+1},m)$ with
\eq{3.6.7}{{\Delta_n}_*(\partial [C])=\partial\xi.}
Write $C_0= C\cap \Bigl((\G_m\setminus\{1\})\times (\P^1\setminus\{1\})^n\Bigr)$. We set
\[\xi:={\Delta_n}_*[C_0]=[\Delta_n(C_0)]\in\Zy_1\bigl((\G_m\setminus\{1\})\times (\P^1\setminus\{1\})^{n+1}\bigr).\]
But $\Delta_n(C_0)$ is closed in $X_{n+1}$ (since a point in $\overline{\Delta_n(C_0)}\setminus \Delta_n(C_0)\subset X_{n+1}$ must have 
the $x$- and the $y_1$-coordinate equal to $1$, which is impossible by the definition of $X_{n+1}$). Thus we may view $\xi$ as a 1-cycle on $X_{n+1}$.
Now $\xi$ lives on $(\G_m\setminus\{1\})\times (\G_m\setminus\{1\})\times(\P^1\setminus\{1\})^n$, hence,  for $j=0,\infty$ 
\[F^j_{n+1,1}.\xi=0\]
and since ${\Delta_n}_*=0$ on $\{1\}\times(\P^1\setminus\{1\})^n$, 
\[F^j_{n+1,i}.\xi={\Delta_n}_*({\Delta_n}^*(F^j_{n+1,i}).[C_0])={\Delta_n}_*(F^j_{n,i-1}.[C])\in\Zy_0(X_n\setminus Y_n)\quad \text{for }i\ge 2.\]
Thus $\xi$ is in good position and satisfies (\ref{3.6.7}). It remains to check the modulus condition. Compactifying and normalizing 
the map $C\to \Delta_n(C_0)$ give us a diagram
\[\xymatrix{ \widetilde{C}\ar[r]^{\tilde{\Delta}_n}\ar[d]^\nu & \widetilde{\Delta_n(C_0)}\ar[d]^{\nu'}\\
             \overline{X}_n\ar[r]^{\bar{\Delta}_n} & \overline{X}_{n+1}.
            }\]
We have $\nu^*F_{n,0}={\tilde{\Delta}_n}^*{\nu'}^*(F^1_{n+1,0})$ and $\nu^*F^1_{n,i}={\tilde{\Delta}_n}^*{\nu'}^*(F^1_{n+1,i+1})$, $i=1,\ldots,n$.
Now we apply $\widetilde{\Delta}_{n*}$ to the modulus condition for $C$ and, as in (ii), we obtain the modulus condition for $\xi$.

(iv). Let $[C]\in\Zy_1(X_n;m)$ be a curve. We have to find a cycle $\xi\in\Zy_1(X_n;rm+r-1)$ with
\eq{3.6.8}{ \varphi_r^*(\partial[C])=\partial\xi.}
We set 
\[\xi:=\varphi_r^*[C]\in\Zy_1(X_n).\]
Then for $j=0,\infty$ and $i=1,\ldots,n$
\[F^j_i.\xi=\varphi_r^*(F^j_i).\varphi_r^*[C]=\varphi_r^*(F^j_i.[C])\in\Zy_0(X_n\setminus Y_n).\]
Hence $\xi$ is in good position and fulfills (\ref{3.6.8}). It remains to show that $\xi$ satisfies the modulus condition.
Let $C'\subset \varphi_r^{-1}(C)$ be an irreducible component. Compactifying and normalizing the map $\varphi_r: C'\to C$, yield
the diagram
\[\xymatrix{\widetilde{C'}\ar[r]^{\tilde{\varphi}_r}\ar[d]^{\nu'} & \widetilde{C}\ar[d]^\nu\\
            \overline{X}_n\ar[r]^{\bar{\varphi}_r}               & \overline{X}_n
           }\]
 Since $\nu'^*F_0$ is the Cartier divisor defined by $\nu'^*x=0$ and $\nu'^*(F^1_j)$ the one defined by $\nu'^*y_i=1$, 
we have by the definition of $\varphi_r$, 
\eq{3.6.9}{r[\nu'^*F_0]=[\tilde{\varphi}_r^*\nu^*(F_0)]=\tilde{\varphi}_r^*[\nu^*F_0]\text{ and } [\nu'^*F^1_i]=\tilde{\varphi}_r^*[\nu^* F^1_i].}
Thus
 \[((mr+r-1)+1)[\nu'^*(F_0)]=(m+1)\tilde{\varphi}_r^*[\nu^*F_0]\le \sum_{i=1}^n{\varphi_r}^*[\nu^*F^1_i]=\sum_{i=1}^n[\nu'^*F^1_i],\]
 hence $\xi\in\Zy_1(X_n;rm+r-1)$.
 
Finally (v). Let $[C']\in \Zy_1(X_n;rm+r-1)$ be a curve. The usual argument shows that 
\[\xi:={\varphi_r}_*([C'])\in \Zy_1(X_n)\]
is in good position and satisfies $\partial \xi={\varphi_r}_*\partial[C']$. Compactifying and normalizing the map $\varphi_r: C'\to C=\xi_\red$ give us
a diagram
\[\xymatrix{\widetilde{C'}\ar[r]^{\tilde{\varphi}_r}\ar[d]^{\nu'} & \widetilde{C}\ar[d]^{\nu}\\
            \overline{X}_n\ar[r]^{\bar{\varphi}_r}               & \overline{X}_n.
           }\]
Applying $\tilde{\varphi}_{r*}$ to the modulus condition of $C'$ and using the relations (\ref{3.6.9}), we obtain
\[\deg\tilde{\varphi}_r (m+1)[\nu^*F_0]=((rm+r-1)+1) \tilde{\varphi}_{r*} [\nu'^* F_0]\le \sum_{i=1}^n \tilde{\varphi}_{r*} [\nu'^* (F^1_i)]
      = \deg\tilde{\varphi}_r\sum_{i=1}^n[\nu^*F^1_i]\]
and we are done.

\end{proof}
Next we investigate the behavior of the additive Chow groups under finite field extensions.

\begin{lem}\label{3.7}
Let $m,n\ge 1$ be integers, $L\supset k$ be a finite field extension and denote by $\pi: \Spec L\to\Spec k$ the induced map. Then
$\pi_*: \Zy_0(X_{n-1,L}\setminus Y_{n-1,L})\to\Zy_0(X_{n-1,k}\setminus Y_{n-1,k})$ induces a group homomorphism (also denoted by $\pi_*$)
\[\pi_*:\Th{n}{L}{m}\longrightarrow\Th{n}{k}{m}\]
and $\pi^*:\Zy_0(X_{n-1,k}\setminus Y_{n-1,k})\to\Zy_0(X_{n-1,L}\setminus Y_{n-1,L})$ induces a group homomorphism
\[\pi^*:\Th{n}{k}{m}\longrightarrow\Th{n}{L}{m},\]
satisfying the following properties
\begin{enumerate}
\item $\pi_*$ and $\pi^*$ commute with $\sD,\sV_r$ and $\sF_r$, all $r$.
\item For $\eta\in\Th{n}{k}{m}$ and $\xi\in\Th{n}{L}{m}$ we have
      \[\pi_*((\pi^*\eta)\st_L\xi)=\eta\st_k\pi_*\xi.\]
\item The following diagrams commute
     \[\xymatrix{\Th{n}{L}{m}\ar[r]^(.6){\theta_L}\ar[d]^{\pi_*} & \BWC{m}{n-1}{L}\ar[d]^{\Tr_{L/k}}\\
                 \Th{n}{k}{m}\ar[r]^(.6){\theta_k}               & \BWC{m}{n-1}{k}
                 }\]
     and
     \[\xymatrix{\Th{n}{L}{m}\ar[r]^(.6){\theta_L} & \BWC{m}{n-1}{L}                                \\
                 \Th{n}{k}{m}\ar[r]^(.6){\theta_k}\ar[u]_{\pi^*} & \BWC{m}{n-1}{k}\ar[u]_{\pi^*},
                 }\]
     where the $\pi^*$ on the right hand side is the natural map induced by $k\subset L$.
\end{enumerate}
\end{lem}
\begin{proof}
First we show, $\pi_*$ is well defined on the additive Chow groups. Therefore take a curve $[C]\in\Zy_1(X_{n,L};m)$ and define 
$\eta:=\pi_*[C]\in\Zy_1(X_{n,k})$. One easily checks that $\eta$ is in good position and satisfies $\partial\eta=\pi_*\partial [C]$.
Compactifying and normalizing the map $C\to \pi(C)$ yield
\[\xymatrix{ \widetilde{C}\ar[r]^{\tilde{\pi}}\ar[d]^{\nu} & \widetilde{\pi(C)}\ar[d]^{\nu'}\\
             \overline{X}_{L,n}\ar[r]^{\bar{\pi}} & \overline{X}_{k,n} .
            }\]
We have 
\eq{3.7.1}{\nu^*F_{0,L}=\tilde{\pi}^*{\nu'}^*(F_{0,k}) \text{ and } \nu^*F^1_{i,L}=\tilde{\pi}^*{\nu'}^*(F^1_{i,k})}
and thus applying $\tilde{\pi}_*$ to the modulus condition for $C$ gives us the modulus condition for $\eta$. 

To show that $\pi^*$ is well defined we take any curve $[C']\in\Zy_1(X_{n,k};m)$ and define $\xi:=\pi^*[C']\in\Zy_1(X_{n,L})$. Again
$\xi$ is in good position and fulfills $\partial\xi=\pi^*\partial[C']$. Now let $C\subset \pi^{-1}(C')$ be an irreducible component.
Compactifying and normalizing the map $\pi:C\to C'$ yield
\[\xymatrix{\widetilde{C}\ar[r]^{\tilde{\pi}}\ar[d]^{\nu} & \widetilde{C'}\ar[d]^{\nu'}\\
             \overline{X}_{L,n}\ar[r]^{\bar{\pi}} & \overline{X}_{k,n} .
           }\]
Since (\ref{3.7.1}) still holds we obtain the modulus condition for $\xi$ by applying $\tilde{\pi}^*$ to the modulus condition for $C'$. 

The proof of (i) and (ii) is immediate (for (ii) use lemma \ref{A9}). To check the commutativity of the first diagram in (iii), 
we take a point $P\in X_{n-1,L}\setminus Y_{n-1,L}$, notice that we have inclusions $k\subset k(\pi(P))\subset k(P)\supset L$,  then
the statement follows from 
\begin{eqnarray*}
\theta_k(\pi_*[P]) & = & [k(P):k(\pi(P))]\Tr_{k(\pi(P))/k}\Bigl(\psi_{n-1}(\pi(P))\Bigr)=\Tr_{k(P)/k}(\psi_{n-1}(P))\\
                   & = & \Tr_{L/k}\Tr_{k(P)/L}(\psi_{n-1}(P)) = \Tr_{L/k}\theta_L([P]).
\end{eqnarray*}
It remains to check the second diagram in (iii). Thus take $Q\in X_{n-1,k}\setminus Y_{n-1,k}$ and denote $E=k(Q)\supset k$.
We write
\[\pi^*[Q]=\sum_i l_i[R_i],\]
where the $R_i\in X_{n-1,L}\setminus Y_{n-1,L}$ are all the points mapped to $Q$ under $\pi$. If we further denote $L_i=k(R_i)$ and 
let $\sigma_i : E \inj L_i$ be the natural inclusion, then 
\[\theta_L(\pi^*[Q])=\sum_i l_i \Tr_{L_i/L}(\psi_{n-1}(R_i))=\sum_i\Tr_{L_i/L}\Bigl(l_i\sigma_i\bigl(\psi_{n-1}(Q)\bigr)\Bigr),\]
which by proposition \ref{2.7.5} equals
\[\pi^*\Tr_{E/k}(\psi_{n-1}(Q))=\pi^*\theta_k([Q])\]
and we are done. 
\end{proof}
\begin{cor}\label{3.8}
The map $\theta : \Th{n}{k}{m}\to \BWC{m}{n-1}{k}$ from theorem \ref{3.2} satisfies
\eq{3.8.1}{\theta\circ\{-\}=[-],\quad\theta\circ\sD=d\circ\theta,\quad \theta\circ\sV_r=\V_r\circ\theta,\quad \theta\circ\sF_r=\F_r\circ\theta,\quad 
                                                                                                                          \text{for all }r\ge 1}
and
\eq{3.8.2}{\theta(\eta\st \xi)=\theta(\eta)\theta(\xi),\quad \text{for all }\eta\in\Th{r}{k}{m}, \xi\in\Th{s}{k}{m}, r,s\ge 1.}
\end{cor}
\begin{proof}
The first equality in (\ref{3.8.1}) follows immediately from the definitions of $\theta$ and $\{-\}$. 
It is enough to prove the other equalities in (\ref{3.8.1})
for $k$-rational points. Indeed, if $P\in X_{n-1,k}$ is a closed point with residue field $k(P)=L$, we denote by $a, b_1,\ldots,b_{n-1}\in L$ the residue
classes of the coordinates in $L$ and by $P'=(a,b_1,\ldots,b_{n-1})\in X_{n-1,L}$. Now let $\pi:\Spec L \to \Spec k $ be the
map induced by $k\subset L$, then we have $\pi_*[P']=[P]\in\Th{n}{k}{m}$. Thus if we know the equalities for $k$-rational points on $X_{n-1,k}$, 
for any field $k$, we know it in general by lemma \ref{3.7} and the fact that on the side of the de Rham-Witt complex the trace commutes with
$\V_r$, $\F_r$ and $d$. Thus it is enough to consider points $P=(1/a, b_1,\ldots, b_{n-1})\in X_{n-1}\setminus Y_{n-1}$ with $a$ and $b_i\in k^\times$.
For such points
\[\theta(\sD[P])=-\theta(\frac{1}{a},\frac{1}{a},b_1,\ldots,b_{n-1})=[a]\frac{d[a]}{[a]}\frac{d[b_1]}{[b_1]}\ldots\frac{d[b_{n-1}]}{[b_{n-1}]}
                = d \theta([P])\]
(here $a\neq 1$ and it is trivially true for $a=1$) and
\[\theta(\sF_r[P])=\theta(\frac{1}{a^r},b_1,\ldots,b_{n-1})=\F_r([a])\frac{d[b_1]}{[b_1]}\ldots\frac{d[b_{n-1}]}{[b_{n-1}]}=\F_r\theta([P]).\]
Furthermore, we recall $\theta(\Div(1-ax^r))=\V_r([a])$ (this follows from (\ref{3.4.3}) and \ref{0.10}) and since
\[ \sV_r([P])=[\Div(1-ax^r)]\times[(b_1,\ldots,b_{n-1})],\]
we obtain
\[\theta(\sV_r[P])=\theta(\Div(1-ax^r))\frac{d[b_1]}{[b_1]}\ldots\frac{d[b_{n-1}]}{[b_{n-1}]}
                  =\V_r([a])\frac{d[b_1]}{[b_1]}\ldots\frac{d[b_{n-1}]}{[b_{n-1}]}=\V_r(\theta([P])).\]
                  
It remains to show $(\ref{3.8.2})$ and since $\st$ is distributive by definition, this is equivalent to 
\[\theta([P]\st[Q])=\theta([P])\theta([Q]),\]
for $P\in X_{r-1}\setminus Y_{r-1}$  and  $Q\in X_{s-1}\setminus Y_{s-1}$ arbitrary points.

{\em First case: $P$ is $k$-rational.} We write $P=(1/a,b_1,\ldots, b_{r-1})$, $a,b_i\in k^\times$,  and denote by $L=k(Q)$ the residue field of $Q$.
Then $\mu_*([P]\times [Q])=[Q']$ is a point whose residue field $k(Q')=L'$ is over $k$ isomorphic to L and the isomorphism is given by
\eq{3.8.3}{L'\stackrel{\simeq}{\longrightarrow} L,\quad x(Q')\mapsto \frac{1}{a} x(Q),\quad y_i(Q')\mapsto y_{i-r-1}(Q), \text{ for } i=r,\ldots, r+s-2.}
Thus
\begin{eqnarray*}
\theta([P]\st[Q]) & = & \Tr_{L'/k}\Bigl(\frac{1}{[x(Q')]}\frac{d[b_1]}{[b_1]}\ldots\frac{d[b_{r-1}]}{[b_{r-1}]}
                                            \frac{d[y_r(Q')]}{[y_r(Q')]}\ldots\frac{d[y_{r+s-2}(Q')]}{[y_{r+s-2}(Q')]}\Bigr)\\
                  & = & \frac{d[b_1]}{[b_1]}\ldots\frac{d[b_{r-1}]}{[b_{r-1}]}
                               \Tr_{L'/k}\Bigl(\frac{1}{[x(Q')]}\frac{d[y_r(Q')]}{[y_r(Q')]}\ldots\frac{d[y_{r+s-2}(Q')]}{[y_{r+s-2}(Q')]}\Bigr),
\end{eqnarray*}
which, by the isomorphism (\ref{3.8.3}) equals
\[\frac{d[b_1]}{[b_1]}\ldots\frac{d[b_{r-1}]}{[b_{r-1}]}
            \Tr_{L/k}\Bigl(\frac{[a]}{[x(Q)]}\frac{d[y_1(Q)]}{[y_1(Q)]}\ldots\frac{d[y_{s-1}(Q)]}{[y_{s-1}(Q)]}\Bigr)
       =\theta([P])\theta([Q]).\]
       
{\em General case.} Now let $P$ and $Q$ be arbitrary. Let $L$ be the residue field of $P$ and $\pi:\Spec L\to \Spec k$ be the map induced by $k\subset L$. 
As we saw above, there is a point $P'\in X_{r-1,L}\setminus Y_{r-1,L}$ with $\pi_*[P']=[P] $. By the first case and lemma \ref{3.7}, (iii) we have 
\[\theta([P']\st\pi^*[Q])=\theta([P'])\theta(\pi^*[Q])=\theta([P'])\pi^*\theta([Q])\quad \text{ in } \Th{r+s-1}{L}{m}.\]
Now we apply $\Tr_{L/k}$ to this and use  lemma \ref{3.7} (ii), (iii) to obtain the assertion.
\end{proof}

Let us recall the cubical definition of {\em Bloch's higher Chow groups over a field on the level of zero cycles} (of course this can be done
more general, see \cite{To92} for the definition of the cubical version of Bloch's higher Chow groups and \cite{Bl86} for a simplicial version).
\begin{defn}[cf. \cite{To92}]\label{3.9}
Let $F$ be a field and denote 
\[c^n(F,n)=\Zy_0\bigl((\P^1_F\setminus\{0,1,\infty\})^n\bigr)\]
and
\[c^n(F,n+1)=\oplus \Z [C],\]
where the sum is over curves $C\subset (\P^1_F\setminus\{1\})^{n+1}$, satisfying $\partial^j_i[C]\in c^n(F,n)$, for $j=0,\infty$ and $i=1,\ldots,n$.
Then we define
\[\Ch^n(F,n)=\frac{c^n(F,n)}{\partial c^n(F,n+1)}.\]
The exterior product of cycles (see \ref{A7}) makes $\bigoplus_n \Ch^n(F,n)$ into a graded ring.
\end{defn}

\begin{lem}\label{3.10}
Let $P_0\in\G_{m,k}$ be a closed point. Denote by 
\[\iota_{P_0}: (\P^1_{k(P_0)}\setminus\{1\})^n= \Spec k(P_0)\times_k(\P^1_k\setminus\{1\})^n\inj\G_{m,k}\times_k(\P^1_k\setminus\{1\})^n\] 
the closed embedding induced by $\Spec k(P_0)\inj \G_m$. Then push-forward yields a homomorphism of groups
\[{\iota_{P_0}}_* :\Ch^{n}(k(P_0),n)\longrightarrow \Th{n+1}{k}{m},\]
for all $n\ge 0$ and $m\ge 1$, and 
\[\bigoplus_{P_0\in\G_m}\Ch^{n}(k(P_0),n))\stackrel{\oplus{\iota_{P_0}}_*}{\longrightarrow}\Th{n+1}{k}{m}\]
is surjective.
\end{lem}
\begin{proof}
To show that ${\iota_{P_0}}_*$ is well defined, let $[C]\in c^n(k(P_0),n+1)$ be a curve and define  $\xi={\iota_{P_0}}_*[C]\in\Zy_1(X_{n+1})$.
Clearly $\xi$ is in good position and since $P_0\neq 0$ the modulus condition is fulfilled trivially  for all $m$. Thus $\xi\in\Zy_1(X_{n+1};m)$ and
$\partial \xi= {\iota_{P_0}}_*\partial[C]$. Hence ${\iota_{P_0}}_*$ is well defined. For the surjectivity, 
take a closed point $P\in\G_m\times(\P^1\setminus\{0,1,\infty\})^n$ and write the maximal ideal of $P$ 
in the form (cf. \cite[proof of Lemma 2]{To92})
\[\mathfrak{m}=(f(x),g_1(x,y_1),\ldots,g_n(x,y_1,\ldots, y_n)),\quad f\in k[x] \text{ irreducible }, g_i\in k[x, y_1,\ldots,y_i].\]
Now let $P_0\in \G_m$ be the point defined by $f(x)=0$ and $P'\in(\P^1_{k(P_0)}\setminus\{0,1,\infty\})^n$ the point defined by the maximal ideal
 $\overline{\mathfrak{m}}=(g_1(\bar{x},y_1),\ldots,g_n(\bar{x},y_1,\ldots, y_n))\subset k(P_0)[y_1,\ldots,y_n]$. Then 
 \[{\iota_{P_0}}_*[P']=[P],\]
 hence the assertion.
\end{proof}

\begin{mkt}\label{3.11}
We recall the basic facts about Milnor $K$-theory, our reference is \cite{BaTa73}.

Let $F$ be a field. Denote by $T_*(F^\times)$ the tensor algebra (over $\Z$) of $F^\times$ and by $R_*(F^\times)\subset T_*(F^\times)$
the homogeneous ideal generated by $ a\otimes (1-a)$, for $a\in F^\times$. Then the {\em Milnor ring of $F$} is the graded algebra
\[K^M_*(F)=T_*(F^\times)/R_*(F^\times).\]
We denote by $K^M_n(F)$ the group of homogeneous elements of degree $n$ and by $\{a_1,\ldots, a_n\}$ the image of $a_1\otimes\ldots \otimes a_n$
in $K^M_n(F)$. We have
\begin{enumerate}
\item $\{a_1,\ldots,a_n\}\cdot\{b_1,\ldots,b_m\}=\{a_1,\ldots a_n,b_1,\ldots,b_m\}$.
\item $\{a_1a_1',a_2\ldots,a_n\}=\{a_1,a_2\ldots,a_n\}+\{a_1',a_2\ldots,a_n\}$.
\item $\{a,b\}=-\{b,a\}$.
\item $\{a,a\}=\{a,-1\}$, in particular $2\{a,a\}=0$.
\end{enumerate} 
Let $L\supset F$ be a finite field extension. Then  Bass-Tate and Kato (see \cite{BaTa73}, \cite{Ka80}) constructed a group homomorphism 
(the construction is due to Bass and Tate, and Kato showed that it is well defined)
\[\Nm_{L/F}:K^M_n(L)\longrightarrow K^M_n(F),\]
 satisfying the following properties
\begin{enumerate}
\item $\Nm_{L/F}(i(x)y)=x\cdot \Nm_{L/F}(y)$, where $x\in K^M_{n-r}(F)$, $y\in K^M_r(L)$ and $i:K^M_r(F)\to K^M_r(L)$ is the natural map.
\item  $\Nm_{L/F}(\{a\})=\{\Nm_{L/F}(a)\}$, for all $a\in F^\times$, where the $\Nm$ on the right hand side is just the usual norm of field extensions.
\item  For field extensions $L\supset K\supset F$ we have 
        \[\Nm_{L/F}=\Nm_{K/F}\circ\Nm_{L/K}\text{ and }\Nm_{L/L}=\id.\]
\end{enumerate} 
\end{mkt}

\begin{thm}[\cite{NeSu89}, \cite{To92}]\label{3.11.5}
The map
\[\rho: K^M_*(F)\stackrel{\simeq}{\longrightarrow} \bigoplus_n\Ch^{n}(F,n),\quad \{a_1,\ldots, a_n\}\mapsto (a_1,\ldots, a_n),\]
is an isomorphism of anti commutative graded rings with inverse
\[\rho^{-1} :\Ch^{n}(F,n)\longrightarrow K^M_n(F),\quad [P]\mapsto \Nm_{F(P)/F}(\{y_1(P),\ldots,y_n(P)\}).\]
\end{thm}
\begin{lem}[\cite{NeSu89}, Lemma 4.7.]\label{3.12}
Let $L\supset F$ be a finite field extension and $\pi:\Spec L\to\Spec F$ the morphism induced by the inclusion. Then
the following diagram commutes
\[\xymatrix{K^M_n(L)\ar[r]^{\rho_L}\ar[d]^{\Nm_{L/F}} & \Ch^{n}(L,n)\ar[d]^{\pi_*}\\
            K^M_n(F)\ar[r]^{\rho_F}                   & \Ch^{n}(F,n).
            }\] 

\end{lem}
\begin{defn}\label{3.12.5}
We denote
\[\Gamma_n=\bigoplus_{P_0\in\G_{m,k}}K^M_n\bigl(k(P_0)\bigr)\]
and write $\{a_1,\ldots, a_n\}_{P_0}$ for the image of $a_1\otimes\ldots\otimes a_n$ in $K^M_n\bigl(k(P_0)\bigr)$,  $P_0\in\G_m$. 
\begin{enumerate}
\item For $P_0,Q_0\in \G_m$ with $[P_0\times Q_0]=\sum_{i=1}^n m_i [R_i]\in\Zy_0(\G_m\times_k\G_m)$ and $r,s\ge 1$, we define
      \[\{a_1,\ldots,a_r\}_{P_0}\st\{b_1,\ldots,b_s\}_{Q_0}=
                                                     \sum_{i=1}^n m_i\,\Nm_{k(R_i)/k(\mu(R_i))}(\{a_1,\ldots,a_r,b_1,\ldots,b_s\}_{R_i}),\]
      where $\mu : \G_m\times\G_m\to \G_m$ is the multiplication and the $a_i$'s and $b_j$'s on the right hand side are viewed as elements in $k(R_i)$
      via the natural inclusions $k(P_0), k(Q_0)\inj k(R_i)$. Bilinear extension yields a map
      \[\st: \Gamma_r\otimes\Gamma_s\longrightarrow \Gamma_{r+s}.\]
      Note that this map is clearly well defined.
\item For $n\ge 1$, we define a group homomorphism 
      \[\sD:\Gamma_n\longrightarrow \Gamma_{n+1}\]
      via
      \[\sD(\{a_1,\ldots,a_n\}_{P_0})=-\{x(P_0),a_1,\ldots,a_n\}_{P_0},\]
      where $x$ is, as before, the coordinate of $\G_m$ and $x(P_0)$ is its residue class in $k(P_0)$.
\item Let $\varphi_r:\G_m\to\G_m$ be the morphism induced by $x\mapsto x^r$. Let $P_0\in\G_m$ be a point and write
       $\varphi_r^*[P_0]=\sum_{i=1}^s m_i [Q_i]\in \Zy_0(\G_m)$. Then we define
             \[\sV_r(\{a_1,\ldots,a_n\}_{P_0})=\sum_{i=1}^s m_i\{\varphi_r^*(a_1),\ldots,\varphi_r^*(a_n)\}_{Q_i},\]
        where $\varphi_r^*(a_j)$ is the image of the $a_j$ under the inclusion $k(P_0)\inj k(Q_i)$. This gives a well defined map
         \[\sV_r:\Gamma_n\longrightarrow \Gamma_n.\]
\item We define
        \[\sF_r:\Gamma_n\longrightarrow \Gamma_n\]
        by
        \[\sF_r(\{a_1,\ldots,a_n\}_{P_0})=\Nm_{k(P_0)/k(\varphi_r(P_0))}(\{a_1,\ldots,a_n\}_{P_0}).\]
\end{enumerate}
\end{defn}

\begin{lem}\label{3.13}
For $\alpha\in\Gamma_r$ and $\beta\in\Gamma_s$, we have
\begin{enumerate}
\item $\alpha\st\beta=(-1)^{rs}\beta\st\alpha$.
\item $\sD(\alpha\st\beta)=\sD(\alpha)\st\beta+(-1)^r\alpha\st\sD(\beta)$.
\item $2\sD\sD(\alpha)=0$.
\item $\sD\sF_r=r\sF_r\sD,\quad r\sD\sV_r=\sV_r\sD$.
\item $\sF_r\sD\sV_r(\{b_1,\ldots,b_s\}_{P_0})=\sD(\{b_1,\ldots,b_s\}_{P_0})-(r-1)\{-1,b_1,\ldots,b_s\}_{P_0}.$
\end{enumerate}
\end{lem}
\begin{proof}
We may assume $\alpha=\{a_1,\ldots,a_r\}_{P_0}\in K^M_r\bigl(k(P_0)\bigr)$ and $\beta=\{b_1,\ldots,b_s\}_{Q_0}\in K^M_s\bigl(k(Q_0)\bigr)$.
Then (i) and (iii) follow immediately from the properties in \ref{3.11}. For (ii) we write $[P_0\times Q_0]=\sum_{i=1}^n m_i[R_i]\in\Zy_0(\G_m\times\G_m)$,
then by definition
\begin{eqnarray*}
\sD(\alpha\st\beta) & = &- \sum_{i=1}^n m_i\,\{x\bigl(\mu(Ri)\bigr)\}\Nm_{k(R_i)/k(\mu(R_i))}(\{a_1,\ldots,a_r,b_1,\ldots,b_s\}_{R_i})\\
                    & = & -\sum_{i=1}^n m_i\,\Nm_{k(R_i)/k(\mu(R_i))}(\{x_1(P_0)x_2(Q_0),a_1,\ldots,a_r,b_1,\ldots,b_s\}_{R_i}) \\
                    & = & \sD(\alpha)\st\beta+(-1)^r \alpha\st\sD(\beta). 
\end{eqnarray*}
The first equation in (iv) follows from
\[\{x(\varphi_r(P_0))\}_{\varphi_r(P_0)}\Nm_{k(P_0)/k(\varphi_r(P_0))}(\{b_1,\ldots,b_s\}_{P_0})=
                                                                        \Nm_{k(P_0)/k(\varphi_r(P_0))}(\{x(P_0)^r,b_1,\ldots,b_s\}_{P_0}).\]
For $P_0\in\G_m$, write $\varphi_r^*[P_0]=\sum_{i=1}^n m_i[Q_i]\in\Zy_0(\G_m)$, then the second equation in (iv) follows from
\[\sum_{i=1}^n m_i r\{x(Q_i),\varphi_r^*(b_1),\ldots,\varphi_r^*(b_s)\}_{Q_i}=
                                                  \sum_{i=1}^n m_i \{\varphi_r^*(x(P_0)),\varphi_r^*(b_1),\ldots,\varphi_r^*(b_s)\}_{Q_i}.\]
In (v) we have
\ml{3.13.1}{\sF_r\sD\sV_r(\{b_1,\ldots,b_s\}_{P_0})=-\sum_{i=1}^n m_i\Nm_{k(Q_i)/k(P_0)}(\{x(Q_i),\varphi_r^*(b_1),\ldots\varphi_r^*(b_s)\}_{Q_i})\\
                                         =-\Bigl(\sum_{i=1}^n m_i \{\Nm_{k(Q_i)/k(P_0)}(x(Q_i))\}_{P_0}\Bigr)\{b_1,\ldots,b_s\}_{P_0}.}
But if we write in $k(P_0)[T]$
\[T^r-x(P_0)=\prod_i f_i(T)^{m_i} \quad \text{with } f_i\in k(P_0)[T] \text{ irreducible of degree } d_i,  \]
then $Q_i=\Spec \frac{k(P_0)[T]}{f_i(T)}$ and 
 \[\Nm_{k(Q_i)/k(P_0)}(x(Q_i))=(-1)^{d_i}f_i(0).\]
Now  (v) follows from (\ref{3.13.1}) and $\prod_i (-1)^{d_im_i}f_i(0)^{m_i}=(-1)^{(r-1)}x(P_0)$.
\end{proof}

\begin{lem}\label{3.14}
Let $\psi_n : \Gamma_n \longrightarrow \Th{n+1}{k}{m}$
be the composition 
\[\Gamma_n\stackrel{\ref{3.11.5}}{\longrightarrow}\bigoplus_{P_0\in\G_m}\Ch^{n}(k(P_0),n)\stackrel{\ref{3.10}}{\longrightarrow} \Th{n+1}{k}{m}.\] 
Then $\psi_n$ is surjective and the following diagrams commute
\begin{enumerate}
\item
\[\xymatrix{ \Gamma_r\otimes\Gamma_s\ar[rr]^(.3){\psi_r\otimes\psi_s}\ar[d]^\st &  & \Th{r+1}{k}{m}\otimes\Th{s+1}{k}{m}\ar[d]^\st\\
             \Gamma_{r+s}\ar[rr]^(.3){\psi_{r+s}}                               &  & \Th{r+s+1}{k}{m}.
            }\]
\item
\[\xymatrix{ \Gamma_n\ar[r]^(.3){\psi_n}\ar[d]^{\sD} & \Th{n+1}{k}{m}\ar[d]^{\sD}\\
             \Gamma_{n+1}\ar[r]^(.3){\psi_{n+1}}     &\Th{n+2}{k}{m}.
            }\]
\item
\[\xymatrix{ \Gamma_n\ar[r]^(.3){\psi_n}\ar[d]^{\sV_r} & \Th{n+1}{k}{m}\ar[d]^{\sV_r}\\
             \Gamma_n\ar[r]^(.2){\psi_n}     &\Th{n+1}{k}{rm+r-1}.
            }\]
\item
\[\xymatrix{ \Gamma_n\ar[r]^(.2){\psi_n}\ar[d]^{\sF_r} & \Th{n+1}{k}{rm+r-1}\ar[d]^{\sF_r}\\
             \Gamma_n\ar[r]^(.3){\psi_n}     &\Th{n+1}{k}{m}.
            }\]
\end{enumerate}
\end{lem}

\begin{proof}
The surjectivity of $\psi_n$ is clear. To prove the commutativity of the first diagram, take points $P_0,Q_0\in\G_m$, 
write $[P_0\times Q_0]=\sum_{i=1}^n m_i [R_i]\in\Zy_0(\G_m\times\G_m)$ and let $\alpha=\{a_1,\ldots,a_r\}_{P_0}$ and $\beta=\{b_1,\ldots,b_s\}_{Q_0}$ be 
elements in $\Gamma_r$ and $\Gamma_s$ respectively. Then
\[\psi_{r+s}(\alpha\st\beta)=\sum_{i=1}^n m_i\, {\iota_{\mu(R_i)}}_*\Nm_{k(R_i)/k(\mu(R_i))}(\{a_1,\ldots,a_r,b_1,\ldots,b_s\}_{R_i}).\]
But it follows from lemma \ref{3.12} that we have in $\Ch^{r+s}\bigl(k\bigl(\mu(R_i)\bigr),r+s\bigr)$
\[\Nm_{k(R_i)/k(\mu(R_i))}(\{a_1,\ldots,a_r,b_1,\ldots,b_s\}_{R_i})={\mu_{|R_i}}_*(a_1,\ldots,a_r,b_1,\ldots,b_s).\]
Hence
\begin{eqnarray*}\psi_{r+s}(\alpha\st\beta) & = & \sum_{i=1}^n m_i\,\mu_*{\iota_{R_i}}_*(a_1,\ldots,a_r,b_1,\ldots,b_s)
                            = \mu_*[{\iota_{P_0}}_*(a_1,\ldots,a_r)\times{\iota_{Q_0}}_*(b_1,\ldots,b_s)]\\
                                            & = & \psi_r(\alpha)\st\psi_s(\beta).
\end{eqnarray*}
The commutativity of the second diagram amounts to show
\[{\Delta_n}_*{\iota_{P_0}}_*(a_1,\ldots,a_n)={\iota_{P_0}}_*(x(P_0),a_1,\ldots,a_n).\]
For $P_0=1$ both sides are zero (the left hand side by definition of ${\Delta_n}_*$ and the right hand side since $\psi_{n+1}$ is well defined)
and for $P_0\neq 0$ the equality holds again by definition of $\Delta_n$.
 Finally diagram (iii) commutes just by the definition of the $\sV_r$'s and diagram (iv) by the definition of the $\sF_r$ together with lemma
\ref{3.12}. 
\end{proof}

Now we are ready to state and prove our main result.
\pagebreak[3]

\begin{thm}\label{3.16}
Let $k$ be a field of characteristic $\neq 2$. Then 
\[\bigl((\oplus_{n\ge 1}\Th{n}{k}{m})_{m\in\N},\, \R,\st,\sD, (\sF_r)_{r\ge 1},(\sV_r)_{r\ge 1})\] 
is a restricted Witt complex (see \ref{1.14.2}) and the natural map induced by the universality of the de Rham-Witt complex
\[\vartheta: \BWC{m}{n-1}{k}\stackrel{\simeq}{\longrightarrow} \Th{n}{k}{m}\]
is an isomorphism for all $m,n\ge 1$ with inverse the map $\theta$ from theorem \ref{3.2}.
\end{thm}
\begin{rmk}\label{3.17}
For $m=1$ the statement becomes
\[\Th{n}{k}{1}\cong\Omega^{n-1}_{k/\Z},\]
which was proved by Bloch and Esnault in \cite{BlEs03} under the additional assumption $1/6\in k$. 
\end{rmk}
\begin{proof}[Proof of the theorem]
First we have to show that $(\bigoplus_{n\ge 1}\Th{n}{k}{m})_{m\in\N}$ is a restricted Witt complex, that is, we have to check all the properties
from definition \ref{1.14.2}. Notice, that everything automatically commutes with the restriction $\R$. Next we observe that $\st$ is associative,
since the multiplication map on $\G_m$ is. We already know from corollary \ref{3.4} that
we have an isomorphism $\Th{1}{k}{m}\cong\BW{m}{k}$. Thus, since $2$ is invertible in $k$, $2$ is also invertible in $\bigoplus_{n\ge 1}\Th{\,n}{k}{m}$.
Hence, by the lemmas \ref{3.13} and \ref{3.14}, $(\bigoplus_{n\ge 1}\Th{n}{k}{m})_{m\in\N}$ is a projective system of differential graded $\Z$-algebras. 
Furthermore it follows directly from the definitions of $\st$, $\sF_r$ and $\sV_r$ that $\sF_r$ is a homomorphism of graded rings, for all $r$, and
we have $\sF_1=\sV_1=\id$, $\sF_r\sF_s=\sF_{rs}$ and $\sV_r\sV_s=\sV_{rs}$, as well as $\sF_r\sV_r=r$ (this is, pulling back and then pushing forward via 
a finite and flat map is just multiplication with the degree). Next we want to show $\sF_r\sV_s=\sV_s\sF_r$, if $(r,s)=1$. By proposition
\ref{A6} this amounts to show that
\[\xymatrix{ \G_m\ar[r]^{\varphi_s}\ar[d]^{\varphi_r} & \G_m\ar[d]^{\varphi_r}\\
             \G_m\ar[r]^{\varphi_s}                   & \G_m
            }\]
is cartesian, where $\varphi_r$ is induced by $x\mapsto x^r$. This means the following: view $k[x_1,1/x_1]$ as a $k[x,1/x]$-module via $x\mapsto x_1^s$
and $k[x_2,1/x_2]$ as a $k[x,1/x]$-module via $x\mapsto x_2^r$, then we have to show
\[A=k[x_1,\tfrac{1}{x_1}]\otimes_{k[x,\tfrac{1}{x}]}k[x_2,\tfrac{1}{x_2}]\longrightarrow k[t,\tfrac{1}{t}],\quad f(x_1)\otimes g(x_2)\mapsto f(t^r)g(t^s)\]
is an isomorphism, if $(r,s)=1$. 

Obviously the map is well defined and surjective. For the injectivity take an element $\alpha=\sum_{i,j\in\Z}\lambda_{i,j} x_1^i\otimes x_2^j\in A$. Then
$\alpha$ goes to zero in $k[t,1/t]$, iff for all $q\in\Z$
\[0=\sum_{ir+js=q}\lambda_{ij}=\sum_{a\in\Z}\lambda_{i_q-as,j_q+ar},\quad i_q,j_q\text{ fixed, with } i_qr+j_qs=q .  \] 
Thus
\[\alpha=\sum_q(\sum_a\lambda_{i_q-as,j_q+ar})x_1^{i_q}\otimes x_2^{j_q}=0,\]
hence the claim. It remains to check the equalities (iii)-(v) from definition \ref{1.14.2}.

Next we show (iii), i.e. $\sV_r(\sF_r([P])\st[Q])=[P]\st\sV_r([Q])$. But again by proposition \ref{A6} this is satisfied, if 
\[\xymatrix{ \G_m\times\G_m\ar[r]^{\id\times\varphi_r}\ar[d]^\mu & \G_m\times\G_m\ar[d]^{\mu\circ(\varphi_r\times\id)}\\
                \G_m\ar[r]^{\varphi_r}                           & \G_m
           }\]
is  cartesian. That is: view $k[x_1,1/x_1,x_2,1/x_2]$ as a $k[z,1/z]$-module via $z\mapsto x_1^rx_2$ and 
$k[x,1/x]$ as a $k[z,1/z]$-module via $z\mapsto x^r$, then 
\[k[x_1,\tfrac{1}{x_1},x_2,\tfrac{1}{x_2}]\otimes_{k[z,\tfrac{1}{z}]}k[x,\tfrac{1}{x}]\longrightarrow k[x_1,\tfrac{1}{x_1},x_2,\tfrac{1}{x_2}],\quad
                                        f(x_1,x_2)\otimes g(x)\mapsto f(x_1,x_2^r)g(x_1x_2)\]
has to be an isomorphism. Indeed the map is clearly well defined, the surjectivity follows from $x_1\otimes 1\mapsto x_1$ and 
$(1/x_1)\otimes x\mapsto x_2$ and the injectivity is proved as above. Hence (iii).

The equalities in definition \ref{1.14.2}, (iv), i.e.
\[\sF_r\sD\sV_r=\sD,\quad \sD\sF_r=r\sF_r\sD,\quad \sV_r\sD=r\sD\sV_r,\]
follow immediately from the lemmas \ref{3.13} and \ref{3.14} (taking into account that $\psi_{s+1}(\{-1,b_1,\ldots,b_s\}_{P_0})=0$, 
since $2$ is invertible in $\bigoplus_{n\ge 1}\Th{n}{k}{m}$).

It remains to prove relation (v), i.e. $\sF_r\sD(\{a\})=\{a\}^{r-1}\sD(\{a\})$ for all $a\in k$. If $a=0$, then by definition $\{a\}=0\in\Th{1}{k}{m}$.
If $a=1$ both sides are zero by definition of $\sD$ and for $a\in k^\times\setminus\{1\}$, both sides equal $(1/a^r,1/a)$.

Thus $\bigoplus_{n\ge 1}\Th{n}{k}{m}$ is a restricted Witt complex and the universality of the de Rham-Witt complex yields a map
\[\vartheta: \BWC{m}{n-1}{k}\longrightarrow \Th{n}{k}{m}.\]
For $\theta$ the map from theorem \ref{3.2}, we have by corollary \ref{3.8}, $\theta\circ\vartheta=\id $. Thus it remains to check that
$\vartheta$ is surjective. Similar to the arguments in  \cite[Lemma 4.7.]{NeSu89} and \cite[Theorem 6.4.]{BlEs03}, this will be done with the use of the 
following lemma.
\begin{lem}\label{3.18}
Let $L\supset k$ be a finite field extension, inducing $\pi :\Spec L\to \Spec k$. Then the following diagram commutes
\[\xymatrix{\BWC{m}{n-1}{L}\ar[r]^\vartheta\ar[d]^{\Tr_{L/k}} & \Th{n}{L}{m}\ar[d]^{\pi_*}\\
            \BWC{m}{n-1}{k}\ar[r]^\vartheta                   & \Th{n}{k}{m}.
            }\]
\end{lem}
Now assuming the lemma holds, we can prove the surjectivity of $\vartheta$. Take a point $P\in X_{n-1,k}\setminus Y_{n-1,k}$ and denote by $L=k(P)$ 
its residue field. Then there is a $L$-rational point $P'=(a,b_1,\ldots,b_{n-1})\in X_{n-1,L}\setminus Y_{n-1,L}$ with $\pi_*[P']=[P]$.
On the other hand we may write (using the notation from definition \ref{3.6}, (i), $\{b\}=\Div(1-bx)$)
\[[P']=(-1)^{n-1}\{\tfrac{b_1\ldots b_{n-1}}{a}\}\st\sD(\{\tfrac{1}{b_1}\})\st\ldots \st\sD(\{\tfrac{1}{b_{n-1}}\})\]
Hence by lemma \ref{3.18}
\[[P]=\vartheta\Bigl(\Tr_{L/k}\bigl((-1)^{n-1}\tfrac{[b_1\ldots b_{n-1}]}{[a]}d(\tfrac{1}{[b_1]})\ldots d(\tfrac{1}{[b_{n-1}]})\bigr)\Bigr)
                                                            \text{ in } \Th{n}{k}{m}.\]
What remains to be done is the \\
\\
{\em Proof of Lemma \ref{3.18}.}
First we prove the assertion on the level of Witt vectors, i.e. $n=1$. Since $\Tr$ and $\pi_*$ commute with $\V_r$ and $\sV_r$ respectively, we are 
reduced to show $\vartheta(\Tr_{L/k}[a])=\pi_*\{a\}$. But by the construction of the trace  on the Witt ring (via norm) and by corollary \ref{3.4}, 
we have
\[\vartheta(\Tr_{L/k}([a]))=\Div(\Nm_{L/k}(1-aT)) \text{ in }\Th{1}{k}{m},\]
which equals $\pi_*\{a\}$, by \cite[Proposition 1.4.]{Fu84}. Thus the lemma is true for $n=1$.
By construction of the trace (theorem \ref{2.7}), it is enough to consider in general the cases, $L\supset k$ separable or purely inseparable of 
degree $p$. 

{\em First case: $L\supset k$ separable.} We have $\BWC{m}{n-1}{L}=\BW{m}{L}\otimes_{\BW{m}{k}}\BWC{m}{n-1}{k}$. Thus it is enough to consider
elements $w\alpha$ with $w\in\BW{m}{L}$ and $\alpha\in\BWC{m}{n-1}{k}$. Hence the assertion follows from
\[\pi_*(\vartheta(w\alpha))=\pi_*(\vartheta(w)\st\pi^*\vartheta(\alpha))=\vartheta(\Tr(w)\alpha),\]
where the last equality is  the case $n=1$ together with lemma \ref{3.7} (ii).

{\em Second case: $L\supset k$ purely inseparable of degree $p$.} We can write $L=k[x]/(x^p-b)$, with $b\in k\setminus k^p$. Denoting by $y$ the image
of $x$ in $L$ it follows from theorem \ref{2.1} and theorem \ref{1.13} that every element in $\BWC{m}{n-1}{L}$ can be written as a sum of elements
of the following form
\begin{enumerate}
\item $\V_r(\alpha[y]^j),\quad \alpha\in\BWC{\lfloor\frac{m}{r}\rfloor}{n-1}{k}, j\in\{0,\ldots,p-1\},\,r\in\{1,\ldots,m\}$,
\item $\V_r(\beta [y]^{j-1}d[y]),\quad \beta\in\BWC{\lfloor\frac{m}{r}\rfloor}{n-2}{k}, j\in\{1,\ldots,p\},\, r\in\{1,\ldots,m\}$,
\item $d\V_r(\beta [y]^j),\quad \beta\in\BWC{\lfloor\frac{m}{r}\rfloor}{n-2}{k}, j\in\{0,\ldots,p-1\},\, r\in\{1,\ldots,m\}$.
\end{enumerate}
Thus it is enough to check the commutativity in the cases (i)-(iii). These are straightforward calculations, similar to the one in the first case.
To prove (ii) one has to use $[y]^{j-1}d[y]=\frac{1}{j}d[y]^j$, for $j<p$, and $[y]^{p-1}d[y]=\F_p(d[y])$. This finishes the proof of the
lemma and hence the proof of theorem \ref{3.16}.
\end{proof}

\begin{rmk}
Notice, that the relations (i)-(iii) and (v) from definition \ref{1.14.2} are already satisfied on the level of cycles and that we explicitly needed
the modulus condition only once, namely in the proof of corollary \ref{3.4}, where we showed the injectivity of the map 
$\theta: \Th{1}{k}{m}\to \BW{m}{k}$.
\end{rmk}
\newpage

%
%
%
%
%
%

\begin{appendix}
\section{Intersection Theory for Cartier Divisors}
In this appendix we give the definition of algebraic $n$-cycles,  flat pull-back, push-forward, pull-back of Cartier divisors and intersection
of those with cycles and we will list some properties. Of course this is all taken from \cite[Chapter 1, 2]{Fu84} (see also \cite[Chapter 2]{VoSuFr00}), 
except that in \cite{Fu84}
everything is done modulo rational equivalence and the push-forward is only defined for proper maps. Whereas here we are defining everything
on the level of cycles and the push-forward is defined whenever the map is proper along the cycle we wish to push-forward 
(see \cite[Chapter 2]{VoSuFr00}).\\
\\ 
All schemes are assumed to be algebraic, i.e. of finite type over a field.\\
\\
Let $k$ be a field and $X$ a $k$-scheme. We define for $n\in\N_0$
\[\Zy_n(X)=\bigoplus \Z \,[V],\]
where the sum is over all $n$-dimensional subvarieties $V\subset X$. If $X$ is equidimensional of pure dimension $d$ we also write
$\Zy^i(X):=\Zy_{d-i}(X)$. An element in $\Zy_n(X)$ is called a $n$-cycle and an element in $\Zy(X)=\bigoplus_{n\ge 0}\Zy_{n}(X)$ just a
cycle. 

Let $X_1,\ldots,X_r$ be the irreducible components of $X$ equipped with the reduced scheme structure
and $\eta_i,$ $i=1,\ldots,r$, their generic points, then the {\em cycle of $X$} is defined to be
\eq{A1}{[X]=\sum_{i=1}^r l_{\sO_{X_i,\eta_i}}(\sO_{X_i,\eta_i})[X_i]\in\Zy(X).} 
If $X$ is a variety with function field $K=k(X)$ and $V\subset X$ is a codimension one subvariety with generic point $\eta$,
then we define the {\em order of vanishing of $f\in K^\times$ along $V$} to be
\eq{A2}{\ord_V(f)=l_{\sO_{X,\eta}}(\sO_{X,\eta}/a)-l_{\sO_{X,\eta}}(\sO_{X,\eta}/b),\quad \text{ with } f=\frac{a}{b},\, a,b\in\sO_{X,\eta}. }

\begin{ex}[ \cite{Fu84}, Example 1.2.3. ]\label{A3}
Let $\tilde{X}\to X$ be the normalization of $X$ in $K$, then
\[\ord_V(f)=\sum_{\tilde{V}\surj V}[k(\tilde{V}):k(V)]\ord_{\tilde{V}}(f),\]
where the sum is over all subvarieties $\tilde{V}\subset \tilde{X}$ which map onto $V$. 
\end{ex}
\begin{pbcd}
Let $X'$ be a variety over $k$, $f: X'\to X$ a morphism of $k$-schemes and $D=\{U_i,g_i\}$ a Cartier divisor on $X$ satisfying $f(X')\not\subset |D|$.
Then $f^*g_i\in k(X')$ is defined and $f^*D=\{f^{-1}(U_i),f^*g_i\}$ is a well defined Cartier divisor on $X$.
We have $|f^*D|=f^{-1}(|D|)$. If $X''$ is a variety and $g: X''\to X'$ is a morphism of $k$ schemes with $g(X'')\not\subset f^{-1}(|D|)$, then
$g^*f^*D=(f\circ g)^*D$. 
\end{pbcd}

\begin{pf}\label{A4}
Let $f: X\to X'$ be a morphism of $k$-schemes. We say $f$ is {\em proper along a cycle} $\alpha\in\Zy(X)$, 
if the restriction of $f$ to the support of $\alpha$ is proper. If $f$ is proper along a $n$-dimensional variety $V\subset X$ we define
\[f_* [V]=\begin{cases}
            [k(V):k(f(V))] [f(V)] & \text{if } [k(V):k(f(V))]<\infty\\
                  0               & \text{else }
           \end{cases} \in \Zy_n(X').\]
We extend $f_*$ additively to all cycles $\alpha$ with $f$  proper along $\alpha$. If $g: X'\to X''$ is proper along $f_*\alpha$, $\alpha\in\Zy_n(X)$, then
$g_*f_*\alpha=(g\circ f)_*\alpha\in\Zy_n(X'')$. For a proper map $f: X\to X'$ we obtain a map $f_*:\Zy_n(X)\to\Zy_n(X')$.
If $f$ is a closed immersion, then $f_*$ is injective, hence we may view cycles  on $X$ as cycles on $X'$.
\end{pf}

\begin{fpb}[\cite{Fu84}, 1.7]\label{A5}
Let $f: X'\to X$ be a flat morphism of $k$-schemes of relative dimension $r$. Then we have a group homomorphism $f^*:\Zy_n(X)\to \Zy_{n+r}(X')$ uniquely
determined by 
\[f^*[V]=[f^{-1}(V)],\quad V\subset X \text{ a variety}.\]
If $Y\subset X$ is any closed subscheme, then $f^*[Y]=[f^{-1}(Y)]$. If $g:X''\to X'$ is flat of relative dimension $s$, then
$g^*f^*=(f\circ g)^*:\Zy_n(X)\to\Zy_{n+s+r}(X)$.
\end{fpb}

\begin{prop}\label{A6}
Let
\[\xymatrix{ X'\ar[r]^{g'}\ar[d]^{f'}   &  X\ar[d]^f\\
             Y'\ar[r]_g                 & Y
            }\]
be a fiber square with $g$ flat of relative dimension $r$. Let $\alpha$ be a cycle on $X$, with $f$ proper along $\alpha$. 
Then $g'$ is flat of relative dimension $r$, $f'$ is proper along $g'^*\alpha$ and
\[f'_*g'^*\alpha=g^*f_*\alpha.\]
\end{prop}
\begin{proof} 
Is the same as in \cite[Proposition 1.7.]{Fu84}. 
\end{proof}

\begin{lem}[\cite{Fu84}, Example 1.7.4.]\label{A6.5}
Let $f: X'\to X$ be a finite and flat $k$-morphism of degree $n$ and $\alpha$ a cycle on $X$. Then
\[f_*f^*\alpha=n\alpha.\] 
\end{lem}

\begin{ep}\label{A7}
Let $X,Y$ be two schemes of finite type over a field $k$. Then the exterior product
\[\Zy_n(X)\otimes_\Z\Zy_r(Y)\stackrel{\times}{\longrightarrow} \Zy_{n+r}(X\times_k Y),\quad \alpha\otimes\beta\mapsto \alpha\times\beta\]
is the bilinear form uniquely determined by
\[[V]\times[W]=[V\times_k W],\quad \text{for subvarieties } V\subset X, W\subset Y.\]
\end{ep}

\begin{prop}[\cite{Fu84}, 1.10]\label{A8}
Let $f: X'\to X$ and $g: Y'\to Y$ be two morphisms of $k$-schemes and $f\times g: X'\times_k Y'\to X\times_k Y$ the induced morphism.
\begin{enumerate}
\item Take $\alpha\in\Zy(X')$ and $\beta\in \Zy(Y')$ with $f$ proper along $\alpha$, $g$ proper along $\beta$. Then $f\times g$ is proper along
       $\alpha\times\beta$ and
       \[(f\times g)_*(\alpha\times\beta)=f_*\alpha\times g_*\beta.\]
\item  If $f$ and $g$ are flat of relative dimension $m$ and $n$, then $f\times g$ is flat of relative dimension $m+n$ and 
        \[(f\times g)^*(\alpha\times\beta)=f^*\alpha\times g^*\beta.\]
\item The exterior product is associative. 
\end{enumerate}
\end{prop}

\begin{lem}\label{A9}
Let $L\supset k$ be a finite field extension, $\pi:\Spec L\to \Spec k$ the induced map, $X$ a $k$-scheme, $\alpha$ a cycle on $X$ and 
$\beta$ a cycle on $X_L=X\times_k \Spec L$. Then
\[\pi_*(\pi^*\alpha\times_L \beta)=\alpha\times_k\pi_*\beta.\]
(Notice the sloppy notation, e.g. the $\pi_*$ on the left hand side is the push-forward from $X_L\times_L X_L=X\times_k X\times_k\Spec L$ to 
$X\times_k X$, etc. .) 
\end{lem}
\begin{proof}
We may assume $\alpha=[V]$ and $\beta=[W]$. Then
\[\pi_*(\pi^*[V]\times_L [W])=\pi_*[V\times_k\Spec L \times_L W]=\pi_*([V]\times_k [W])=[V]\times_k \pi_*[W],\]
where the last equality holds by \ref{A8}, (i).
\end{proof}

\begin{icd}[\cite{Fu84}, Remark 2.3.]\label{A10}
Let $X$ be a $k$-scheme and $D=\{U_i,g_i\}$ a Cartier divisor on $X$ such that 
\[(\ast)\quad {\sO_X(D)_|}_{|D|} \text{ is trivial }\]
(e.g. there is an open set $U\subset X$ with $|D|\subset U$ and $D_{|U}$ is a principal divisor). Let $j:V\inj X$ be a $n$-dimensional subvariety. 
Then we define
\[D.[V]=\begin{cases}
         [j^*D] & \text{if } V\not\subset |D|\\
           0    & \text{else }
         \end{cases}\in\Zy_{n-1}(|D|\cap V)\]
where $[j^*D]$ is the Weil divisor on $V$ associated to the Cartier divisor $j^*D$, i.e.
\[[j^*D]=\sum_{Z\in\Zy^1(V)}\ord_Z(j^*D)[Z],\]
with $\ord_Z(j^*D)=\ord_Z(j^*g_i)$ for any $i$ with $Z\cap j^{-1}(U_i)\neq\emptyset$. By linearity we obtain a map
\[\Zy_n(X)\longrightarrow\Zy_{n-1}(|D|),\quad \alpha\mapsto D.\alpha.\]
(Notice, that the assumption $(\ast)$ can be dropped if one defines the
intersection cycle modulo rational equivalence. Then $D.[V]$ does not need to be zero for $V\subset |D|$.)
\end{icd}

\begin{prop}\label{A11}
Let $X$ be a $k$-scheme and  $D$ a Cartier divisor satisfying $(\ast)$.
\begin{enumerate}
\item For $\alpha,\alpha'\in\Zy_n(X)$
      \[D.(\alpha+\alpha')=D.\alpha+D.\alpha'\quad \text{in } \Zy_{n-1}(|D|\cap(|\alpha|\cup|\alpha'|)).\]
\item If $D'$ is a Cartier divisor satisfying $(\ast)$ and  $\alpha\in \Zy_n(X)$, then
      \[(D+D').\alpha=D.\alpha+D'.\alpha\quad \text{in }\Zy_{n-1}((|D|\cup|D'|)\cap |\alpha|).\]
\item Let $X'$ be a variety over $k$, $\alpha\in\Zy_n(X')$, $f: X'\to X$ a $k$-morphism which is proper along $\alpha$ and satisfies 
      $f(X')\not\subset |D|$. Then $f^*D$ satisfies $(\ast)$, $g=f_{|f^{-1}(|D|)\cap |\alpha|}$ is proper along $f^*D.\alpha$ and
      \[g_*(f^*D.\alpha)=D.f_*\alpha\quad \text{in }\Zy_{n-1}(|D|\cap|\alpha|).\]
\item Let $f: X'\to X$ be a flat morphism of relative dimension $r$ with $f(X')\not\subset |D|$ and $\alpha\in \Zy_n(X)$.
      Then $f^*D$ satisfies $(\ast)$, $g=f_{|f^{-1}(|D|\cap |\alpha|)}$ is flat of relative dimension $r$ and
      \[ f^*D.f^*\alpha=g^*(D.\alpha)\quad \text{in }\Zy_{n+r-1}(f^{-1}(|D|\cap|\alpha|)).\]
      In particular $[f^*D]=f^*[D]$.
\item Let $Y$ be a $k$-scheme, $\alpha\in\Zy_n(X)$, $\beta\in \Zy_r(Y)$ and $p:X\times_k Y\to X$ the projection. Then
       \[(p^*D).(\alpha\times\beta)=(D.\alpha\times\beta)\quad\text{in }\Zy_{n+r-1}((|D|\cap |\alpha|)\times |\beta|).\]
\end{enumerate} 

\end{prop}
\begin{proof}
Although the statement differs slightly from \cite[Proposition 2.3., Example 2.3.1.]{Fu84}, the proof is almost the same and therefore skipped.
\end{proof}

\end{appendix}

\end{document}